\documentclass[a4paper,11pt,reqno]{amsart}
\usepackage[utf8]{inputenc}
\usepackage[margin=0.92in]{geometry}
\usepackage{array,longtable,booktabs}
\usepackage{amsmath,amssymb,amsthm,mathtools}
\usepackage{microtype}
\usepackage{enumitem}
\usepackage{graphicx}
\usepackage[dvipsnames]{xcolor}
\usepackage{tikz}
\usepackage{tikz-cd}
\usetikzlibrary{arrows.meta,positioning,calc}
\usepackage[hidelinks]{hyperref}

\allowdisplaybreaks
\setlist[enumerate]{topsep=3pt,itemsep=2pt,parsep=0pt}
\setlist[itemize]{topsep=3pt,itemsep=2pt,parsep=0pt}

\theoremstyle{plain}
\newtheorem{theorem}{Theorem}[section]
\newtheorem{proposition}[theorem]{Proposition}
\newtheorem{lemma}[theorem]{Lemma}
\newtheorem{corollary}[theorem]{Corollary}
\theoremstyle{definition}
\newtheorem{definition}[theorem]{Definition}

\newcommand{\C}{\mathbb C}
\newcommand{\Z}{\mathbb Z}
\newcommand{\R}{\mathbb R}
\newcommand{\Q}{\mathbb Q}
\newcommand{\Pj}{\mathbb P}
\newcommand{\cC}{\mathcal C}
\newcommand{\cE}{\mathcal E}
\newcommand{\cD}{\mathcal D}
\newcommand{\Db}{\mathrm D^{\mathrm b}}
\newcommand{\Dsg}{\mathrm D_{\mathrm{sg}}}
\DeclareMathOperator{\MF}{MF}
\DeclareMathOperator{\Ku}{Ku}
\DeclareMathOperator{\Perf}{Perf}
\DeclareMathOperator{\Ind}{Ind}
\DeclareMathOperator{\Hom}{Hom}
\DeclareMathOperator{\Ext}{Ext}

\DeclareMathOperator{\im}{im}

\DeclareMathOperator{\Cone}{Cone}

\newcommand{\gr}{\mathrm{gr}}
\newcommand{\num}{\mathrm{num}}
\newcommand{\id}{\mathrm{id}}
\newcommand{\OO}{\mathcal O}

\title{Gepner stability conditions on hypersurface categories}
\author[A. Aravena]{Anibal Aravena}
\author[P. Liu]{Peize Liu}
\author[J. Tevelev]{Jenia Tevelev}
\date{1 November 2026}
\hypersetup{pdftitle={Bridgeland stability conditions on hypersurface categories},
 pdfauthor={Anibal Aravena, Peize Liu, and Jenia Tevelev},pdfsubject={Bridgeland stability and hypersurface categories}}
\begin{document}

\begin{abstract}
We construct numerical Bridgeland stability conditions on the category of
graded matrix factorizations of every smooth complex projective
hypersurface of degree $d\ge2$.  These conditions are of Gepner type
$(\tau,2/d)$, where $\tau$ is the grading shift autoequivalence.
By Orlov's theorem, this gives
 stability conditions of Gepner type on all Calabi--Yau hypersurfaces
and on Kuznetsov components of all Fano hypersurfaces.
These conditions satisfy a mass--Hom bound and the support property
on the full numerical Grothendieck group.  Semistable objects of
fixed numerical class and phase admit proper good moduli spaces.
\end{abstract}

\maketitle

\section{Introduction}

Let $X\subset\Pj^{N-1}$ be a smooth complex hypersurface of degree
$d\ge2$.  When $N>d$, the derived category $\Db(X)$ contains an
exceptional collection of line bundles
$\OO_X,\OO_X(1),\ldots,\OO_X(N-d-1)$.
The remaining part is called {\em the Kuznetsov component} of the Fano hypersurface:
\begin{equation}\label{eq:intro-sod}
 \Db(X)=\langle\Ku(X),\OO_X,\OO_X(1),\ldots,
                                     \OO_X(N-d-1)\rangle.
\end{equation}
When $N=d$, the hypersurface is Calabi--Yau, and the category of
interest is $\Db(X)$ itself.  Both cases admit a common description
in terms of the homogeneous equation $W$ defining $X$: Orlov's
theorem~\cite{Orlov} identifies them with the category $\MF^{\gr}(W)$
of graded matrix factorizations.

We write $[n]$ for cohomological shift and $(a)$ for the internal
grading shift of a graded module, with $M(a)_j=M_{j+a}$.
We denote the grading-shift autoequivalence 
of the category of graded matrix factorizations
by $\tau$, so that $\tau^aE=E(a)$ and $\tau^d\simeq[2]$.

With a suitable normalization of Orlov's equivalence, the same
autoequivalence on $\Ku(X)$ when $N>d$, or on $\Db(X)$ when
$N=d$, is given by
\[
 \tau(E)\simeq\Cone\left(
   \operatorname{RHom}(\OO_X,E(1))\otimes_{\C}\OO_X
       \xrightarrow{\mathrm{ev}}E(1)\right),
 \qquad E(1)=E\otimes\OO_X(1),
\]
where the arrow is the evaluation map
\cite[Definition~6.11 and Theorem~6.13]{BFKHodge}.
In the Fano case, this is the \emph{rotation functor}: tensoring
by $\OO_X(1)$ followed by left mutation through $\OO_X$.
In the Calabi--Yau case, it is tensoring by $\OO_X(1)$ followed
by the spherical twist about $\OO_X$. 

It is natural to look for a Bridgeland stability
condition in which this fractional shift acts by a constant rotation
of phases.  Since $[1]$ rotates the phase by $1$,
the desired phase shift is $2/d$.

\begin{definition}\label{def:phase-shift-gepner}
Let $\mathcal D$ be a triangulated category with a stability
condition $\sigma=(Z,\mathcal P)$ and 
an exact
autoequivalence
$\Phi:\mathcal D\xrightarrow{\sim}\mathcal D$, and let $s\in\R$.
Following~\cite[Definition~1.1]{TodaMF}, we say that the stability condition $\sigma$ is
\emph{of Gepner type $(\Phi,s)$} if, for every $\phi\in\R$ and every $E\in\mathcal D$,
\[
 \Phi\mathcal P(\phi)=\mathcal P(\phi+s)\qquad\hbox{\rm and}\qquad
 Z(\Phi E)=e^{i\pi s}Z(E).
\]
In particular, $\Phi$ preserves semistability and translates phases by $s$.
\end{definition}

\begin{theorem}\label{thm:main}
Let $d\ge2$. Let
$X=\{W=0\}\subset\Pj^{N-1}$, $N\ge d$,
be a smooth hypersurface of degree~$d$.  The category
$\MF^{\gr}(W)$ of graded matrix factorizations admits a numerical Bridgeland stability
condition $\sigma=(Z,\mathcal P)$ satisfying the support property
on the full numerical Grothendieck group $K_{\num}(\MF^{\gr}(W))$.
It is of Gepner type $(\tau,2/d)$, where $\tau$ is the internal
grading shift; explicitly,
\begin{equation}\label{eq:main-gepner}
 Z(\tau E)=e^{2\pi i/d}Z(E),\qquad
 \tau\mathcal P(\phi)=\mathcal P\left(\phi+\frac2d\right).
\end{equation}
For every fixed object $A$, there is a constant $C_A>0$ such that
\begin{equation}\label{eq:main-mass-hom}
 \sum_{n\in\Z}\dim_{\C}\Hom(A,E[n])
       \le C_A m_\sigma(E)\qquad\text{for every }E.
\end{equation}
Here $m_\sigma(E)$ is the sum of the absolute values of the central
charges of its Harder--Narasimhan factors.
By~Orlov's theorem~\cite{Orlov}, $\MF^{\gr}(W)$ is equivalent to
$\Ku(X)$ for $N>d$ and to $\Db(X)$ for $N=d$.
\end{theorem}

Theorem~\ref{thm:all-degrees} also covers the case $N<d$,
when 
$\MF^{\gr}(W)$ contains $\Db(X)$.

The bound~\eqref{eq:main-mass-hom} applies existence of proper moduli spaces.

\begin{corollary}\label{cor:proper-gepner-moduli}
Let $\sigma=(Z,\mathcal P)$ be a stability condition constructed in
Theorem~\ref{thm:main} on $\mathcal C=\Ku(X)$ for $N>d$, or
on $\mathcal C=\Db(X)$ for $N=d$.
Fix $v\in K_{\num}(\mathcal C)$ and $\phi\in\R$ with
$Z(v)\in\R_{>0}e^{i\pi\phi}$.
The stack of semistable objects of class $v$ and phase $\phi$
is an algebraic stack of finite type over $\C$, has affine diagonal,
and admits a proper good moduli algebraic space $M_\sigma(v,\phi)$.
\end{corollary}

The closed points of this space classify S-equivalence classes
of semistable objects.  The corollary follows from the
mass--Hom bound by the theorem of Halpern-Leistner and
Robotis~\cite[Theorem~2.35]{HLR}; we give the precise application
at the end of Section~\ref{sec:hyperplane}.

For Fano hypersurfaces, even the existence of any stability condition
on the Kuznetsov component is a already a new result. Building on
his joint work with Macr\`i, Perry, Stellari, and Zhao~\cite{LMPSZ}
on rigidity of stability conditions for products of curves,
Li proves that $\Db(Y)$ admits stability conditions for every
smooth projective complex variety $Y$~\cite{LiExistence}.
This theorem does not settle existence of stability conditions on admissible subcategories.
Indeed, there can be no  existence theorem
in this generality, since nonzero phantom subcategories
have vanishing Grothendieck group and hence admit no stability
conditions~\cite[Definition~1.8 and Corollary~1.15]{GorchinskiyOrlov}.
For cubic hypersurfaces, several constructions were already
available. Bayer--Lahoz--Macr\`i--Stellari constructed numerical
stability conditions on the Kuznetsov component of every smooth
cubic fourfold~\cite{BLMS}, without asserting Gepner symmetry.
The stronger requirement of Gepner symmetry was established
by Toda for cubic surfaces, cubic threefolds, and a general
cubic fourfold containing a plane
\cite[Theorems~1.2 and~1.3]{TodaOK}.
More recently, Liu constructed stability conditions with the
Gepner property on the Kuznetsov component of every smooth
cubic fivefold~\cite[Corollary~5.12]{LiuCubic}.

For Calabi--Yau hypersurfaces, the category in question is the
full derived category, so Li's theorem supplies stability
conditions in every dimension. The additional challenge is
to construct one with Gepner symmetry. Toda did this for plane
cubics and quartic K3 surfaces~\cite[Theorem~1.3]{TodaMF}.
Li proved existence of  numerical Bridgeland stability conditions on smooth
quintic threefolds
by proving a suitable Bogomolov--Gieseker inequality
\cite{LiQuintic}. The search for Gepner symmetry led Toda
to formulate a stronger inequality~\cite{TodaQuintic},
subsequently proved by Xu for slope-semistable
sheaves~\cite{XuQuintic}. The recent OpenAI release includes
a manuscript asserting the existence of Gepner-type stability
conditions on every smooth quintic threefold~\cite{OAIQuintic}.

Our theorem gives a uniform construction of stability conditions
with Gepner symmetry on the Kuznetsov components of all smooth Fano
hypersurfaces and on the derived categories of all smooth Calabi--Yau
hypersurfaces. 

The proof of Theorem~\ref{thm:main}
has three steps.  We start with
the category of graded matrix factorizations of $x^d$, which is equivalent
to the derived category of representations of the  quiver
$A_{d-1}$.  The grading shift corresponds to $S^{-1}[1]$, where
$S$ is the Serre functor~\cite{Takahashi,KST}.
Qiu~\cite[Theorem~4.7]{QiuGepner} proves that the stability condition
of Gepner type $(S^{-1}[1],2/d)$ on this category is unique up to
the natural $\C$-action, as part of his treatment of Dynkin quivers.
These quiver categories
provide the building blocks of our construction.
Let $Q_d$ be the quiver $A_{d-1}$ with the following orientation: 
every edge is directed from
an odd-numbered vertex to an even-numbered vertex:
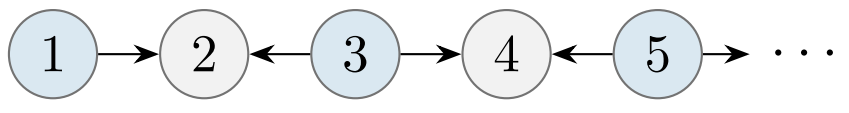
\begin{figure}[htbp]
 \begin{tikzpicture}[baseline=(v1.base),>=Stealth,
  vertex/.style={circle,draw=black!55,inner sep=3pt,font=\small}]
 \node[vertex,fill=MidnightBlue!10] (v1) at (0,0) {$1$};
 \node[vertex,fill=black!5] (v2) at (1,0) {$2$};
 \node[vertex,fill=MidnightBlue!10] (v3) at (2,0) {$3$};
 \node[vertex,fill=black!5] (v4) at (3,0) {$4$};
 \node[vertex,fill=MidnightBlue!10] (v5) at (4,0) {$5$};
 \node (dots) at (5,0) {$\cdots$};
 \draw[->] (v1)--(v2);
 \draw[->] (v3)--(v2);
 \draw[->] (v3)--(v4);
 \draw[->] (v5)--(v4);
 \draw[->] (v5)--(dots);
 \end{tikzpicture}
\caption{Alternating quiver $Q_d$.}
\label{fig:Alternatibg-quiver}
\end{figure}

Write $\mathcal Q_d=\Db(\operatorname{rep}_{\C}Q_d)$.
Given a Bridgeland stability condition on an enhanced idempotent-complete 
smooth and proper triangulated 
category $\mathcal C$,
we construct one on the tensor product $\mathcal C\otimes\mathcal Q_d$ in Section~\ref{sec:all-degree}.
Its objects may be viewed as quiver diagrams, with an object of
$\mathcal C$ at each vertex and a morphism along each arrow.  We weigh central charges with the positive weights
\[
 r_j=\frac{\sin(j\pi/d)}{\sin(\pi/d)},\qquad 1\le j\le d-1,
\]
which form an eigenvector of the adjacency matrix of the underlying
unoriented graph.
This eigenvector identity turns  reflection functors of the quiver  into rotations
of the central charge.

We iterate this construction to obtain a stability condition on
$\mathcal Q_d^{\otimes M}$, the bounded derived category of
finite-dimensional modules over the tensor product of $M$ copies
of the path algebra of~$Q_d$.
The Thom--Sebastiani equivalence
\cite[Corollary~2.40]{BFK} identifies this tensor product with the
category of matrix factorizations of the Fermat polynomial
$$F=y_1^d+\ldots+y_M^d,$$
graded by the group
$L=\Z^M/\langle d(e_i-e_j)\rangle$, where $\deg(y_i)=e_i$.
The differences $e_i-e_j$ generate a finite subgroup
$H\simeq(\Z/d\Z)^{M-1}$.  Grading shifts by this subgroup preserve the stability
condition.  We then form the category of $H$-equivariant objects
in the sense of \cite[Definition~3.5]{Elagin}, which can be identified with the 
category of 
matrix factorizations $\MF^{\gr}(F)$ with its usual integer grading.

Finally, every degree-$d$ form over $\C$ is a sum of $d$th powers
of linear forms.  This realizes every smooth hypersurface as a linear
section of a sufficiently large Fermat hypersurface.
Proposition~\ref{prop:weighted-restriction} transfers stability condition one
hyperplane section at a time, preserving the Gepner symmetry,
until we reach the given hypersurface.
For Calabi--Yau Fermat polynomials, Okada~\cite[Sections~2--4]{OkadaMirror}
studied mirror stability conditions on tensor products of type~$A$
quiver categories and their interpretation through homological mirror symmetry.

As we were completing this paper, Wu's
preprint~\cite[Theorem~1.1]{WuKuznetsov} appeared on the arXiv.
Using stability conditions on products of elliptic curves, Wu
constructs numerical Gepner-type stability conditions on the Kuznetsov
components of smooth Fano hypersurfaces of degrees $3$, $4$, and $6$
and dimension at least three.  We also began with products of elliptic
curves, before realizing that tensor products of quiver categories
allow the construction to work in every degree.  The main work went
into constructing stability conditions on these tensor products.

\subsection*{Acknowledgements}
This work was supported in part by a grant of access to OpenAI models
through the ChatGPT for Academic Researchers program.
Codex Astra helped us find several steps in the construction and was
also used to edit the manuscript.
Any remaining errors are, of course, our own.
J.T. was partially supported by NSF grant DMS-2401387.

\section{Generalities on Categories}
\label{sec:background}

We will use enhanced models of triangulated categories.  
These models will allow us to construct
reflection functors, define tensor products, and
study finite group actions on categories.  
Stability conditions are defined on the underlying
triangulated categories.

\subsection{Stability conditions}

For a triangulated category $\mathcal C$, write $K_0(\mathcal C)$ for the
abelian group generated by its objects with relations
$[B]=[A]+[C]$ for distinguished triangles $A\to B\to C\to A[1]$.
The small triangulated categories on which we construct stability conditions
are $\C$-linear and proper, which means that the total space
$\bigoplus_k\Hom(E,F[k])$ is a finite dimensional $\C$-vector space for any two objects.
Thus their Euler pairing
\[
 \chi(E,F)=\sum_k(-1)^k\dim_{\C}\Hom(E,F[k])
\]
is defined.  The numerical Grothendieck group $K_{\num}(\mathcal C)$
is the quotient of $K_0(\mathcal C)$ by the left radical of the Euler pairing, which consists of
classes $v$ such that $\chi(v,w)=0$ for every $w$.
The categories considered here have Serre functors, recalled in
Section~\ref{subsec:categorical-tensor-products}; Serre duality
identifies the left and right radicals of the Euler pairing.
A numerical class map $v$ is a homomorphism from $K_0(\mathcal C)$ to a
finite-rank free abelian group $\Lambda$ that factors through $K_{\num}(\mathcal C)$.

A \emph{slicing} is a family of full additive subcategories
$\mathcal P(\phi)\subset\cC$, indexed by $\phi\in\R$, such that
$\mathcal P(\phi+1)=\mathcal P(\phi)[1]$, there are no morphisms
from a larger phase to a smaller phase, and every nonzero object has a
finite Harder--Narasimhan filtration by distinguished triangles whose nonzero factors lie
in $\mathcal P(\phi_1),\ldots,\mathcal P(\phi_s)$ with
$\phi_1>\cdots>\phi_s$.  For an interval $I\subset\R$, we denote by $\mathcal P(I)$  the
finite extension closure of subcategories $\mathcal P(\phi)$ with phases in $I$.

A \emph{stability condition} consists of a slicing and a central charge
$Z:K_0(\mathcal C)\to\C$ such that
\[
 Z(E)\in\R_{>0}e^{i\pi\phi}\quad
 \text{for }0\ne E\in\mathcal P(\phi).
\]
Each $\mathcal P(\phi)$ is an abelian category; its nonzero objects are called
\emph{semistable}, and its simple objects are called \emph{stable}.
Throughout the paper we require the \emph{support property}
with respect to a specified numerical class map
$v:K_0(\mathcal C)\to\Lambda$ and a factorization
$Z=\overline Z\circ v$.  This means that
\begin{equation}\label{eq:background-support}
 \|v(E)\|\le C|Z(E)|\qquad\text{for all semistable }E,
\end{equation}
for some norm on $\Lambda_{\R}$ and constant $C>0$;
see \cite[Definition~2.3 and Proposition~5.5]{BLMS}.

For an object $E$ with HN factors $A_i$, its \emph{mass} is defined as
$m_\sigma(E)=\sum_i|Z(A_i)|$.

Mass is unchanged by shifts and is additive under direct sums.
It also satisfies the \emph{mass triangle inequality}: for every
distinguished triangle
\[
 A\longrightarrow B\longrightarrow C\longrightarrow A[1],
\]
one has~\cite[Proposition~3.3]{Ikeda}
\begin{equation}\label{eq:mass-triangle-inequality}
 m_\sigma(B)\le m_\sigma(A)+m_\sigma(C).
\end{equation}

We will prove the following strong boundedness property for the
conditions constructed here.
We say that $\sigma$ has a \emph{total mass--Hom bound} if, for
every fixed $A$, there is a constant $C_A>0$ such that
\begin{equation}\label{eq:total-mass-hom}
 \sum_{n\in\Z}\dim_{\C}\Hom(A,E[n])
 \le C_A m_\sigma(E)\qquad\text{for every }E\in\mathcal C.
\end{equation}
In particular, \(\dim_{\mathbb C}\Hom(A,E)\le C_A m_\sigma(E)\) for every fixed \(A\) and every \(E\).

A t-structure is determined
by its \emph{aisle}, a full additive subcategory $U\subset\mathcal C$
such that $U[1]\subset U$ and every object $E$ admits a distinguished
{\em truncation triangle} $A\to E\to B\to A[1]$ with $A\in U$ and $B\in U^\perp$.
Here $U^\perp=\{B:\Hom(A,B)=0\text{ for all }A\in U\}$.
We put
$$\mathcal C^{\le n}=U[-n]\quad\hbox{\rm and}\quad
\mathcal C^{\ge n}=U^\perp[1-n].$$
The t-structure is \emph{bounded}
if every object belongs to $\mathcal C^{\ge a}\cap\mathcal C^{\le b}$
for some  $a\le b$.  Its \emph{heart}
$\mathcal A=U\cap U^\perp[1]$ is an abelian category, and truncation triangles define
cohomology objects $H^j_{\mathcal A}(E)\in\mathcal A$.

A slicing determines the aisle $U=\mathcal P((0,\infty))$ and the heart
$\mathcal A=\mathcal P((0,1])$.  Conversely, a stability function
on the heart of a bounded t-structure, taking nonzero objects into
$\{z:\operatorname{Im}z>0\}\cup\R_{<0}$ and satisfying the HN
property, determines a slicing~\cite[Proposition~5.3]{Bridgeland}.

We will repeatedly use tilting of hearts.  A \emph{torsion pair}
$(\mathcal T,\mathcal F)$ consists of full subcategories of
$\mathcal A$ such that $\Hom(\mathcal T,\mathcal F)=0$ and every
$E\in\mathcal A$ has an exact sequence
$0\to T\to E\to F\to0$ with $T\in\mathcal T$ and $F\in\mathcal F$.
The tilted heart is the extension closure
\[
 \mathcal A'=\langle\mathcal F[1],\mathcal T\rangle\subset\cC.
\]
Equivalently, $E\in\mathcal A'$ precisely when
$H^j_{\mathcal A}(E)=0$ for $j\ne-1,0$,
$H^{-1}_{\mathcal A}(E)\in\mathcal F$, and
$H^0_{\mathcal A}(E)\in\mathcal T$.
Its~aisle is $U'=\{E\in U:H^0_{\mathcal A}(E)\in\mathcal T\}$,
so $U[1]\subset U'\subset U$.

\subsection{Enhancements.}\label{subsec:categorical-tensor-products}
We follow
\cite[Sections~2.2--2.3 and~4.5--4.6]{Keller}
for dg categories and enhancements.
As a motivating example, let $E,F$ be bounded complexes of coherent
sheaves on a smooth projective variety $X$.  Suitable resolutions
produce a complex $\operatorname{RHom}(E,F)$ whose $k$-th cohomology group
is $\Hom_{\Db(X)}(E,F[k])$. A~\emph{differential graded
category}, or dg category, retains these complexes.  Thus every pair
of objects has a cochain complex of morphisms $\operatorname{Hom}^{\bullet}(E,F)$,
and composition is compatible with the differential:
$
 d(gf)=d(g)f+(-1)^{\deg g}g\,d(f)$.
Concrete enhancement models are described
in \cite[Section~5]{BLL}. 
The homotopy category $H^0(\mathcal E)$ of a dg category $\mathcal E$
has the same objects and
morphism spaces $H^0\operatorname{Hom}^{\bullet}(E,F)$.
For complexes the morphism differential is
$d(f)=d_Ff-(-1)^{\deg f}fd_E$; degree-zero cycles are chain maps,
and two chain maps have the same class in $H^0$ precisely when they
are homotopic.

A dg category is \emph{pretriangulated} if shifts and cones of closed
degree-zero maps are represented, up to homotopy, by objects of the
category.  Its homotopy category is then triangulated.  An
\emph{enhancement} of a triangulated category $\mathcal C$ is such a
dg category $\mathcal E$, together with an exact equivalence
$H^0(\mathcal E)\simeq\mathcal C$.
We use additive dg models, so finite direct sums exist at
the dg level.  
Adding direct summands is called \emph{idempotent completion}:
every endomorphism $e$ with $e^2=e$ becomes the projection onto a
direct summand.  Idempotent completion preserves the triangulated
structure~\cite[Definitions~1.1--1.2 and Theorem~1.5]{BalmerSchlichting}.
Our small triangulated categories are idempotent complete, with
additive pretriangulated dg enhancements.  
Alongside enhancements of the bounded derived category of coherent
sheaves, we will use the
concrete models supplied by quiver representations and matrix
factorizations.

A \emph{dg functor} acts by chain maps on morphism complexes and
preserves identities and composition.  It is a \emph{quasi-equivalence}
if these chain maps induce isomorphisms in every cohomological degree and the induced
functor on $H^0$ is essentially surjective.  An \emph{enhanced}
functor is supplied with a dg lift, after replacing the models by
quasi-equivalent ones when necessary.

An exact functor $F$ is \emph{conservative} if
$F(E)\simeq0$ implies $E\simeq0$.  Equivalently, a morphism is an
isomorphism whenever its image under $F$ is an isomorphism.

Let $\mathcal E$ be the chosen dg enhancement of a triangulated category $\mathcal C$.
A \emph{right dg module} over $\mathcal E$ is a dg functor from
$\mathcal E^{\mathrm{op}}$ to the dg category of cochain complexes
of complex vector spaces~\cite[Section~3.1]{Keller}.
The basic examples are the representable right dg modules
$h_E=\operatorname{Hom}^{\bullet}_{\mathcal E}(-,E)$.
The category $\Ind(\mathcal C)$ is obtained from
the homotopy category of dg modules by inverting maps $M\to N$ for
which $M(E)\to N(E)$ is a quasi-isomorphism for every object $E$
\cite[Section~3.2]{Keller}.
A dg module is called \emph{perfect} if its image in $\Ind(\mathcal C)$
belongs to the smallest full triangulated subcategory containing
the representable modules and closed under direct summands
\cite[Section~3.5]{Keller}.   We denote their category by $\Perf(\mathcal E)$.
The dg Yoneda embedding
$E\mapsto h_E$ therefore induces an equivalence
$\mathcal C\simeq\Perf(\mathcal E)$~\cite[Section~4.6]{Keller}.

An object $G\in\mathcal C$ is called a \emph{split generator} if every object
of $\mathcal C$ can be obtained from $G$ by finitely many shifts,
cones, finite direct sums, and direct summands.
An object $E$ of $\Ind(\mathcal C)$
is called  \emph{compact} if $\Hom(E,-)$ commutes with
arbitrary coproducts.  The compact objects are precisely the perfect
modules~\cite[Corollary~3.7]{Keller}, so the compact objects of
$\Ind(\mathcal C)$ recover $\mathcal C$.
For $\mathcal C=\Db(X)$ with $X$ smooth projective,
$\Ind(\mathcal C)$ is the unbounded derived category of
quasi-coherent sheaves on $X$~\cite[Section~3.6]{Keller}.

For dg enhancements $\mathcal E,\mathcal F$ of $\mathcal C,\mathcal D$,
their tensor product $\mathcal E\otimes_{\C}\mathcal F$ is the dg
category whose objects are pairs $E\boxtimes F$ and whose morphism
complexes are
\[
 \operatorname{Hom}^{\bullet}(E\boxtimes F,E'\boxtimes F')
 =\operatorname{Hom}^{\bullet}_{\mathcal E}(E,E')\otimes_{\C}
  \operatorname{Hom}^{\bullet}_{\mathcal F}(F,F').
\]
The differential and composition use the usual signs for tensor
products of complexes.  Over $\C$, this computes the derived dg
tensor product~\cite[Sections~2.3 and~4.3]{Keller}.
The \emph{Morita tensor product} of the enhanced triangulated
categories is, by definition,
$
 \mathcal C\otimes\mathcal D
 :=\Perf(\mathcal E\otimes_{\C}\mathcal F)$.
Here we take perfect right dg modules over the dg category just
defined.  Thus $\mathcal C\otimes\mathcal D$ is generated by the
representable modules of the pairs $E\boxtimes F$ under shifts,
cones, finite direct sums, and direct summands.  
For smooth projective varieties, this construction recovers the geometric product:
$\Db(X)\otimes\Db(Y)\simeq\Db(X\times Y)$
\cite[Theorem~5.5]{BLL}.

A dg category $\mathcal E$ is \emph{proper} if
$\bigoplus_k H^k\operatorname{Hom}^{\bullet}(E,F)$ is finite
dimensional for every pair of objects.  A dg \emph{bimodule} is a dg
functor of two variables, contravariant in one and covariant in the
other.  In our right-module convention, the \emph{diagonal bimodule}
is the right module over $\mathcal E^{\mathrm{op}}\otimes\mathcal E$
given by $(E,F)\mapsto\operatorname{Hom}^{\bullet}(F,E)$.
The category is \emph{smooth} if this bimodule is perfect
\cite[Section~4.7]{Keller}.  This is the categorical counterpart of
the perfectness of $\OO_{\Delta}$ on $X\times X$ for a smooth variety.

If $\mathcal E$ is smooth and proper then $\mathcal C$ has a split
generator $G$~\cite[Lemma~2.5]{EfimovCompletion}, and a dg module
$M$ is perfect if and only if
$\bigoplus_k H^k(M(G))$ is finite dimensional
\cite[Lemma~2.8 and Corollary~2.9]{ToenVaquie}.
Since $G$ split-generates $\mathcal C$, finiteness at $G$ is equivalent to
finiteness at every object of $\mathcal C$.

Smoothness and properness also give Serre duality: $\mathcal C$
has an exact autoequivalence $S_{\mathcal C}$, called its
\emph{Serre functor}, with natural isomorphisms
$
 \Hom(E,F)^*\simeq\Hom(F,S_{\mathcal C}E)$.
Existence follows from the preceding perfectness criterion and
\cite[Propositions~3.3 and~3.4(a)]{BondalKapranov}.
For $\Db(X)$ with $X$ smooth projective, Serre duality gives
$S_X=(-)\otimes\omega_X[\dim X]$
\cite[Example~3.2(2)]{BondalKapranov}.

\begin{lemma}\label{lem:mass-hom-generator}
Let $G$ split-generate $\mathcal C$.
It suffices to check~\eqref{eq:total-mass-hom} with $A=G$ and
$E$ semistable.  If $K_{\num}(\mathcal C)$ is a finite-rank lattice,
any numerical stability condition with a total mass--Hom bound
satisfies the support property on the full numerical lattice.
\end{lemma}

\begin{proof}
For every fixed $A$ there is $b_A>0$ such that
\[
 \sum_{n\in\Z}\dim_{\C}\Hom(A,E[n])
 \le b_A\sum_{n\in\Z}\dim_{\C}\Hom(G,E[n])
 \qquad\text{for all }E.
\]
Indeed, the objects $A$
admitting such a constant form a subcategory closed under shifts,
cones, finite sums, and direct summands, and this subcategory
contains $G$.  Applying the bound for $G$ to the HN factors of
$E$ and adding proves it for arbitrary $E$.

For the second assertion, choose objects $A_1,\ldots,A_s$ whose
Euler functionals separate numerical classes.  Then
\[
 \|u\|^2=\sum_i\chi(A_i,u)^2
\]
is a norm on $K_{\num}(\mathcal C)_\R$, and
\[
 \|[E]\|^2
 \le\sum_i\Bigl(\sum_{n\in\Z}\dim_{\C}\Hom(A_i,E[n])\Bigr)^2
 \le\Bigl(\sum_i C_{A_i}^2\Bigr)m_\sigma(E)^2.
\]
For semistable $E$, its mass is $|Z(E)|$, proving the assertion.
This is the Euler-pairing argument of~\cite[Theorem~2.1]{Cheng},
in its total-Hom formulation.
\end{proof}

For later use, this finite-rank hypothesis holds for $\Db(X)$
when $X$ is a smooth projective complex variety
\cite[Section~1.3]{Bridgeland}.  Indeed, the
Hirzebruch--Riemann--Roch theorem factors the Euler pairing through
the finite-dimensional space $H^{\mathrm{even}}(X,\Q)$;
finitely many integral Euler
functionals therefore embed $K_{\num}(X)$ into $\Z^s$.
The same hypothesis holds for admissible subcategories of $\Db(X)$,
and for categories obtained from them by adjoining finitely many
exceptional objects in a semiorthogonal decomposition.
For an admissible embedding $i$, its right adjoint $i^!$ is exact
and satisfies $i^!i\simeq\id$
\cite[Proposition~1.16]{BLL}.  The adjunction identity
$\chi(iA,B)=\chi(A,i^!B)$ and Serre duality show that both functors
descend to numerical groups.  Thus the numerical group of an
admissible subcategory is a direct summand of the ambient numerical
group.  More generally, the projection functors associated with a
semiorthogonal decomposition split the numerical group into those
of its components, \cite[Section~1.1]{OrlovQuivers}.

\subsection{Matrix factorizations}

We follow the grading and translation conventions of
\cite[Section~2.2, especially Definition~2.6]{TodaMF}.
Let $S=\C[x_1,\ldots,x_N]$ be positively graded, and let $W\in S$
be a nonzero homogeneous polynomial of degree $d$.  For a graded module we use
$M(a)_j=M_{j+a}$.  A \emph{graded matrix factorization} of $W$ is a
pair of finite-rank graded free $S$-modules and degree-zero maps
\begin{equation}\label{eq:graded-factorization}
 P_0\xrightarrow{p_0}P_1\xrightarrow{p_1}P_0(d)
\end{equation}
such that $p_1p_0=W\operatorname{id}_{P_0}$ and
$p_0(d)p_1=W\operatorname{id}_{P_1}$.  Modulo $W$, we can interpret 
\eqref{eq:graded-factorization}
as a
two-periodic complex, with an internal twist by $d$ every two steps.

A morphism $f:P\to Q$ consists of degree-zero maps
$f_i:P_i\to Q_i$ making both squares commute:
\[
 \begin{tikzcd}[column sep=large,row sep=large]
 P_0 \arrow[r,"p_0"] \arrow[d,"f_0"'] &
 P_1 \arrow[r,"p_1"] \arrow[d,"f_1"'] &
 P_0(d) \arrow[d,"f_0(d)"] \\
 Q_0 \arrow[r,"q_0"'] & Q_1 \arrow[r,"q_1"'] & Q_0(d).
 \end{tikzcd}
\]

Two morphisms $f,g:P\to Q$ are \emph{homotopic} if there are
degree-zero maps $h_0:P_0\to Q_1(-d)$ and $h_1:P_1\to Q_0$,
shown by the dashed arrows below, satisfying the two identities
that follow:
\[
 \begin{tikzcd}[column sep=large,row sep=large]
 & P_0 \arrow[r,"p_0"] \arrow[dl,dashed,"h_0"'] &
 P_1 \arrow[r,"p_1"] \arrow[dl,dashed,"h_1"'] &
 P_0(d) \arrow[dl,dashed,"h_0(d)"'] \\
 Q_1(-d) \arrow[r,"q_1(-d)"'] &
 Q_0 \arrow[r,"q_0"'] & Q_1 \arrow[r,"q_1"'] & Q_0(d).
 \end{tikzcd}
\]
\[
 \begin{split}
 f_0-g_0&=q_1(-d)h_0+h_1p_0,\\
 f_1-g_1&=q_0h_1+h_0(d)p_1.
 \end{split}
\]

These morphisms and homotopies arise from the standard dg category
of matrix factorizations~\cite[Definitions~2.7 and~2.12]{BFK}.
Its degree-zero cycles are morphisms of factorizations, and its
degree-zero boundaries are nullhomotopic morphisms.  Thus $H^0$
gives morphisms modulo homotopy.
We write
$\MF^{\gr}(W)$ for the resulting
homotopy category, with direct summands added.  Its translation functor
$[1]$ sends \eqref{eq:graded-factorization} to
$
 P_1\xrightarrow{-p_1}P_0(d)
       \xrightarrow{-p_0(d)}P_1(d)$.
Applying this translation twice shifts the internal grading by $d$.
The \emph{internal grading shift} $\tau$ sends a factorization $P$
to $P(1)$, obtained by shifting both modules and both maps by one.
More generally, $\tau^aP=P(a)$.  Thus
\begin{equation}\label{eq:mf-periodicity}
 \tau=(1),\qquad \tau^d=(d)\simeq[2].
\end{equation}
This explains the phase increment $2/d$ in Theorem~\ref{thm:main}:
the grading shift is a categorical $d$th root of the double shift.
There is another useful description.  Put $A=S/(W)$ and let
$\Db_{\gr}(A)$ be the bounded derived category of finitely generated
graded $A$-modules.  The \emph{graded singularity category} is the
Verdier quotient
\[
 \Dsg^{\gr}(A)=\Db_{\gr}(A)/\Perf^{\gr}(A).
\]
Here $\Perf^{\gr}(A)$ consists of complexes quasi-isomorphic to
bounded complexes of finitely generated graded projective modules.
The module
$\operatorname{coker}(p_0)$
of a matrix factorization
 is annihilated by $W$, hence is a graded
$A$-module.  Sending a factorization to this cokernel induces the
equivalence
$\MF^{\gr}(W)\simeq\Dsg^{\gr}(A)$, with idempotent completions
included~\cite[Theorem~2.7]{TodaMF}.

Suppose now that $\deg x_i=1$ for every $i$ and that the projective
hypersurface $X=(W=0)\subset\Pj^{N-1}$ is smooth.  Orlov's
theorem~\cite[Theorem~3.11]{Orlov} compares its derived category with
$\MF^{\gr}(W)$, according to the sign of $N-d$.
If $N>d$, it identifies $\MF^{\gr}(W)$ with the residual category
$\Ku(X)$.  With the normalization used here, the semiorthogonal
decomposition is
\[
 \Db(X)=\langle\Ku(X),\OO_X,\OO_X(1),\ldots,
                                      \OO_X(N-d-1)\rangle.
\]
Equivalently,
\[
 \Ku(X)=\{E\in\Db(X):\Hom^{\bullet}(\OO_X(j),E)=0
                         \text{ for }0\le j<N-d\}.
\]
For $N=d$, $\MF^{\gr}(W)$ is equivalent to the whole of $\Db(X)$.
For $N<d$, it admits a semiorthogonal decomposition into a copy of
$\Db(X)$ and $d-N$ exceptional objects.

When $N>d$, the Orlov equivalence can be chosen so that $\tau$
corresponds to the \emph{rotation functor}: tensoring by $\OO_X(1)$
followed by left mutation through $\OO_X$
\cite[Definition~6.11 and Theorem~6.13]{BFKHodge}.

The weighted version used in Section~\ref{sec:hyperplane}
replaces $X$ by the projective hypersurface stack
$[(W^{-1}(0)\setminus\{0\})/\C^*]$, where $\C^*$ acts with the
prescribed positive weights.  The three cases are then determined
by the sign of $\sum_i\deg x_i-d$
\cite[Theorem~2.5 and Remark~3.12]{Orlov}.

\subsection{Thom--Sebastiani theorem}

Let $H$ be a finite group and $\mathcal C$ an enhanced triangulated
category.  A \emph{coherent left action} of $H$ on $\mathcal C$
consists of exact autoequivalences $F_h:\mathcal C\to\mathcal C$
for $h\in H$ and natural isomorphisms
$\alpha_{g,h}:F_gF_h\xrightarrow{\sim}F_{gh}$.
We normalize $F_1=\operatorname{id}$ and
$\alpha_{1,h}=\alpha_{h,1}=\operatorname{id}$.
Associativity is the identity
$
 \alpha_{gh,k}\circ(\alpha_{g,h}F_k)
 =\alpha_{g,hk}\circ(F_g\alpha_{h,k})$
between the two maps $F_gF_hF_k\to F_{ghk}$,
$g,h,k\in H$.
Fix an additive pretriangulated dg enhancement $\mathcal E$ of
$\mathcal C$, so that $H^0(\mathcal E)\simeq\mathcal C$.
We assume that the $F_h$ lift to dg autoequivalences of
$\mathcal E$ and that the $\alpha_{g,h}$ lift to closed degree-zero
natural isomorphisms satisfying the same identities
\cite[Definition~8.1]{Elagin}.

An \emph{equivariant object} is a pair $(E,\rho)$, where
$E\in\mathcal C$ and the isomorphisms
$\rho_h:F_hE\xrightarrow{\sim}E$ satisfy
\[
 \rho_1=\operatorname{id}_E,\qquad
 \rho_{gh}\circ\alpha_{g,h,E}=\rho_g\circ F_g(\rho_h).
\]
A morphism $(E,\rho)\to(E',\rho')$ is a morphism $f:E\to E'$
in $\mathcal C$ satisfying
$f\circ\rho_h=\rho'_h\circ F_h(f)$ for every $h\in H$.
These objects and morphisms form the \emph{equivariant category}
$\mathcal C^H$; compare
\cite[Definitions~3.1 and~3.5]{Elagin}, written here for left actions.
The dg category $\mathcal E^H$ is defined similarly, using closed
degree-zero equivariance isomorphisms and complexes of equivariant
maps.  Its category of perfect dg modules
provides an enhancement of $\mathcal C^H$, and hence a triangulated
structure~\cite[Theorem~8.7]{Elagin}.  This uses our standing
assumption that $\mathcal C$ is idempotent complete.

The forgetful functor $U:\mathcal C^H\to\mathcal C$ sends
$(E,\rho)$ to $E$.  In the other direction, induction sends
$X\in\mathcal C$ to
$ I(X)=\bigoplus_{h\in H}F_hX$.
The maps $\alpha_{g,h}$ permute these summands and define its
equivariant structure.  Since $H$ is finite, $I$ is both a left
and a right adjoint to $U$~\cite[Lemma~3.8]{Elagin}.
The adjunction map $I(U(E,\rho))\to(E,\rho)$ has components
$\rho_h$.  Over $\C$, it has an equivariant section with components
$|H|^{-1}\rho_h^{-1}$.  Thus every equivariant object is a direct
summand of an induced object~\cite[Lemma~3.9]{Elagin}.
We say that the action preserves a stability condition
$\sigma=(Z,\mathcal P)$ on $\mathcal C$ if
\[
 Z\circ(F_h)_*=Z,\qquad F_h\mathcal P(\phi)=\mathcal P(\phi)
 \quad(h\in H,\ \phi\in\R).
\]
For an equivariant object $(E,\rho)$, write $E_j$ for the
successive terms of its HN filtration.  Each term is a canonical
phase truncation of $E$, retaining the HN factors above a chosen
phase.  Since $F_h$ preserves phases, $\rho_h$ induces a unique
isomorphism $\rho_h^j:F_hE_j\to E_j$ compatible with the maps to $E$.
The two compositions
$\rho_{gh}^j\circ\alpha_{g,h,E_j}$,
 $\rho_g^j\circ F_g(\rho_h^j)$
are induced by the same map $F_gF_hE\to E$, so they agree by
uniqueness of the truncation.  The semistable factors inherit
compatible structures in the same way.  Thus the HN filtration
lifts to $\mathcal C^H$.  The induced charge and slicing are
\[
 Z^H=Z\circ U_*,\qquad
 \mathcal P^H(\phi)=\{(E,\rho):E\in\mathcal P(\phi)\}.
\]
Lemma~\ref{lem:equivariant-support} proves the support property
and numerical property; compare
\cite[Theorem~2.16 and Remark~2.17]{LMPSZ}.

The \emph{Thom--Sebastiani theorem} relates a tensor product of
matrix-factorization categories to the sum of their polynomials.
The sum carries a grading that records the degree in each set
of variables.
Let $f(x)$ and $g(y)$ be degree-$d$ homogeneous polynomials in
disjoint sets of variables, each with an isolated critical point
at the origin.  Give the variables $x_i$ degree $e_1$ and the
variables $y_j$ degree $e_2$ in the group
\[
 L=(\Z e_1\oplus\Z e_2)/\langle d(e_1-e_2)\rangle
       \simeq\Z\oplus\Z/d.
\]
Put $S=\C[x_1,\ldots,x_s,y_1,\ldots,y_t]$ and
$\eta=de_1=de_2\in L$, the common degree of $f$ and $g$.
An $L$-graded $S$-module is a decomposition
$M=\bigoplus_{\ell\in L}M_\ell$ with
$S_aM_b\subseteq M_{a+b}$.  Its grading shift is
$M(a)_\ell=M_{\ell+a}$ for $a\in L$.
An \emph{$L$-graded matrix factorization} of $f+g$ consists of
two finite-rank graded free modules $P_0,P_1$, each a finite direct
sum of shifts $S(a)$, and degree-zero maps
\[
               P_0\xrightarrow{p_0}P_1
                    \xrightarrow{p_1}P_0(\eta)
\]
such that
$
 p_1p_0=(f+g)\operatorname{id}_{P_0},$
 $p_0(\eta)p_1=(f+g)\operatorname{id}_{P_1}$.
Morphisms, homotopies, and the dg enhancement are defined as
above, with grading shifts in $L$ and $d$ replaced by $\eta$.
The homotopy category, with direct summands added, is denoted
$\MF^L(f+g)$.  The cohomological shift and the grading shift satisfy
$[2]\simeq(\eta)$, just as for integer gradings.

The Thom--Sebastiani theorem~\cite[Theorem~1.1]{BFK} gives an
equivalence
\[
 \MF^{\gr}(f)\otimes\MF^{\gr}(g)\simeq\MF^L(f+g).
\]
In the tensor product, degrees $a$ and $b$ combine to give $ae_1+be_2$.  Thus shifting the first or second
factor by one shifts the $L$-grading by $e_1$ or $e_2$, respectively.
Section~\ref{sec:thom-sebastiani} treats the passage from this
grading to ordinary integer grading.

\section{Stability with quiver coefficients}
\label{sec:all-degree}

Throughout this section, $\mathcal C$ is an
idempotent-complete triangulated category with a smooth proper dg
enhancement $\cE$.
Given a stability condition on $\cC$, we will study
stability conditions on its tensor product $\cC_Q$ with a quiver
category $\Db(\operatorname{rep}_{\C}Q)$,
where $Q$ is a
finite acyclic quiver.  

An object of $\cC_Q$ is given by {\em a quiver diagram}. 
It specifies an object of $\cC$ at each vertex of $Q$ 
and a morphism along each arrow.  To form cones of maps between
such diagrams functorially, we work with the dg modules introduced
in Section~\ref{subsec:categorical-tensor-products}.
A \emph{dg quiver diagram} consists of a right dg $\mathcal E$-module
$M_i$ at each vertex and a chain map $M_a:M_i\to M_j$ along each
arrow $a:i\to j$.  Morphisms are families of chain maps commuting
with the arrows.  We identify homotopic morphisms and invert maps
that are quasi-isomorphisms at every vertex, after evaluation at
every object of $\mathcal E$.
Let $\C Q$
be the path category of the quiver: its objects are vertices and its morphisms
are linear combinations of directed paths, all in degree zero.
A dg quiver diagram is precisely a right dg module over
$\mathcal E\otimes(\C Q)^{\mathrm{op}}$. 
For the dg-module background, see
\cite[Sections~3.1--3.2 and~3.5]{Keller}.

This construction belongs to the general theory of coherent quiver
diagrams studied by Groth--\v{S}\v{t}ov\'i\v{c}ek
\cite[Sections~2--3]{GrothStovicek}; see also their treatment of
type~$A$ in~\cite{GrothStovicekA}.
In the literature on stability conditions, Kawatani
\cite[Section~2.2]{KawataniMorphisms} treats the case of $A_2$
and explains why forming diagrams in the enhancement retains
more information than forming them in the triangulated category.
Haiden--Katzarkov--Simpson
\cite[Sections~5.1 and~5.3]{HKScoefficients} use $A_n$-diagrams
with categorical coefficients to describe Fukaya categories of
marked disks and their tensor-product models.

Let $P_i$ be the standard projective representation of the quiver $Q$:
its vector space at vertex $j$ has a basis indexed by paths from $i$ to $j$. 
For an arrow \(a:j\to k\), the linear map
$
(P_i)_j\longrightarrow(P_i)_k
$
sends the basis vector indexed by a path \(p:i\to j\) to the basis vector indexed by the path
$
a\circ p:i\longrightarrow j\longrightarrow k
$. For any dg quiver diagram $M$, the standard quiver resolution
\cite[Section~1]{CrawleyBoeveyQuivers} gives a short exact sequence
\begin{equation}\label{eq:quiver-projective-resolution}
 0\longrightarrow\bigoplus_{a:i\to j}P_j\otimes_{\C}M_i
 \longrightarrow\bigoplus_iP_i\otimes_{\C}M_i
 \longrightarrow M\longrightarrow0.
\end{equation}
To describe the maps, fix a vertex $v$.  The middle term at $v$
contains one copy of $M_i$ for each directed path $p:i\to v$.
On this copy, the right-hand map is $M_p:M_i\to M_v$, the
composition of the arrow maps along $p$.  The path of length zero
at $v$ contributes a copy of $M_v$ on which this map is the
identity, so the right-hand map is surjective.
The left-hand term at $v$ contains one copy of $M_i$ for each
arrow $a:i\to j$ followed by a path $p:j\to v$.
This copy maps into two summands of the middle term: by the
identity into the copy of $M_i$ indexed by the concatenated path
$p\circ a$, and by $-M_a$ into the copy of $M_j$ indexed by $p$.
Their images in $M_v$ cancel, since $M_{p\circ a}=M_p\circ M_a$.
Thus the left-hand map records the relation between two ways of
following the same path: either keep its first arrow as part of
the path, or apply that arrow to the coefficient object first.

The resolution gives a criterion for perfectness in terms of the
vertex modules.

\begin{lemma}\label{lem:quiver-diagram-perfectness}
Let $Q$ be a finite acyclic quiver and $M$ a dg quiver diagram
of shape $Q$.  Then $M$ is perfect as a right module over
$\mathcal E\otimes(\C Q)^{\mathrm{op}}$ if and only if every
vertex module $M_i$ is perfect over $\mathcal E$.
\end{lemma}

\begin{proof}
Perfect modules are obtained from representable modules
by taking shifts, cones, finite direct sums, and direct summands.
For an object $A\in\mathcal E$, write
$h_A=\operatorname{Hom}^{\bullet}_{\mathcal E}(-,A)$.
The objects of the tensor-product dg category
$\mathcal E\otimes(\C Q)^{\mathrm{op}}$ are pairs $(A,j)$,
consisting of an object $A$ of $\mathcal E$ and a vertex $j$ of $Q$.
The right module represented by $(A,j)$ assigns to another object
$(B,v)$ the complex
\[
\operatorname{Hom}^{\bullet}_{\mathcal E\otimes(\C Q)^{\mathrm{op}}}
\bigl((B,v),(A,j)\bigr)
=\operatorname{Hom}^{\bullet}_{\mathcal E}(B,A)
  \otimes_{\C}\operatorname{Hom}_{\C Q}(j,v)
=\bigoplus_{p:j\to v}h_A(B).
\]
Thus, viewed as a quiver diagram of $\mathcal E$-modules, this
representable module has one copy of $h_A$ at vertex $v$ for each
path $j\to v$.  Its arrow maps concatenate paths, acting by the
identity on $h_A$.  This is exactly the diagram
$P_j\otimes_{\C}h_A$.

Since $P_j\otimes_{\C}-$ preserves shifts, cones, finite direct sums, and direct summands, it sends
any perfect $\mathcal E$-module to a perfect module over this
tensor product.  Thus, if every $M_i$ is perfect, so is each
term $P_j\otimes_{\C}M_i$ in the resolution \eqref{eq:quiver-projective-resolution}.
Both sums are finite, because $Q$ has finitely many vertices and
arrows.  The resolution identifies $M$, in the derived category,
with the cone of the map between these two perfect modules.
It follows that $M$ is perfect as well.

Conversely, evaluate a representable diagram at a vertex $v$:
$
 (P_j\otimes_{\C}h_A)_v
   =\bigoplus_{p:j\to v}h_A$.
There is one summand for each directed path from $j$ to $v$.
Since $Q$ is finite and acyclic, there are only finitely many such
paths, so this is a perfect $\mathcal E$-module.
Evaluation at $v$  preserves shifts, cones, finite direct
sums, and direct summands.  It therefore sends every perfect
module over $\mathcal E\otimes(\C Q)^{\mathrm{op}}$ to a perfect
$\mathcal E$-module.  In particular, if $M$ is perfect, then every
vertex module $M_v$ is perfect.
\end{proof}

Lemma~\ref{lem:quiver-diagram-perfectness} is a special case of
To\"en--Vaqui\'e's criterion
\cite[Definition~2.7 and Lemma~2.8(3)]{ToenVaquie}.
For dg categories $T$ and $T'$, they call a right module over
$T\otimes^{\mathbf L}T'$ \emph{pseudo-perfect relative to $T'$}
if its evaluation at every object of $T$ is a perfect $T'$-module.
When $T$ is smooth and proper, this condition is equivalent
to perfectness over $T\otimes^{\mathbf L}T'$.
In our setting, take $T=(\C Q)^{\mathrm{op}}$ and $T'=\mathcal E$.
Evaluation at an object of $T$ is simply evaluation at a vertex
of $Q$, so relative pseudo-perfectness means  that each
$M_i$ is perfect.  The quiver factor is proper because
its path spaces are finite-dimensional, and smooth because its
diagonal has a finite resolution by finitely generated projective
bimodules.  The explicit quiver
argument above makes this general criterion concrete.  It also
shows that no smoothness or properness assumption on $\mathcal E$
is needed for this particular lemma.

Thus quiver diagrams with perfect vertex modules give the category
\[
 \mathcal C_Q=
 \Perf\bigl(\mathcal E\otimes(\C Q)^{\mathrm{op}}\bigr)
 \simeq\mathcal C\otimes\Db(\operatorname{rep}_{\C}Q).
\]

\begin{lemma}\label{lem:quiver-coefficient-heart}
A bounded t-structure on $\mathcal C$ induces a bounded t-structure
on $\mathcal C_Q$ with
\[
 \mathcal C_Q^{\le0}=\{M:M_i\in\mathcal C^{\le0}\},\qquad
 \mathcal C_Q^{\ge0}=\{M:M_i\in\mathcal C^{\ge0}\}.
\]
If $\mathcal H$ is the heart on $\mathcal C$, the induced heart is
naturally equivalent to $\operatorname{Rep}(Q,\mathcal H)$.
Its kernels, cokernels, and exact sequences are computed at each vertex.
We call this the standard heart on $\mathcal C_Q$ associated with
$\mathcal H$.
\end{lemma}

\begin{proof}
Extend the t-structure on $\mathcal C$ to $\Ind(\mathcal C)$
using~\cite[Lemma~2.1.1]{Polishchuk}.
By~\cite[Theorem~A]{SSVHearts}, the extended t-structure induces
a pointwise t-structure on $\Ind(\mathcal C_Q)$, viewed as the
category of quiver diagrams with coefficients in
$\Ind(\mathcal C)$.
Its heart consists of quiver diagrams in the heart of the extended
t-structure.
Since truncations preserve $\mathcal C$, this t-structure
restricts to diagrams with perfect vertices, hence to
$\mathcal C_Q$.
It is bounded because $Q$ has finitely many vertices.
Its heart is $\operatorname{Rep}(Q,\mathcal H)$, where kernels,
cokernels, and exact sequences are computed at each vertex.
\end{proof}

\begin{proposition}[Takahashi]\label{prop:one-variable-model}
Let $d\ge2$, and let $Q$ be the quiver $A_{d-1}$ with any orientation.
Write $\operatorname{rep}_{\C}Q$ for its category of finite-dimensional
complex representations.  There is an enhanced equivalence
\begin{equation}\label{eq:one-variable-equivalence}
 \Phi:\MF^{\gr}(x^d)\xrightarrow{\sim}
                 \Db(\operatorname{rep}_{\C}Q)
\end{equation}
that sends the internal grading shift $\tau$ of $\MF^{\gr}(x^d)$ to
$S_Q^{-1}[1]$, where $S_Q$ is the Serre functor on
$\Db(\operatorname{rep}_{\C}Q)$.\break
In particular,
\begin{equation}\label{eq:one-variable-serre}
 S_{\MF}\simeq\tau^{-1}[1],\qquad S_{\MF}^{d}\simeq[d-2].
\end{equation}
\end{proposition}

\begin{proof}
We recall Takahashi's construction~\cite{Takahashi}
for the linear orientation of the quiver, with all arrows pointing in the same direction.
Put $R=\C[x]$, with $\deg x=1$.  The matrix factorizations
\[
 E_\ell=\bigl(R\xrightarrow{x^\ell}R(\ell)
                 \xrightarrow{x^{d-\ell}}R(d)\bigr),
                  \qquad 1\le\ell\le d-1,
\]
split-generate $\MF^{\gr}(x^d)$ and satisfy
\begin{equation}\label{eq:one-variable-exceptional}
 \Hom(E_k,E_\ell[n])=
 \begin{cases}
  \C,&n=0\text{ and }k\le\ell,\\
  0,&\text{otherwise}.
 \end{cases}
\end{equation}
For $k\le\ell$, the maps $f_{k\ell}=(1,x^{\ell-k})$ span
$\Hom(E_k,E_\ell)$ and satisfy
$f_{\ell m}f_{k\ell}=f_{km}$.  

Thus $(E_1,\ldots,E_{d-1})$ is a full strong exceptional collection.
Put $T=\bigoplus_{\ell=1}^{d-1}E_\ell$ and
$B=\operatorname{End}(T)$.  The composition rule identifies
$B$ with the path algebra of the quiver $1\to2\to\cdots\to d-1$.
The~equivalence theorem for full strong exceptional collections
\cite[Corollary~1.9]{OrlovQuivers}, with its dg enhancement
\cite[Theorem~3.8(b) and Section~3.8]{Keller}, gives an enhanced equivalence
\[
 \Psi:\MF^{\gr}(x^d)\xrightarrow{\sim}\Perf(B),
 \qquad E\longmapsto\operatorname{RHom}(T,E).
\]
The same theorem identifies $\Perf(B)$ with the bounded derived
category of finite-dimensional right $B$-modules.
The right action is by precomposition.  Thus, writing
$V_i=\operatorname{RHom}(E_i,E)$, the arrows reverse direction.
For example, when $d=5$ the two chains are
\[
 \begin{tikzcd}[column sep=large,row sep=large]
 E_1 \arrow[r,"f_{12}"] & E_2 \arrow[r,"f_{23}"] &
 E_3 \arrow[r,"f_{34}"] & E_4 \\
 V_1 & V_2 \arrow[l,"f_{12}^{*}"'] &
 V_3 \arrow[l,"f_{23}^{*}"'] & V_4 \arrow[l,"f_{34}^{*}"'].
 \end{tikzcd}
\]
Here $f_{i,i+1}^{*}$ means precomposition with $f_{i,i+1}$.
Relabeling vertex $i$ as $d-i$ makes the bottom chain point to the
right and gives \eqref{eq:one-variable-equivalence} for the linear
orientation.

Takahashi proves the triangulated equivalence
\eqref{eq:one-variable-equivalence} in~\cite[Theorem~3.1]{Takahashi}.
His trace pairing identifies $F=\tau^{-1}[1]$ with the Serre
functor $S_{\MF}$~\cite[Lemma~3.5 and Remark~3.6]{Takahashi}.
On $\Perf(B)$, the Serre functor is given by~\cite[Section~5.5]{KellerQuivers}
\[
 N\longmapsto N\otimes_B^{\mathbf L}D(B),
 \qquad D(B)=\Hom_{\C}(B,\C).
\]
Under $\Psi$, the functor $F$ becomes
$N\mapsto N\otimes_B^{\mathbf L}M$, where
$M=\operatorname{RHom}(T,FT)$ and the left $B$-action is induced
by $F$ on endomorphisms of $T$.
To compare the functors at the enhanced level, it therefore suffices
to identify $M$ with $D(B)$.
Serre duality for $T[-n]$ and $T$ gives
\[
 H^n(M)=\Hom(T,FT[n])
       \simeq\Hom(T[-n],FT)
       \simeq\Hom(T,T[-n])^*.
\]
By~\eqref{eq:one-variable-exceptional}, this vanishes for $n\ne0$,
while $H^0(M)\simeq\operatorname{End}(T)^*=D(B)$.
Naturality makes this identification compatible with both $B$-actions.
Truncation therefore gives $M\simeq D(B)$ as bimodules, proving
the required compatibility with the enhancements.

To pass to arbitrary orientations of the quiver, we use the derived
Bernstein--Gelfand--Ponomarev reflection functors.
Let $i$ be a source of $Q$, so all arrows incident to $i$ point
outwards, and let $Q'$ be the quiver obtained by reversing these arrows.
For a complex $V$ of $Q$-representations, the reflection replaces
$V_i$ by $\Cone(V_i\to\bigoplus_{i\to j}V_j)$ and leaves the
other vertex complexes unchanged.  The reversed arrows are the
natural maps from the neighboring complexes $V_j$ to this cone.
At an interior source vertex of $Q$ with outgoing arrows $a$ and $b$, put $u=(a,b)$ and $C=\Cone(u)$.
The reflection has the following local picture:
\[
 \begin{tikzcd}[column sep=1.8em]
 V_{i-1} & V_i \arrow[l,"a"'] \arrow[r,"b"] & V_{i+1}
 \end{tikzcd}
 \quad\xmapsto{\ G\ }\quad
 \begin{tikzcd}[column sep=1.8em]
 V_{i-1} \arrow[r,"j_-"] & C & V_{i+1} \arrow[l,"j_+"'].
 \end{tikzcd}
\]
Both new arrows come from the same distinguished triangle:
\begin{equation}\label{wefvwetgwetgwetg}
 V_i\xrightarrow{\ u\ }V_{i-1}\oplus V_{i+1}
 \xrightarrow{\ (j_-,j_+)\ }C\longrightarrow V_i[1].
\end{equation}
At an endpoint there is only one neighboring term.
In $Q'$, the vertex $i$ is a \emph{sink}: all adjacent arrows point
towards it.  For a complex $W$ of $Q'$-representations, the inverse
is given at $i$ by
\begin{equation}\label{eq:inverse-vertex-reflection}
 (G^{-1}W)_i\simeq
 \Cone\left(\bigoplus_{j\to i}W_j\longrightarrow W_i\right)[-1].
\end{equation}
Rotating the triangle \eqref{wefvwetgwetgwetg} shows that this recovers $V_i$
and its outgoing arrows.  The two constructions give an equivalence of categories
\[
 G:\Db(\operatorname{rep}_{\C}Q)\xrightarrow{\sim}
                      \Db(\operatorname{rep}_{\C}Q').
\]
For the classical construction and its realization as a derived
tensor product with a complex of bimodules, see
\cite[Section~9.1]{KellerQuivers}.  This realization also gives
$G$ a dg lift.
We can also see explicitly that this equivalence intertwines the Serre functors.
For $X,Y\in\Db(\operatorname{rep}_{\C}Q)$, there are natural
quasi-isomorphisms
\[
 \operatorname{RHom}_{Q'}(GX,GY)^\vee
 \simeq\operatorname{RHom}_{Q}(X,Y)^\vee
 \simeq\operatorname{RHom}_{Q}(Y,S_QX)
 \simeq\operatorname{RHom}_{Q'}(GY,GS_QX).
\]
Here $\vee$ denotes the dual complex of vector spaces.
The first and last maps come from the dg lift of $G$, and the
middle map is Serre duality.  Thus $GS_QX$ and $S_{Q'}GX$
represent the same dg module, giving $GS_Q\simeq S_{Q'}G$
at the enhanced level.

Every orientation of a path can be made linear: move the leftmost
left-pointing arrow to the left endpoint by sink reflections,
reverse it by reflecting at that endpoint, and repeat.
For example, consider a path on four vertices.  In the following
diagrams, a label $i$ denotes reflection at the sink $i$, and the
blue arrow is the left-pointing arrow being moved:
\[
 (1\to2\mathrel{\color{RoyalBlue}\leftarrow}3\leftarrow4)
 \xmapsto{\,2\,}
 (1\mathrel{\color{RoyalBlue}\leftarrow}2\to3\leftarrow4)
 \xmapsto{\,1\,}
 (1\to2\to3\leftarrow4).
\]
Repeat with the remaining
left-pointing arrow:
\[
 (1\to2\to3\mathrel{\color{RoyalBlue}\leftarrow}4)
 \xmapsto{\,3\,}
 (1\to2\mathrel{\color{RoyalBlue}\leftarrow}3\to4)
 \xmapsto{\,2\,}
 (1\mathrel{\color{RoyalBlue}\leftarrow}2\to3\to4)
 \xmapsto{\,1\,}
 (1\to2\to3\to4).
\]
Each pass removes one left-pointing arrow, so the procedure terminates.
Consequently, $\tau$ corresponds to $S_Q^{-1}[1]$ for every
orientation.  Finally, $\tau^d\simeq[2]$ gives \eqref{eq:one-variable-serre}.
\end{proof}

\begin{figure}[htbp]
\centering
\begin{tikzpicture}[x=1cm,y=1cm,>=Stealth,
  vertex/.style={circle,draw=black!55,minimum size=6.8mm,
                 inner sep=1pt,font=\small},
  source/.style={vertex,fill=MidnightBlue!10},
  target/.style={vertex,fill=black!5},
  targetsum/.style={rectangle,rounded corners=2pt,draw=black!55,
                   fill=black!5,minimum width=23mm,
                   minimum height=7mm,inner sep=3pt,font=\small},
  every edge/.style={draw,->}]
\node at (1.725,1.95) {$Q_5\bigsqcup Q_5^{\rm op}$};
\node[source] (qX1) at (0,.55) {$X_1$};
\node[target] (qY2) at (1.15,.55) {$Y_2$};
\node[source] (qX3) at (2.3,.55) {$X_3$};
\node[target] (qY4) at (3.45,.55) {$Y_4$};
\node[target] (qY1) at (0,-.55) {$Y_1$};
\node[source] (qX2) at (1.15,-.55) {$X_2$};
\node[target] (qY3) at (2.3,-.55) {$Y_3$};
\node[source] (qX4) at (3.45,-.55) {$X_4$};
\path (qX1) edge (qY2) (qX3) edge (qY2) (qX3) edge (qY4)
      (qX2) edge (qY1) (qX2) edge (qY3) (qX4) edge (qY3);
\node[font=\small] at (1.725,-1.45) {$X=(X_1,X_2,X_3,X_4)$};
\node[font=\small] at (1.725,-1.95) {$Y=(Y_1,Y_2,Y_3,Y_4)$};
\node at (4.65,0) {$\Longleftrightarrow$};
\node at (5.9,1.95) {$X$};
\node at (9.1,1.95) {$FY$};
\node[source] (tX1) at (5.9,1.35) {$X_1$};
\node[source] (tX2) at (5.9,.45) {$X_2$};
\node[source] (tX3) at (5.9,-.45) {$X_3$};
\node[source] (tX4) at (5.9,-1.35) {$X_4$};
\node[targetsum] (tY1) at (9.1,1.35) {$Y_2$};
\node[targetsum] (tY2) at (9.1,.45) {$Y_1\oplus Y_3$};
\node[targetsum] (tY3) at (9.1,-.45) {$Y_2\oplus Y_4$};
\node[targetsum] (tY4) at (9.1,-1.35) {$Y_3$};
\path (tX1) edge node[above,font=\small] {$m_1$} (tY1)
      (tX2) edge node[above,font=\small] {$m_2$} (tY2)
      (tX3) edge node[above,font=\small] {$m_3$} (tY3)
      (tX4) edge node[above,font=\small] {$m_4$} (tY4);
\node[font=\small] at (7.5,-1.95) {$m=(m_1,m_2,m_3,m_4):X\longrightarrow FY$};
\end{tikzpicture}
\caption{A quiver diagram for $d=5$.}
\label{fig:quiver-tuples}
\end{figure}
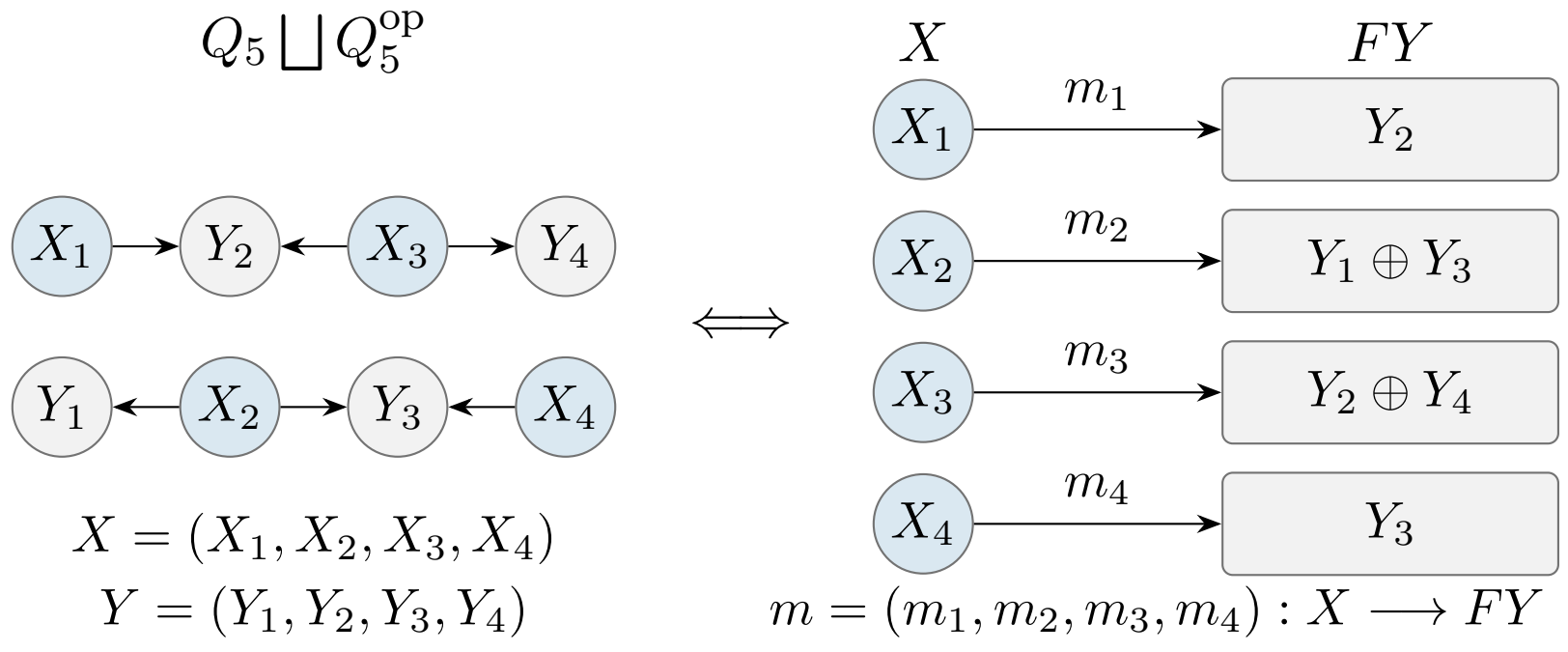

\begin{definition}
Fix $d\ge2$.  
Let $Q_d$ be the $A_{d-1}$ quiver on vertices $1,\ldots,d-1$,
with each edge oriented from its odd endpoint to its even endpoint, as in Figure~\ref{fig:Alternatibg-quiver}.
For $d=2$, the quiver consists of a single vertex.
Let $Q_d^{\rm op}$ be the quiver obtained by reversing every arrow, and set
\[
 \mathcal Q_d=\Db(\operatorname{rep}_{\C}Q_d),\qquad
 \widehat{\mathcal Q}_d
 =\Db\bigl(\operatorname{rep}_{\C}(Q_d\bigsqcup Q_d^{\rm op})\bigr).
\]
In $Q_d$, write $x_i$ for vertex $i$ when $i$ is odd and $y_i$
when $i$ is even; use the opposite labels in $Q_d^{\rm op}$.
Thus the arrows
are precisely $x_i\to y_j$ for $|i-j|=1$.
For $d=2$, we designate the first isolated vertex $x_1$ as a source
and the second $y_1$ as a target.

Let $\mathcal B=\mathcal C^{\oplus(d-1)}$.  
The \emph{adjacency functor} $F:\mathcal B\to\mathcal B$ collects
the neighboring entries of a tuple:
\begin{equation}\label{eq:adjacency-functor}
 (FB)_i=B_{i-1}\oplus B_{i+1},\qquad \hbox{\rm where}\quad B_0=B_d=0.
\end{equation}
The functor $F$ is exact and self-adjoint:
$
 \Hom_{\mathcal B}(FB,B')\simeq\Hom_{\mathcal B}(B,FB')$.

For a quiver diagram, write
$X_i$ for the object at the vertex labeled $x_i$ and $Y_i$ for the object
at the vertex labeled $y_i$.  These form two tuples $X=(X_i)$ and $Y=(Y_i)$
in $\mathcal B$.
The arrows leaving $x_i$ end at its neighboring vertices
$y_{i-1}$ and $y_{i+1}$.  Their maps therefore combine into a
single map $m_i:X_i\to(FY)_i$.  In this notation, the entire
quiver diagram is written as
\[
 E=(X\xrightarrow{m}FY),\qquad
 m_i:X_i\longrightarrow Y_{i-1}\oplus Y_{i+1}.
\]
Figure~\ref{fig:quiver-tuples} illustrates this description for $d=5$,
showing how the objects on the two paths are collected into tuples
and how their arrows give the components of $m$.
We interpret these maps using dg modules and chain maps,
so that their cones are functorial.
\end{definition}

\begin{definition}\label{def:reflection-path-reversal}
For a quiver diagram $E=(X\xrightarrow{m}FY)$, put
\[
 C_i=\Cone\bigl(X_i\xrightarrow{m_i}Y_{i-1}\oplus Y_{i+1}\bigr),
 \qquad C=(C_i)=\Cone(m).
\]
The \emph{reflection functor} $R$ on
$\mathcal C\otimes\widehat{\mathcal Q}_d$ applies source reflection
at every vertex $x_i$.  These reflections can be performed simultaneously,
since no two sources are joined by an arrow.
We then relabel the old target $y_i$ as the new source $x_i$, and
place the cone $C_i$ at the new target $y_i$.
Thus the new source tuple is $Y$, the new target tuple is $C$, and
$$
 RE=\bigl(Y\xrightarrow{m'}FC\bigr).$$
For $|i-j|=1$, the map along the new arrow $x_i\to y_j$ is
$
 Y_i\longrightarrow Y_{j-1}\oplus Y_{j+1}\longrightarrow C_j$,
where the first map is the inclusion of a summand and the second
is the canonical map to the cone.  With this relabeling, $R$ exchanges the
two quiver components for $d\ge3$.

The \emph{path-reversal functor} $\Pi$ reverses the vertex indices
while preserving the source and target labels.  Write
$(JB)_i=B_{d-i}$ for reversal of a tuple.  Then
$
 \Pi E=\bigl(JX\xrightarrow{Jm}JFY\simeq FJY\bigr)$.
The identification $JFY\simeq FJY$ exchanges the two neighboring
summands.  Thus $\Pi E$ has $X_{d-i}$ at $x_i$ and $Y_{d-i}$ at
$y_i$, and its map along $x_i\to y_j$ is the original map
$X_{d-i}\to Y_{d-j}$.  On morphisms, $\Pi$ applies the same
permutation of indices.  It exchanges the two components of the quiver when
$d$ is odd and preserves them when $d$ is even.

For $d=2$, the adjacency functor is zero, so these formulas become
$R(X,Y)=(Y,X[1])$ and $\Pi=\id$.
Taking $\mathcal C=\Perf(\C)$ gives the corresponding functors on
$\widehat{\mathcal Q}_d$ itself.
\end{definition}

Figure~\ref{fig:doubled-path} illustrates the two functors for $d=5$.

\begin{figure}[htbp]
\centering
\begin{tikzpicture}[scale=1.3,x=.88cm,y=.72cm,>=Stealth,
  vertex/.style={circle,draw=black!55,minimum size=5.8mm,
                 inner sep=1pt,font=\small},
  source/.style={vertex,fill=MidnightBlue!10},
  target/.style={vertex,fill=black!5},
  every edge/.style={draw,->}]
\begin{scope}
 \node at (1.5,.68) {$E=(X,Y)$};
 \node[source] (a1) at (0,0) {$X_1$};
 \node[target] (a2) at (1,0) {$Y_2$};
 \node[source] (a3) at (2,0) {$X_3$};
 \node[target] (a4) at (3,0) {$Y_4$};
 \node[target] (a5) at (0,-1) {$Y_1$};
 \node[source] (a6) at (1,-1) {$X_2$};
 \node[target] (a7) at (2,-1) {$Y_3$};
 \node[source] (a8) at (3,-1) {$X_4$};
 \path (a1) edge (a2) (a3) edge (a2) (a3) edge (a4)
       (a6) edge (a5) (a6) edge (a7) (a8) edge (a7);
\end{scope}
\begin{scope}[xshift=3.65cm]
 \node at (1.5,.68) {$RE=(Y,C)$};
 \node[source] (b1) at (0,0) {$Y_1$};
 \node[target] (b2) at (1,0) {$C_2$};
 \node[source] (b3) at (2,0) {$Y_3$};
 \node[target] (b4) at (3,0) {$C_4$};
 \node[target] (b5) at (0,-1) {$C_1$};
 \node[source] (b6) at (1,-1) {$Y_2$};
 \node[target] (b7) at (2,-1) {$C_3$};
 \node[source] (b8) at (3,-1) {$Y_4$};
 \path (b1) edge (b2) (b3) edge (b2) (b3) edge (b4)
       (b6) edge (b5) (b6) edge (b7) (b8) edge (b7);
\end{scope}
\begin{scope}[xshift=7.3cm]
 \node at (1.5,.68) {$\Pi E$};
 \node[source] (c1) at (0,0) {$X_4$};
 \node[target] (c2) at (1,0) {$Y_3$};
 \node[source] (c3) at (2,0) {$X_2$};
 \node[target] (c4) at (3,0) {$Y_1$};
 \node[target] (c5) at (0,-1) {$Y_4$};
 \node[source] (c6) at (1,-1) {$X_3$};
 \node[target] (c7) at (2,-1) {$Y_2$};
 \node[source] (c8) at (3,-1) {$X_1$};
 \path (c1) edge (c2) (c3) edge (c2) (c3) edge (c4)
       (c6) edge (c5) (c6) edge (c7) (c8) edge (c7);
\end{scope}
\end{tikzpicture}
\[
 C_i=\Cone\bigl(X_i\longrightarrow Y_{i-1}\oplus Y_{i+1}\bigr).
\]
\caption{Quivers for $d=5$. 
Both functors $R$ and $\Pi$ exchange the two components.}
\label{fig:doubled-path}
\end{figure}

\begin{lemma}\label{lem:all-degree-reflection-periodicity}
The functor $R$ is an enhanced autoequivalence, and $\Pi$ is an
enhanced involution commuting with $R$.
There are natural isomorphisms of enhanced functors
\[
 R^2\simeq S_{\widehat{\mathcal Q}_d}^{-1}[1],\qquad
 R^d\simeq\Pi[1],
\]
on $\widehat{\mathcal Q}_d$ and on
$\mathcal C\otimes\widehat{\mathcal Q}_d$, where the Serre functor
acts on the quiver factor only.
\end{lemma}

\begin{figure}[htbp]
\centering
\begin{tikzpicture}[x=1cm,y=1cm,>=Stealth,
  vertex/.style={circle,draw=black!55,minimum size=7.2mm,
                 inner sep=1pt,font=\small},
  source/.style={vertex,fill=MidnightBlue!10},
  target/.style={vertex,fill=black!5},
  every edge/.style={draw,->}]
\begin{scope}[shift={(0,0)}]
 \node[font=\small] at (1.725,.85) {$E=(X,Y)$};
 \node[source] (a0s1) at (0,0) {$X_1$};
 \node[target] (a0t2) at (1.15,0) {$Y_2$};
 \node[source] (a0s3) at (2.3,0) {$X_3$};
 \node[target] (a0t4) at (3.45,0) {$Y_4$};
 \node[target] (a0t1) at (0,-1.1) {$Y_1$};
 \node[source] (a0s2) at (1.15,-1.1) {$X_2$};
 \node[target] (a0t3) at (2.3,-1.1) {$Y_3$};
 \node[source] (a0s4) at (3.45,-1.1) {$X_4$};
 \path (a0s1) edge (a0t2) (a0s3) edge (a0t2)
       (a0s3) edge (a0t4) (a0s2) edge (a0t1)
       (a0s2) edge (a0t3) (a0s4) edge (a0t3);
\end{scope}
\begin{scope}[shift={(5.2,0)}]
 \node[font=\small] at (1.725,.85) {$RE=(Y,C)$};
 \node[source] (a1s1) at (0,0) {$Y_1$};
 \node[target] (a1t2) at (1.15,0) {$C_2$};
 \node[source] (a1s3) at (2.3,0) {$Y_3$};
 \node[target] (a1t4) at (3.45,0) {$C_4$};
 \node[target] (a1t1) at (0,-1.1) {$C_1$};
 \node[source] (a1s2) at (1.15,-1.1) {$Y_2$};
 \node[target] (a1t3) at (2.3,-1.1) {$C_3$};
 \node[source] (a1s4) at (3.45,-1.1) {$Y_4$};
 \path (a1s1) edge (a1t2) (a1s3) edge (a1t2)
       (a1s3) edge (a1t4) (a1s2) edge (a1t1)
       (a1s2) edge (a1t3) (a1s4) edge (a1t3);
\end{scope}
\begin{scope}[shift={(10.4,0)}]
 \node[font=\small] at (1.725,.85) {$R^2E=(C,D)$};
 \node[source] (a2s1) at (0,0) {$C_1$};
 \node[target] (a2t2) at (1.15,0) {$D_2$};
 \node[source] (a2s3) at (2.3,0) {$C_3$};
 \node[target] (a2t4) at (3.45,0) {$D_4$};
 \node[target] (a2t1) at (0,-1.1) {$D_1$};
 \node[source] (a2s2) at (1.15,-1.1) {$C_2$};
 \node[target] (a2t3) at (2.3,-1.1) {$D_3$};
 \node[source] (a2s4) at (3.45,-1.1) {$C_4$};
 \path (a2s1) edge (a2t2) (a2s3) edge (a2t2)
       (a2s3) edge (a2t4) (a2s2) edge (a2t1)
       (a2s2) edge (a2t3) (a2s4) edge (a2t3);
\end{scope}
\draw[->,thick,MidnightBlue] (4,-0.55)--(4.65,-0.55)
  node[midway,above=3pt] {$R$};
\draw[->,thick,MidnightBlue] (9.2,-0.55)--(9.85,-0.55)
  node[midway,above=3pt] {$R$};

\end{tikzpicture}
\caption{$R^2\simeq S^{-1}[1]$.}
\label{fig:inverse-coxeter}
\end{figure}
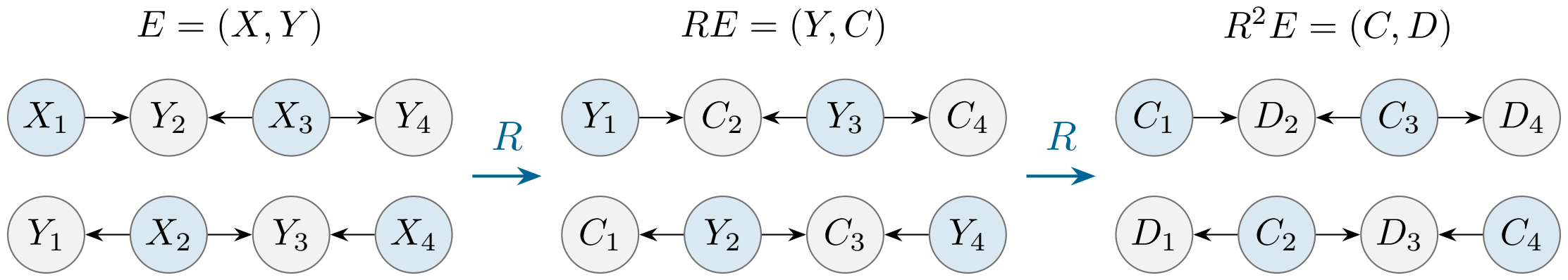

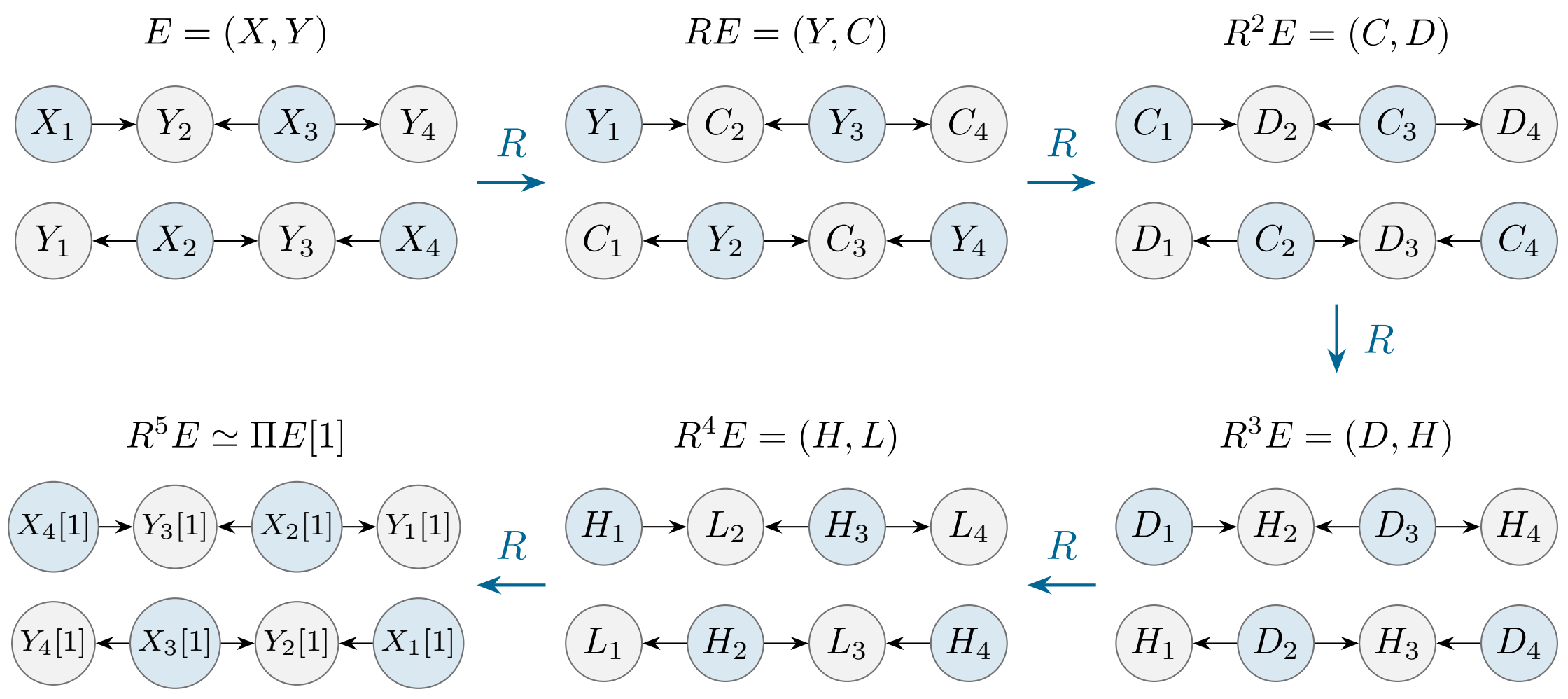
\begin{figure}[htbp]
\centering
\begin{tikzpicture}[x=1cm,y=1cm,>=Stealth,
  vertex/.style={circle,draw=black!55,minimum size=7.2mm,
                 inner sep=1pt,font=\small},
  source/.style={vertex,fill=MidnightBlue!10},
  target/.style={vertex,fill=black!5},
  every edge/.style={draw,->}]
\begin{scope}[shift={(0,-4)}]
 \node[font=\small] at (1.725,.85) {$E=(X,Y)$};
 \node[source] (b0s1) at (0,0) {$X_1$};
 \node[target] (b0t2) at (1.15,0) {$Y_2$};
 \node[source] (b0s3) at (2.3,0) {$X_3$};
 \node[target] (b0t4) at (3.45,0) {$Y_4$};
 \node[target] (b0t1) at (0,-1.1) {$Y_1$};
 \node[source] (b0s2) at (1.15,-1.1) {$X_2$};
 \node[target] (b0t3) at (2.3,-1.1) {$Y_3$};
 \node[source] (b0s4) at (3.45,-1.1) {$X_4$};
 \path (b0s1) edge (b0t2) (b0s3) edge (b0t2)
       (b0s3) edge (b0t4) (b0s2) edge (b0t1)
       (b0s2) edge (b0t3) (b0s4) edge (b0t3);
\end{scope}
\begin{scope}[shift={(5.2,-4)}]
 \node[font=\small] at (1.725,.85) {$RE=(Y,C)$};
 \node[source] (b1s1) at (0,0) {$Y_1$};
 \node[target] (b1t2) at (1.15,0) {$C_2$};
 \node[source] (b1s3) at (2.3,0) {$Y_3$};
 \node[target] (b1t4) at (3.45,0) {$C_4$};
 \node[target] (b1t1) at (0,-1.1) {$C_1$};
 \node[source] (b1s2) at (1.15,-1.1) {$Y_2$};
 \node[target] (b1t3) at (2.3,-1.1) {$C_3$};
 \node[source] (b1s4) at (3.45,-1.1) {$Y_4$};
 \path (b1s1) edge (b1t2) (b1s3) edge (b1t2)
       (b1s3) edge (b1t4) (b1s2) edge (b1t1)
       (b1s2) edge (b1t3) (b1s4) edge (b1t3);
\end{scope}
\begin{scope}[shift={(10.4,-4)}]
 \node[font=\small] at (1.725,.85) {$R^2E=(C,D)$};
 \node[source] (b2s1) at (0,0) {$C_1$};
 \node[target] (b2t2) at (1.15,0) {$D_2$};
 \node[source] (b2s3) at (2.3,0) {$C_3$};
 \node[target] (b2t4) at (3.45,0) {$D_4$};
 \node[target] (b2t1) at (0,-1.1) {$D_1$};
 \node[source] (b2s2) at (1.15,-1.1) {$C_2$};
 \node[target] (b2t3) at (2.3,-1.1) {$D_3$};
 \node[source] (b2s4) at (3.45,-1.1) {$C_4$};
 \path (b2s1) edge (b2t2) (b2s3) edge (b2t2)
       (b2s3) edge (b2t4) (b2s2) edge (b2t1)
       (b2s2) edge (b2t3) (b2s4) edge (b2t3);
\end{scope}
\draw[->,thick,MidnightBlue] (4,-4.55)--(4.65,-4.55)
  node[midway,above=3pt] {$R$};
\draw[->,thick,MidnightBlue] (9.2,-4.55)--(9.85,-4.55)
  node[midway,above=3pt] {$R$};
\draw[->,thick,MidnightBlue] (12.125,-5.7)--(12.125,-6.35)
  node[midway,right=3pt] {$R$};
\begin{scope}[shift={(10.4,-7.8)}]
 \node[font=\small] at (1.725,.85) {$R^3E=(D,H)$};
 \node[source] (b3s1) at (0,0) {$D_1$};
 \node[target] (b3t2) at (1.15,0) {$H_2$};
 \node[source] (b3s3) at (2.3,0) {$D_3$};
 \node[target] (b3t4) at (3.45,0) {$H_4$};
 \node[target] (b3t1) at (0,-1.1) {$H_1$};
 \node[source] (b3s2) at (1.15,-1.1) {$D_2$};
 \node[target] (b3t3) at (2.3,-1.1) {$H_3$};
 \node[source] (b3s4) at (3.45,-1.1) {$D_4$};
 \path (b3s1) edge (b3t2) (b3s3) edge (b3t2)
       (b3s3) edge (b3t4) (b3s2) edge (b3t1)
       (b3s2) edge (b3t3) (b3s4) edge (b3t3);
\end{scope}
\begin{scope}[shift={(5.2,-7.8)}]
 \node[font=\small] at (1.725,.85) {$R^4E=(H,L)$};
 \node[source] (b4s1) at (0,0) {$H_1$};
 \node[target] (b4t2) at (1.15,0) {$L_2$};
 \node[source] (b4s3) at (2.3,0) {$H_3$};
 \node[target] (b4t4) at (3.45,0) {$L_4$};
 \node[target] (b4t1) at (0,-1.1) {$L_1$};
 \node[source] (b4s2) at (1.15,-1.1) {$H_2$};
 \node[target] (b4t3) at (2.3,-1.1) {$L_3$};
 \node[source] (b4s4) at (3.45,-1.1) {$H_4$};
 \path (b4s1) edge (b4t2) (b4s3) edge (b4t2)
       (b4s3) edge (b4t4) (b4s2) edge (b4t1)
       (b4s2) edge (b4t3) (b4s4) edge (b4t3);
\end{scope}
\begin{scope}[shift={(0,-7.8)}]
 \node[font=\small] at (1.725,.85) {$R^5E\simeq\Pi E[1]$};
 \node[source,font=\scriptsize] (b5s1) at (0,0) {$X_4[1]$};
 \node[target,font=\scriptsize] (b5t2) at (1.15,0) {$Y_3[1]$};
 \node[source,font=\scriptsize] (b5s3) at (2.3,0) {$X_2[1]$};
 \node[target,font=\scriptsize] (b5t4) at (3.45,0) {$Y_1[1]$};
 \node[target,font=\scriptsize] (b5t1) at (0,-1.1) {$Y_4[1]$};
 \node[source,font=\scriptsize] (b5s2) at (1.15,-1.1) {$X_3[1]$};
 \node[target,font=\scriptsize] (b5t3) at (2.3,-1.1) {$Y_2[1]$};
 \node[source,font=\scriptsize] (b5s4) at (3.45,-1.1) {$X_1[1]$};
 \path (b5s1) edge (b5t2) (b5s3) edge (b5t2)
       (b5s3) edge (b5t4) (b5s2) edge (b5t1)
       (b5s2) edge (b5t3) (b5s4) edge (b5t3);
\end{scope}
\draw[->,thick,MidnightBlue] (9.85,-8.35)--(9.2,-8.35)
  node[midway,above=3pt] {$R$};
\draw[->,thick,MidnightBlue] (4.65,-8.35)--(4,-8.35)
  node[midway,above=3pt] {$R$};
\end{tikzpicture}
\caption{$R^dE\simeq\Pi E[1]$.}
\label{fig:coxeter-periodicity}
\end{figure}

\begin{proof}
Reversing the path twice fixes every vertex,
and reversal commutes with the cone construction.  Thus
$\Pi^2\simeq\id$ and $\Pi R\simeq R\Pi$.
Each source reflection is an enhanced equivalence
\cite[Theorem~9.19 and Definition~9.20]{GrothStovicek}; see also
\cite[Sections~1--2]{DJWReflection}.
The inverse sink reflection recovers the source object from its
cone triangle, as explained in the proof of
Proposition~\ref{prop:one-variable-model}.  Hence $R$ is an enhanced
autoequivalence.  

The \emph{Coxeter functor} is the composition of one derived sink
reflection at each vertex, ordered so that each vertex is a sink
when it is reflected.  For a path quiver $Q$, this composition is
naturally isomorphic to $S_Q[-1]$
\cite[Definitions~5.8 and~5.13, Theorem~10.9]{GrothStovicekA}.
Its inverse uses source reflections.  For $d\ge3$, two applications
of $R$ perform exactly this inverse composition on each component:
first they reflect at the original sources, then at the original
targets, which have become sources.  Hence $R^2\simeq S^{-1}[1]$.
Figure~\ref{fig:inverse-coxeter} illustrates these two rounds
for $d=5$.

Proposition~\ref{prop:one-variable-model} also gives
$S^d\simeq[d-2]$, so the first identity implies
\[
 R^{2d}\simeq(S^{-1}[1])^d
          \simeq S^{-d}[d]\simeq[2].
\]
Thus $R^d[-1]$ is an involution.  To prove the second identity,
we must identify this involution with the path reversal $\Pi$.
We keep track of the vertex labels
under successive reflections.  For $d=3$, each component is a single arrow.
Starting with $A\to B$ and writing $C=\Cone(A\to B)$, the three
reflections follow its distinguished triangle:
\[
 (A\to B)\longmapsto(B\to C)\longmapsto(C\to A[1])
              \longmapsto(A[1]\to B[1]),
\]
up to natural isomorphisms of diagrams.  The final arrow lies in
the opposite component, with the vertex indices reversed, as
prescribed by $\Pi$.
For a general path, Groth--\v{S}\v{t}ov\'i\v{c}ek organize
successive reflections into a diagram of cones
\cite[Theorem~4.14 and Proposition~5.4]{GrothStovicekA}.
Their Corollary~4.11 identifies the shift $[1]$ with a flip of this
diagram.  For our alternating paths with $d-1$ vertices, $d$ rounds
of source reflection reach this flipped diagram with the vertex
order reversed.  Thus the objects at $x_i$ and $y_i$ become
$X_{d-i}[1]$ and $Y_{d-i}[1]$, with the corresponding shifted
arrow maps.  This gives $R^d\simeq\Pi[1]$.
Figure~\ref{fig:coxeter-periodicity} illustrates this identity
for $d=5$.
For $d=2$, both identities follow directly from
$R(X,Y)=(Y,X[1])$ and $\Pi=\id$.

These identities are realized by complexes of bimodules
\cite[Section~10]{GrothStovicekA}.
Tensoring them with the identity bimodule of $\mathcal C$
gives the same  identities for quiver diagrams with
coefficients in $\mathcal C$.
\end{proof}

\begin{theorem}
\label{thm:all-degree-coefficients}
Let $\mathcal C$ be an idempotent-complete smooth proper enhanced
triangulated category with a stability condition
$\sigma=(Z,\mathcal P)$ supported on a finite-rank numerical quotient
$v:K_0(\mathcal C)\to\Lambda$.  

The category $\mathcal C\otimes\widehat{\mathcal Q}_d$ has a
numerical stability condition $\Sigma=(\mathcal Z,\mathcal P_\Sigma)$
supported on $\Lambda^{\oplus 2(d-1)}$ and of Gepner type
$(\Phi,1/d)$ for each of the following autoequivalences $\Phi$:
\begin{enumerate}
\item[(i)] the reflection functor $R$;
\item[(ii)] the vertexwise action of every exact enhanced
  autoequivalence $\Psi:\mathcal C\to\mathcal C$ for which
  $\sigma$ is of Gepner type $(\Psi,1/d)$.
\end{enumerate}

Moreover, for every quiver diagram $E$,
\[
 \sum_w m_\sigma(E_w)\le B_d m_\Sigma(E),\qquad
 B_d=\bigl(1+2\cos(\pi/d)\bigr)^d\sin(\pi/(2d))^{-4},
\]
where the sum runs over all vertices.
If $\sigma$ has a total mass--Hom bound, so does $\Sigma$.
\end{theorem}

The proof occupies the rest of the section.  Fix
$\sigma=(Z,\mathcal P)$ as in the theorem.  We first choose the
central charge $\mathcal Z$ and tilt the standard heart of $\mathcal C\otimes\widehat{\mathcal Q}_d$ to place its
charges in the upper half-plane.  Write $U$ for the tilted aisle.
The phase bound in~\eqref{eq:all-degree-torsion-mass} gives
$RU\subset U$, so the aisles $R^jU$, $j\in\Z$, form a nested
sequence.  Successive truncations at these aisles give factors in
$R^jU\cap(R^{j+1}U)^\perp$, whose nonzero charges lie in sectors
of angle $\pi/d$.
The mass bound in~\eqref{eq:all-degree-torsion-mass}, together
with the reflection estimate~\eqref{eq:all-degree-reflection-mass},
yields~\eqref{eq:all-degree-block-mass}: the weighted sum of the
vertex masses of each factor is bounded by a constant times the
absolute value of its central charge.
Finally, we show that Gepner type
$(\Psi,1/d)$ is preserved under the vertexwise action of
an autoequivalence $\Psi:\mathcal C\to\mathcal C$ for which
  $\sigma$ is of Gepner type $(\Psi,1/d)$.

\subsection{The tilted heart}

Let $\mathcal A=\mathcal P((0,1])$ be the heart on $\mathcal C$.
Give the tuple category $\mathcal B=\mathcal C^{\oplus(d-1)}$
the product slicing $\mathcal P_{\mathcal B}$:
a tuple is semistable of phase $\phi$ when each nonzero entry belongs
to $\mathcal P(\phi)$.  Its heart is
$\mathcal A_{\mathcal B}=\mathcal A^{\oplus(d-1)}$.
The adjacency functor $F$ preserves each slice: if every
$B_j$ belongs to $\mathcal P(\phi)$, then so does
$(FB)_j=B_{j-1}\oplus B_{j+1}$, because $\mathcal P(\phi)$
is closed under finite direct sums.  Here, as before, $B_0=B_d=0$.
Define the positive weights and weighted charge by
\[
 r_j=\frac{\sin(j\pi/d)}{\sin(\pi/d)},\qquad
 z(B)=\sum_{j=1}^{d-1}r_jZ(B_j).
\]
The identities 
$$r_{j-1}+r_{j+1}=\lambda r_j,\qquad
 \lambda=2\cos(\pi/d),\quad\hbox{\rm with}\quad r_0=r_d=0,$$
give $z(FB)=\lambda z(B)$.
We will also use the symmetry $r_j=r_{d-j}$.

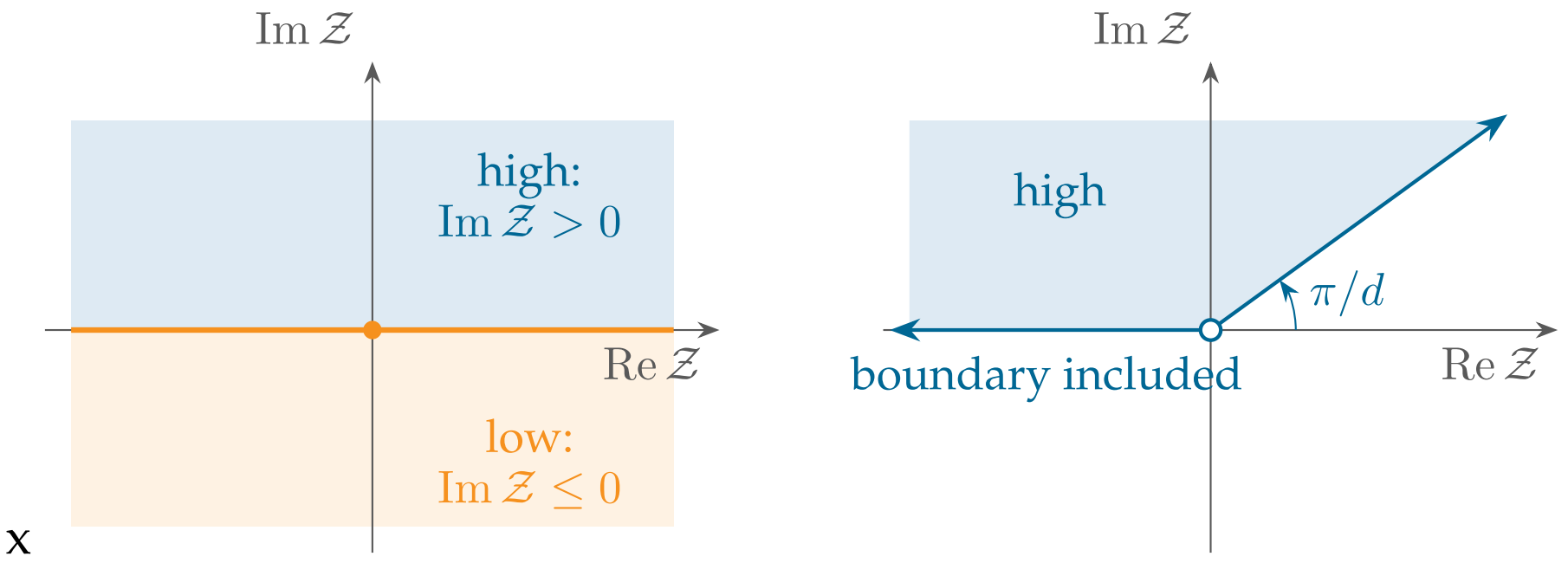
\begin{figure}[htbp]x
\centering
\begin{tikzpicture}[>=Stealth,font=\small]
 \begin{scope}
  \fill[MidnightBlue!9] (-2.3,0) rectangle (2.3,1.6);
  \fill[BurntOrange!10] (-2.3,-1.5) rectangle (2.3,0);
  \draw[->,black!65] (-2.5,0)--(2.65,0)
       node[below left] {$\operatorname{Re}\mathcal Z$};
  \draw[->,black!65] (0,-1.7)--(0,2.05)
       node[above left] {$\operatorname{Im}\mathcal Z$};
  \draw[BurntOrange,very thick] (-2.3,0)--(2.3,0);
  \fill[BurntOrange] (0,0) circle (2pt);
  \node[MidnightBlue,align=center] at (1.2,1.03)
       {high:\\[-1pt] $\operatorname{Im}\mathcal Z>0$};
  \node[BurntOrange,align=center] at (1.2,-1.02)
       {low:\\[-1pt] $\operatorname{Im}\mathcal Z\le0$};
 \end{scope}
 \begin{scope}[xshift=6.4cm]
  \fill[MidnightBlue!9]
       (-2.3,0)--(-2.3,1.6)--({1.6/tan(36)},1.6)--(0,0)--cycle;
  \draw[->,black!65] (-2.5,0)--(2.65,0)
       node[below left] {$\operatorname{Re}\mathcal Z$};
  \draw[->,black!65] (0,-1.7)--(0,2.05)
       node[above left] {$\operatorname{Im}\mathcal Z$};
  \draw[->,MidnightBlue,thick] (0,0)--(-2.45,0);
  \draw[->,MidnightBlue,thick] (0,0)--(36:2.8);
  \draw[->,MidnightBlue] (.65,0) arc[start angle=0,end angle=36,radius=.65];
  \node[MidnightBlue] at (1.04,.28) {$\pi/d$};
  \node[MidnightBlue] at (-1.15,1.03) {high};
  \node[MidnightBlue,below] at (-1.25,-.06) {boundary included};
  \draw[MidnightBlue,fill=white,thick] (0,0) circle (2.2pt);
 \end{scope}
\end{tikzpicture}
\caption{Charges of high/low objects.  Left: $p>0$. Right: $p=0$.}
\label{fig:high-low-charges}
\end{figure}

By Lemma~\ref{lem:quiver-coefficient-heart},
$\mathcal H=\operatorname{Rep}(Q_d\bigsqcup Q_d^{\rm op},\mathcal A)$ is the heart
of a bounded t-structure on $\mathcal C\otimes\widehat{\mathcal Q}_d$.
Its kernels, cokernels, and exact sequences are computed at each
of the $2(d-1)$ vertices.
Put
\[
 t=1-\frac1d,\qquad q=e^{-\pi it}=-e^{\pi i/d}.
\]
Then $q^2+\lambda q+1=0$.
Define the central charge $\mathcal Z$ and the real-valued
additive function $p$ by
\[
 \mathcal Z(X,Y)=z(X)+qz(Y),\qquad
 p(X,Y)=\operatorname{Im}z(X)+\operatorname{Im}z(Y).
\]
The torsion pair will be determined by $\operatorname{Im}\mathcal Z$
and $p$.
Multiplication by $q$ rotates the charges of the $Y$-vertices
clockwise by $\pi t$, so $\mathcal Z$ need not take $\mathcal H$
into the upper half-plane.  We correct this by tilting the standard
heart.
Following the terminology of~\cite[Section~3]{OAIQuintic},
we call a nonzero object of $\mathcal H$ \emph{high} if
$p=0$ or $\operatorname{Im}\mathcal Z>0$, and \emph{low} if $p>0$ and $\operatorname{Im}\mathcal Z\le0$.
Since $p\ge0$ on $\mathcal H$, these conditions partition its nonzero
objects.  The boundary case $p=0$ deserves attention: all nonzero
vertex objects then have phase one, so
\[
 z(X)=-a,\qquad z(Y)=-b,\qquad
 \mathcal Z(X,Y)=-a+b e^{\pi i/d},
 \qquad a,b\ge0,
\]
with $a,b$ not both zero.  In particular, $\operatorname{Im}\mathcal Z\ge0$ in this case.
Figure~\ref{fig:high-low-charges} shows both cases separately.

\begin{proposition}\label{prop:all-degree-torsion}
There is a torsion pair $(\mathcal T,\mathcal F)$ in $\mathcal H$, where
\[
 \begin{split}
 \mathcal T&=\{E:\text{every nonzero quotient of $E$ in $\mathcal H$
                              is high}\},\\
 \mathcal F&=\{E:\text{every nonzero subobject of $E$ in $\mathcal H$
                              is low}\}.
 \end{split}
\]
Neither rationality of the charge $Z$ nor noetherianity of the heart
$\mathcal A$ is assumed.
\end{proposition}

\begin{proof}
For a tuple $Y\in\mathcal A_{\mathcal B}$, let $Y_{>t}\subseteq Y$
denote the subobject determined by all HN factors of phase greater
than $t$, and set $Y_{\le t}=Y/Y_{>t}$.
Every object $E=(X,Y,m)\in\mathcal F$ satisfies
\begin{equation}\label{eq:all-degree-torsion-constraints}
 Y\in\mathcal P_{\mathcal B}((0,t]),\qquad m:X\hookrightarrow FY.
\end{equation}
Indeed, $(0,Y_{>t})$ and $(\ker m,0)$ are subobjects of $E$
and are high whenever nonzero, so both vanish.
Since $F$ preserves the product slicing and
$\mathcal P_{\mathcal B}((0,t])$ is closed under subobjects,
\eqref{eq:all-degree-torsion-constraints} puts all vertex phases in $(0,t]$.  Every nonzero subobject
of a diagram satisfying~\eqref{eq:all-degree-torsion-constraints} therefore has $p>0$.

Given an object $E=(X,Y,m)$ in $\mathcal H$, we denote by $\overline E$ its quotient
obtained by first taking the quotient by the subobject $(0,Y_{>t})$ and then by the subobject
$\bigl(\ker(X\to F(Y_{\le t})),0\bigr)$. Thus $\overline E$ satisfies \eqref{eq:all-degree-torsion-constraints}.
The first subobject has $\operatorname{Im}\mathcal Z>0$ unless it
vanishes, since rotation by $q$ places its phases in $(0,1/d]$;
the second satisfies $\operatorname{Im}\mathcal Z=p\ge0$. Thus
\[
\operatorname{Im}\mathcal Z(\overline E)\le \operatorname{Im}\mathcal Z(E),
\]
with equality only if both removed subobjects have $p=0$. In particular,
\[
\operatorname{Im}\mathcal Z(\overline E)=\operatorname{Im}\mathcal Z(E)\quad\Longrightarrow\quad
p(\overline E)=p(E).
\]
\smallskip
\noindent\emph{Finitely many possible charges.}
Fix $0\ne E\in\mathcal H$, and let $l>0$ be the smallest HN
phase among its nonzero vertex objects.  Every quotient $Q$ of
$E$ satisfying~\eqref{eq:all-degree-torsion-constraints} has vertex phases in $[l,t]$: quotients
cannot lower the least HN phase.  If $l>t$, only the zero quotient
occurs.  Otherwise set
\[
 \epsilon_E=\min_{\phi\in[l,t]}\sin(\pi\phi)>0.
\]
For a support constant $C_\sigma$ and every vertex $w$, the
support property and the phase bounds give
\[
 \|v(Q_w)\|\le C_\sigma m_\sigma(Q_w)
 \le C_\sigma\epsilon_E^{-1}\operatorname{Im}Z(Q_w)
 \le C_\sigma\epsilon_E^{-1}\operatorname{Im}Z(E_w).
\]
The last inequality follows from the nonnegative imaginary charge
of the kernel of $E_w\to Q_w$.  Thus the vertex classes lie in a
bounded subset of the lattice $\Lambda^{\oplus2(d-1)}$.
Only finitely many classes occur, so $\operatorname{Im}\mathcal Z(Q)$ and $p(Q)$ take only
finitely many values.

\smallskip
\noindent\emph{The torsion sequence.}
Among the quotients $E\twoheadrightarrow Q$ satisfying~\eqref{eq:all-degree-torsion-constraints},
including zero, choose one minimizing $\operatorname{Im}\mathcal Z(Q)$ and, subject to this,
maximizing $p(Q)$.  Put $K=\ker(E\to Q)$.

First, $Q\in\mathcal F$.  Every nonzero subobject $A\subset Q$
has $p(A)>0$.  If $\operatorname{Im}\mathcal Z(A)>0$, the quotient $\overline{Q/A}$ of $E$
satisfies~\eqref{eq:all-degree-torsion-constraints} and
\[
 \operatorname{Im}\mathcal Z(\overline{Q/A})\le \operatorname{Im}\mathcal Z(Q/A)
 =\operatorname{Im}\mathcal Z(Q)-\operatorname{Im}\mathcal Z(A)<\operatorname{Im}\mathcal Z(Q),
\]
contradicting minimality.  Thus every such $A$ is low.

Next, suppose that $K\notin\mathcal T$.  Choose a nonzero low
quotient $K\twoheadrightarrow L$, put $K'=\ker(K\to L)$, and
write $P=E/K'$.  Then
\[
 0\longrightarrow L\longrightarrow P\longrightarrow Q
   \longrightarrow0.
\]
The quotient $\overline P$ of $E$ satisfies~\eqref{eq:all-degree-torsion-constraints}, so minimality gives
\[
 \operatorname{Im}\mathcal Z(Q)\le \operatorname{Im}\mathcal Z(\overline P)\le \operatorname{Im}\mathcal Z(P)\\
 =\operatorname{Im}\mathcal Z(Q)+\operatorname{Im}\mathcal Z(L)\le \operatorname{Im}\mathcal Z(Q).
\]
Equality holds throughout.  The construction of $\overline P$
therefore preserves $p$, giving
\[
 p(\overline P)=p(P)=p(Q)+p(L)>p(Q),
\]
contrary to the choice of $Q$.  Hence $K\in\mathcal T$.

Finally, a nonzero image of a morphism from $\mathcal T$ to
$\mathcal F$ would be both high and low.  Thus
$\Hom(\mathcal T,\mathcal F)=0$, and
$0\to K\to E\to Q\to0$ is the required torsion sequence.
\end{proof}

We now tilt the heart $\mathcal H$ at the torsion pair
$(\mathcal T,\mathcal F)$, as recalled in Section~\ref{sec:background}.
The tilted heart $\mathcal A'=\langle\mathcal F[1],\mathcal T\rangle$
consists precisely of objects $E$ fitting into a distinguished triangle
\[
 F_E[1]\longrightarrow E\longrightarrow T_E\longrightarrow F_E[2],
 \qquad F_E\in\mathcal F,\quad T_E\in\mathcal T.
\]
Thus the tilt keeps objects of $\mathcal T$ in degree zero and
moves objects of $\mathcal F$ to degree $-1$.  This shift reverses
their central charges: $\mathcal Z(F_E[1])=-\mathcal Z(F_E)$.
Let $V$ be the aisle of $\mathcal H$. The aisle of the tilted
heart is
$U=\{E\in V:H^0_{\mathcal H}(E)\in\mathcal T\}$ and we have
$V[1]\subset U\subset V$.

\subsection{A mass estimate}

The central charge $\mathcal Z(E)$ combines the charges of the source
vertices with the rotated charges of the target vertices, so cancellation
may reduce its magnitude.  For $E\in\mathcal F$, the condition on
subobjects controls this loss.
To make this precise, we measure the size of a diagram before
combining its charges: for any quiver diagram $E=(X,Y,m)$, put
\[
 M(E)=\sum_{j=1}^{d-1}r_j\bigl(m_\sigma(X_j)+m_\sigma(Y_j)\bigr).
\]
Evaluation at a vertex commutes with shifts and preserves
distinguished triangles.  Shifting an object of $\mathcal C$
shifts its HN factors without changing the absolute values of
their central charges.  Thus $M(E[n])=M(E)$ for every integer $n$.
Similarly, a distinguished triangle of quiver diagrams
\[
 E_1\longrightarrow E_2\longrightarrow E_3\longrightarrow E_1[1]
\]
gives a distinguished triangle at each vertex.  Applying the
mass triangle inequality~\eqref{eq:mass-triangle-inequality} at each vertex and summing
over all source and target vertices with the positive weights
$r_j$, we obtain
\begin{equation}\label{eq:all-degree-mass-triangle}
 M(E_2)\le M(E_1)+M(E_3).
\end{equation}

\begin{lemma}\label{lem:all-degree-mass}
Let $c=\cos(\pi t/2)=\sin(\pi/(2d))>0$.
For every nonzero $E\in\mathcal F$, using $\arg\in(-\pi,\pi]$,
\begin{equation}\label{eq:all-degree-torsion-mass}
 \arg \mathcal Z(E)/\pi\in(-t,0],\qquad
 |\mathcal Z(E)|\ge c^3M(E).
 \end{equation}
\end{lemma}

\begin{proof}
We apply the following elementary geometric observation three times.
If vectors $\zeta_\nu\in\C$ lie in a sector of angle at most
$\pi t<\pi$, write $\theta$ for the angle of its bisecting ray.
The orthogonal projection of each $\zeta_\nu$ onto
this ray has length at least $c|\zeta_\nu|$, where
$c=\cos(\pi t/2)$.  Thus
\begin{equation}\label{eq:sector-projection-estimate}
 \left|\sum_\nu\zeta_\nu\right|
 \ge\operatorname{Re}\left(e^{-i\theta}\sum_\nu\zeta_\nu\right)
 \ge c\sum_\nu|\zeta_\nu|,
\end{equation}
and the sum lies in the same sector.  In particular, it is
nonzero if at least one vector is nonzero.
We organize the proof into three claims, illustrated in
Figure~\ref{fig:phase-pairing}.

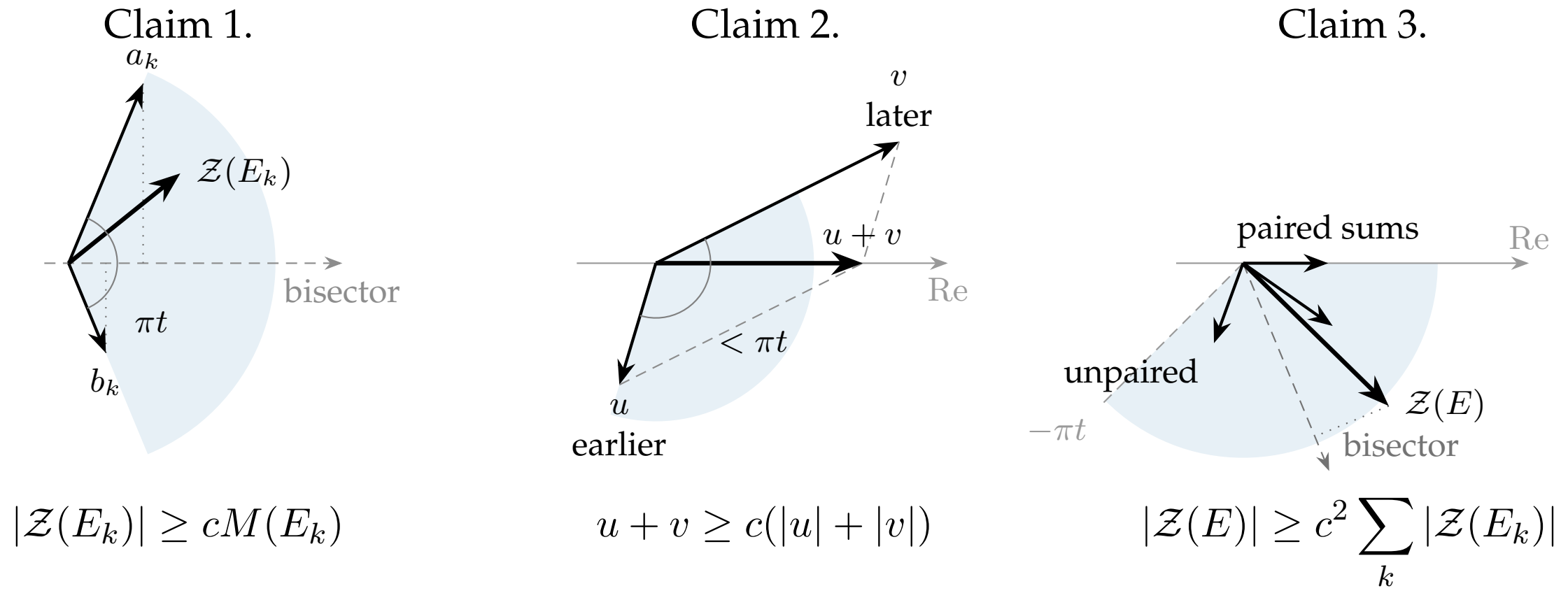
\begin{figure}[htbp]
\centering
\begin{minipage}[t]{.32\linewidth}
\centering
\textup{Claim 1.}\par\smallskip
\begin{tikzpicture}[>=Stealth,x=1.15cm,y=1.15cm,font=\footnotesize]
\path[use as bounding box] (-.85,-1.85) rectangle (2.6,1.75);
\coordinate (O) at (0,0);
\coordinate (A) at (67.5:1.6);
\coordinate (B) at (-67.5:.8);
\coordinate (W) at ($(A)+(B)$);
\fill[MidnightBlue!6] (O)--(-67.5:1.7)
 arc[start angle=-67.5,end angle=67.5,radius=1.7]--cycle;
\draw[densely dashed,black!45,->] (-.2,0)--(2.25,0)
 node[below] {bisector};
\draw[dotted,black!55] (A)--(A |- O) (B)--(B |- O);
\draw[->,thick,black] (O)--(A) node[above] {$a_k$};
\draw[->,thick,black] (O)--(B) node[below] {$b_k$};
\draw[->,very thick] (O)--(W) node[right] {$\mathcal Z(E_k)$};
\draw[black!50] (-67.5:.4) arc[start angle=-67.5,end angle=67.5,radius=.4];
\node at (.68,-.48) {$\pi t$};
\end{tikzpicture}\par
\(\displaystyle |\mathcal Z(E_k)|\ge cM(E_k)\)
\end{minipage}\hfill
\begin{minipage}[t]{.32\linewidth}
\centering
\textup{Claim 2.}\par\smallskip
\begin{tikzpicture}[>=Stealth,x=1.15cm,y=1.15cm,font=\footnotesize]
\path[use as bounding box] (-.85,-1.85) rectangle (2.6,1.75);
\coordinate (O) at (0,0);
\coordinate (U) at (-.3,-1);
\coordinate (V) at (2,1);
\coordinate (S) at (1.7,0);
\fill[MidnightBlue!6] (O)--(-106.699:1.3)
 arc[start angle=-106.699,end angle=26.565,radius=1.3]--cycle;
\draw[->,black!40] (-.65,0)--(2.4,0) node[below] {$\operatorname{Re}$};
\draw[->,thick,black] (O)--(U)
 node[below,align=center] {$u$\\earlier};
\draw[->,thick,black] (O)--(V)
 node[above,align=center] {$v$\\later};
\draw[densely dashed,black!45] (U)--(S)--(V);
\draw[->,very thick] (O)--(S) node[above] {$u+v$};
\draw[black!50] (-106.699:.45)
 arc[start angle=-106.699,end angle=26.565,radius=.45];
\node at (.8,-.65) {$<\pi t$};
\end{tikzpicture}\par
\(\displaystyle u+v\ge c(|u|+|v|)\)
\end{minipage}\hfill
\begin{minipage}[t]{.32\linewidth}
\centering
\textup{Claim 3.}\par\smallskip
\begin{tikzpicture}[>=Stealth,x=1.15cm,y=1.15cm,font=\footnotesize]
\path[use as bounding box] (-.85,-1.85) rectangle (2.6,1.75);
\coordinate (O) at (0,0);
\coordinate (A) at (.7,0);
\coordinate (B) at (-35:.9);
\coordinate (C) at (-110:.7);
\coordinate (W) at ($(A)+(B)+(C)$);
\coordinate (P) at (.590,-1.425);
\fill[MidnightBlue!6] (O)--(-135:1.6)
 arc[start angle=-135,end angle=0,radius=1.6]--cycle;
\draw[->,black!40] (-.55,0)--(2.35,0)
 node[above] {$\operatorname{Re}$};
\draw[densely dashed,black!40] (O)--(-135:1.65)
 node[below left] {$-\pi t$};
\draw[densely dashed,black!55,->] (O)--(-67.5:1.85)
 node[above right] {bisector};
\draw[->,thick,black] (O)--(A)
 node[above] {paired sums};
\draw[->,thick,black] (O)--(B);
\draw[->,thick,black] (O)--(C)
 node[below left] {unpaired};
\draw[dotted,black!60] (W)--(P);
\draw[->,very thick] (O)--(W) node[right] {$\mathcal Z(E)$};
\end{tikzpicture}\par
\(\displaystyle |\mathcal Z(E)|\ge c^2\sum_k|\mathcal Z(E_k)|\)
\end{minipage}
\caption{Three applications of the projection estimate.}
\label{fig:phase-pairing}
\end{figure}

\smallskip

\begin{samepage}
\begingroup\itshape
\noindent Claim 1.
There is a filtration in $\mathcal H$
\begin{equation}\label{qerfqrfqrgqrwgqwrgw}
 0=E^{(0)}\subset E^{(1)}\subset\cdots\subset E^{(s)}=E
\end{equation}
with nonzero factors $E_k=E^{(k)}/E^{(k-1)}$ such that all
nonzero vertex objects of $E_k$ are semistable of one phase
$\phi_k$, where $t\ge\phi_1>\cdots>\phi_s>0$.  Moreover,
\begin{equation}\label{asdfasfasfasgqwgqw}
 M(E)=\sum_k M(E_k),
 \end{equation}
 \begin{equation}\label{asvasfvasfasf}
\arg\mathcal Z(E_k)/\pi\in[\phi_k-t,\phi_k]\subset(-t,1),\qquad
 |\mathcal Z(E_k)|\ge cM(E_k).
 \end{equation}
\par\endgroup
\end{samepage}

\noindent\emph{Proof of Claim~1.}
Since $E\in\mathcal F$, the observation \eqref{eq:all-degree-torsion-constraints}
says that each target object $Y_j$ has HN phases in $(0,t]$
and that $m:X\to FY$ is injective.  Thus each source object
$X_j$ is a subobject of $Y_{j-1}\oplus Y_{j+1}$.
The latter also has HN phases in $(0,t]$, and
$\mathcal P((0,t])$ is closed under subobjects in the heart
$\mathcal A$.  Hence all vertex HN phases of $E$ lie in $(0,t]$.
Let $\phi_1>\cdots>\phi_s$ be the distinct phases occurring
among the HN factors of all vertex objects of $E$; thus
$0<\phi_k\le t$ for every $k$.
For each $k$, take at every vertex the initial HN subobject
whose factors have phase at least $\phi_k$.
These HN truncations are functorial~\cite[Section~3]{Bridgeland},
so every arrow map restricts to a map between the chosen subobjects.
They therefore form a quiver subdiagram $E^{(k)}\subset E$,
giving the  filtration \eqref{qerfqrfqrgqrwgqwrgw}.

At each vertex, the quotient $E_k=E^{(k)}/E^{(k-1)}$ is the
HN factor of phase $\phi_k$, or zero if that phase does not
occur there.  Thus all its nonzero vertex objects are
semistable of the same phase~$\phi_k$.
Every HN factor at every vertex occurs exactly once among
these quotients.  Since mass is the sum of the absolute values
of their central charges, summing with the weights $r_j$ gives \eqref{asdfasfasfasgqwgqw}.

Write $\mathcal Z(E_k)=a_k+b_k$, where $a_k$ is the sum of
the source charges and $b_k$ is the sum of the rotated target
charges, with the weights $r_j$ included.  If nonzero, $a_k$
has phase $\phi_k$ and $b_k$ has phase $\phi_k-t$; moreover,
$|a_k|+|b_k|=M(E_k)$.  These two rays bound a sector of angle
$\pi t$, so our geometric observation shows that
$\mathcal Z(E_k)\ne0$ and proves \eqref{asvasfvasfasf}.

\smallskip
\begingroup\itshape
\noindent Claim 2.
The charges $\mathcal Z(E_k)$ can be split and regrouped into
finitely many nonzero vectors whose arguments lie in
$(-\pi t,0]$,
whose sum remains $\mathcal Z(E)$,
 and whose total length is at least
\[
 c\sum_k|\mathcal Z(E_k)|.
\]
\par\endgroup

\noindent\emph{Proof of Claim~2.}
An individual factor $E_k$ need not belong to $\mathcal F$,
so $\mathcal Z(E_k)$ may have positive imaginary part.  But each 
part $E^{(j)}$ of the filtration~\eqref{qerfqrfqrgqrwgqwrgw} is a subobject of $E\in\mathcal F$, and hence
\begin{equation}\label{qwergweg2ergwegweq}
 \sum_{k\le j}\operatorname{Im}\mathcal Z(E_k)=\operatorname{Im}\mathcal Z(E^{(j)})\le0
 \qquad(1\le j\le s).
\end{equation}

Start with the ordered list
$\mathcal Z(E_1),\ldots,\mathcal Z(E_s)$.
Whenever the list contains a term with positive imaginary part,
choose the first such term and an earlier term with negative
imaginary part; such an earlier
term exists by the partial-sum inequality \eqref{qwergweg2ergwegweq}.  If their imaginary
parts cancel, remove the two terms from the ordered list and record them as a pair.
Otherwise, split the term whose imaginary part has larger
magnitude into two positive real multiples, one with imaginary
part opposite to that of the other chosen term.  Remove this
piece and the other chosen term from the ordered list as a pair, keeping the second
piece in its original position.
The partial-sum inequalities \eqref{qwergweg2ergwegweq} remain valid.  Each operation removes at least one entry.
Repeating therefore gives finitely many pairs and leaves only
unpaired vectors with nonpositive imaginary part.

Each recorded pair consists of positive real multiples $u,v$
of $\mathcal Z(E_k),\mathcal Z(E_l)$ for some $k<l$, with
$\operatorname{Im}u=-\operatorname{Im}v<0$.
Positive scaling preserves arguments, so Claim~1 gives
\[
       0<\arg\mathcal Z(E_l)/\pi-\arg\mathcal Z(E_k)/\pi
           \le\phi_l-(\phi_k-t)<t.
\]
The last inequality uses $\phi_l<\phi_k$.
Thus the sector between their rays has angle less than $\pi t$
and contains the positive real ray but not the negative real
ray.  Their real sum lies in this sector, and the projection
estimate makes it nonzero and hence positive, with
\[
 u+v>0,\qquad |u+v|\ge c(|u|+|v|).
\]

Replace each pair by its positive real sum.  By Claim~1, the
unpaired vectors have phases in $(-t,0]$, so all resulting
vectors lie in this sector.  Splitting preserved the total
length, each pair retains at least a factor $c$, and unpaired
vectors keep their full length.  The total length is therefore
at least $c\sum_k|\mathcal Z(E_k)|$, while the sum remains
$\mathcal Z(E)$.  This proves Claim~2.

\smallskip
\begingroup\itshape
\noindent Claim 3.
The total charge satisfies $ \arg\mathcal Z(E)/\pi\in(-t,0]$ and
\begin{equation}\label{qsavzdfgadrgadrh}
 |\mathcal Z(E)|\ge c^2\sum_k|\mathcal Z(E_k)|.
\end{equation}
\par\endgroup

\noindent\emph{Proof of Claim~3.}
Apply the projection estimate to the vectors supplied by
Claim~2.  They lie in a sector of angle $\pi t$, so the length
of their sum is at least $c$ times their total length.  This
gives \eqref{qsavzdfgadrgadrh} and, in particular,
$\mathcal Z(E)\ne0$.  Their sum lies in the same sector;
because every vector lies strictly above its lower boundary,
the phase of the sum is strictly greater than $-t$ as well.

Combining Claims~1 and~3 now gives
\[
 |\mathcal Z(E)|\ge c^2\sum_k|\mathcal Z(E_k)|
       \ge c^3\sum_kM(E_k)=c^3M(E).
\]
Together with the phase bound in Claim~3, this proves the lemma.
\end{proof}

\begin{lemma}\label{lem:all-degree-stability-function}
Every nonzero object $E$ of the tilted heart $\mathcal A'$
has central charge in the upper half-plane or on the negative
real axis:
\[
 \mathcal Z(E)\in\{z\in\C:\operatorname{Im}z>0\}\cup\R_{<0}.
\]
\end{lemma}

\begin{proof}
We first check the two kinds of objects used to form the tilted
heart.  Every nonzero $T\in\mathcal T$ is high, so either
$\operatorname{Im}\mathcal Z(T)>0$ or $p(T)=0$.
In the latter case, all nonzero vertex objects have phase~$1$.
Writing $T=(X,Y,m)$, we therefore have
$
 z(X)=-a$, $z(Y)=-b$,
 $\mathcal Z(T)=-a+b e^{\pi i/d}$,
where $a,b\ge0$ and $a+b>0$.
Since $\sin(\pi/d)>0$, this charge has positive imaginary
part if $b>0$, and is negative real if $b=0$.

For a nonzero $F\in\mathcal F$,
Lemma~\ref{lem:all-degree-mass} gives $\mathcal Z(F)\ne0$ and
$\arg\mathcal Z(F)/\pi\in(-t,0]$.
Shifting changes the sign of the charge, so
\[
 \mathcal Z(F[1])=-\mathcal Z(F),\qquad
 \arg\mathcal Z(F[1])/\pi\in(1-t,1].
\]
Thus $\mathcal Z(F[1])$ also lies in the upper half-plane or
on the negative real axis.

Finally, every nonzero $E\in\mathcal A'$ fits into the
canonical exact sequence in $\mathcal A'$
\[
 0\longrightarrow F_E[1]\longrightarrow E\longrightarrow T_E
 \longrightarrow0,
 \qquad F_E\in\mathcal F,\quad T_E\in\mathcal T.
\]
Both factors have nonnegative imaginary charge.  If
$\operatorname{Im}\mathcal Z(E)=0$, additivity forces both
imaginary parts to vanish.  Each nonzero factor then has
negative real charge, and at least one factor is nonzero.
Their sum is therefore negative real as well, proving the assertion.
\end{proof}

\subsection{Action of the reflection functor}

Path reversal functor $\Pi$ permutes the vertex objects, so it preserves
the  heart $\mathcal H$ and its aisle $V$.
The symmetry $r_j=r_{d-j}$ shows that it also preserves the
weighted charges $z(X)$ and $z(Y)$, and hence $\mathcal Z$ and $p$.
It therefore preserves high and low objects, the torsion pair
$(\mathcal T,\mathcal F)$, and the tilted aisle $U$.
Using $R^d\simeq\Pi[1]$ from
Lemma~\ref{lem:all-degree-reflection-periodicity}, we obtain
\begin{equation}\label{eq:all-degree-aisle-periodicity}
 R^dU=\Pi U[1]=U[1].
\end{equation}

For $E=(X,Y,m)$, Definition~\ref{def:reflection-path-reversal}
gives $RE=(Y,\Cone(m))$.  Additivity of $z$ and the identity
$z(FY)=\lambda z(Y)$ give
\[
 z(\Cone(m))=z(FY)-z(X)=\lambda z(Y)-z(X).
\]
Since $1+\lambda q=-q^2$ and $-q=e^{\pi i/d}$, it follows that
\begin{equation}\label{eq:all-degree-charge-rotation}
 \mathcal Z(RE)=z(Y)+q\bigl(\lambda z(Y)-z(X)\bigr)
      =-q\bigl(z(X)+qz(Y)\bigr)=e^{\pi i/d}\mathcal Z(E).
\end{equation}

\begin{lemma}\label{lem:all-degree-R-aisle}
The reflection functor sends the tilted aisle into itself:
$RU\subset U$.
\end{lemma}

\begin{proof}
First, $R$ sends the heart $\mathcal H$ into its aisle $V$.
Indeed, the cone of a morphism in the abelian category
$\mathcal A$ has $\mathcal A$-cohomology given by the kernel
in degree $-1$ and the cokernel in degree $0$, with all other
cohomology zero.  The formula for $R$ therefore shows that,
for $E\in\mathcal H$, every vertex object of $RE$ has
$\mathcal A$-cohomology in nonpositive degrees.

Every object of $V$ is a finite extension of objects in
$\mathcal H[n]$ with $n\ge0$.  Their images under $R$ lie in
$R(\mathcal H)[n]\subset V[n]\subset V$.
Since $R$ is exact and $V$ is closed under extensions, we obtain
$RV\subset V$.  Commuting $R$ with the shift then gives
\[
 R(V[1])=(RV)[1]\subset V[1]\subset U.
\]

Every object of $U$ is an extension of an object of $\mathcal T$
by one of $V[1]$.  Because $U$ is closed under extensions,
it remains to show that $RE\in U$ for $E\in\mathcal T$.
We already know that $RE\in V$, so the definition of $U$
reduces this to showing that $H^0_{\mathcal H}(RE)\in\mathcal T$.

Write $E=(X,Y,m)$.  By
Definition~\ref{def:reflection-path-reversal}, the source vertex
$x_j$ of $RE$ carries $Y_j$, and the target vertex $y_k$ carries
the cone of $m_k:X_k\to Y_{k-1}\oplus Y_{k+1}$.
Cohomology with respect to $\mathcal H$ is computed at each
vertex.  Since the $X_k$ and $Y_j$ belong to $\mathcal A$,
taking zeroth cohomology leaves $Y_j$ at the source and gives
$C_k:=\operatorname{coker}(m_k)$ at the target.  Hence
\[
 A=H^0_{\mathcal H}(RE)=(Y,C),\qquad
 C=\operatorname{coker}(m:X\to FY).
\]
The structure map $Y\to FC$ is adjoint to the quotient
$FY\twoheadrightarrow C$ under $F\dashv F$.
Indeed, taking zeroth cohomology of the reflected arrows gives
\[
 Y_j\hookrightarrow Y_{k-1}\oplus Y_{k+1}
     \twoheadrightarrow C_k,\qquad |j-k|=1.
\]
For fixed $j$, these are the components of $Y_j\to(FC)_j$.
For fixed $k$, they combine into the cokernel map
$(FY)_k\twoheadrightarrow C_k$.
This regrouping of the same arrows by source or by target
is precisely the adjunction.
If $A\notin\mathcal T$, its torsion decomposition gives a nonzero quotient
$A\twoheadrightarrow Q=(Y',C')\in\mathcal F$.
Since quotients are computed at the vertices, we have a
commutative square
\[
 \begin{tikzcd}[column sep=large,row sep=large]
 FY \arrow[r,two heads] \arrow[d,two heads]
     & C \arrow[d,two heads]\\
 FY' \arrow[r] & C'.
 \end{tikzcd}
\]
The top and right maps are epimorphisms, so the bottom map
is an epimorphism as well.  At~each~target vertex $y_k$, the
incoming arrows therefore give an epimorphism
$Y'_{k-1}\oplus Y'_{k+1}\twoheadrightarrow C'_k$.
The inverse-reflection formula~\eqref{eq:inverse-vertex-reflection}
replaces $C'_k$ by
\[
 \Cone\bigl(Y'_{k-1}\oplus Y'_{k+1}\to C'_k\bigr)[-1]
 \simeq\ker\bigl(Y'_{k-1}\oplus Y'_{k+1}\to C'_k\bigr).
\]
Indeed, because the map is an epimorphism, the shifted cone
has cohomology only in degree zero, where it is the kernel.
After relabeling, these kernels are the new source objects,
and the $Y'_j$ are the new target objects.  The reversed arrows
are the kernel inclusions followed by the projections to the
neighboring $Y'_j$.
All these vertex objects belong to $\mathcal A$, so
\[
 R^{-1}Q\simeq\bigl(\ker(FY'\to C'),Y'\bigr)\in\mathcal H.
\]

Apply $R^{-1}$ to the composition $RE\to A\twoheadrightarrow Q$.
The resulting map $E\to R^{-1}Q$ has target component
$Y\twoheadrightarrow Y'$, so its cokernel is supported only
at the sources.  Thus, for some $J\in\mathcal A_{\mathcal B}$,
we have an exact sequence in $\mathcal H$
\[
 E\longrightarrow R^{-1}Q
   \longrightarrow(J,0)\longrightarrow0.
\]
Every diagram supported only at the sources belongs to
$\mathcal T$: its nonzero quotients satisfy
$\operatorname{Im}\mathcal Z=p\ge0$, and hence are high.
The image of $E\to R^{-1}Q$ also belongs to $\mathcal T$,
being a quotient of $E$.  Closure of $\mathcal T$ under
extensions therefore gives $R^{-1}Q\in\mathcal T$.

The charge rotation~\eqref{eq:all-degree-charge-rotation} gives
\[
 \mathcal Z(R^{-1}Q)=e^{-\pi i/d}\mathcal Z(Q).
\]
Since $0\ne Q\in\mathcal F$,
\eqref{eq:all-degree-torsion-mass} places this charge in the
sector of phases $(-1,-1/d]$, where the imaginary part is
strictly negative.  This contradicts
Lemma~\ref{lem:all-degree-stability-function}, since
$R^{-1}Q\in\mathcal T\subset\mathcal A'$.
Thus $A\in\mathcal T$ and $RE\in U$.
\end{proof}

\begin{lemma}\label{lem:all-degree-reflection-mass}
For every quiver diagram $E$, one has
\begin{equation}\label{eq:all-degree-reflection-mass}
 M(RE)\le(1+\lambda)M(E),\qquad \lambda=2\cos(\pi/d).
\end{equation}
For $d=2$, equality holds.
\end{lemma}

\begin{proof}
Write $E=(X,Y,m)$.  The reflection functor gives
$RE=(Y,\Cone(m))$.  At each vertex $j$, rotate the cone
triangle to obtain
\[
 (FY)_j\longrightarrow\Cone(m_j)\longrightarrow X_j[1]
       \longrightarrow(FY)_j[1].
\]
The mass triangle inequality~\eqref{eq:mass-triangle-inequality},
together with invariance under shifts and additivity under direct
sums, gives
\[
 m_\sigma(\Cone(m_j))
 \le m_\sigma((FY)_j)+m_\sigma(X_j[1])
 =m_\sigma(Y_{j-1})+m_\sigma(Y_{j+1})+m_\sigma(X_j).
\]
Multiplying by $r_j$ and summing, we use
$r_{j-1}+r_{j+1}=\lambda r_j$ to collect the neighboring target
contributions.  Adding the masses of the source objects $Y_j$
of $RE$, we obtain
\[
 M(RE)\le\sum_j r_j m_\sigma(X_j)
              +(1+\lambda)\sum_j r_j m_\sigma(Y_j)
      \le(1+\lambda)M(E).
\]
For $d=2$, reflection is $(X,Y)\mapsto(Y,X[1])$, and preserves
$M$ exactly.
\end{proof}

\subsection{Reflection blocks}

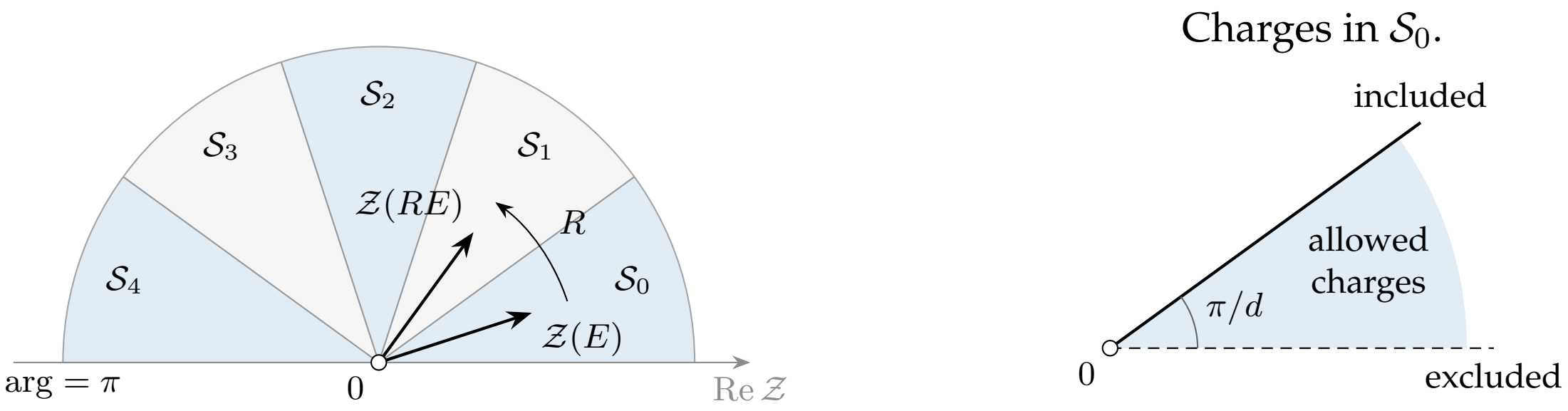
\begin{figure}[htbp]
\centering
\begin{minipage}[c]{.59\linewidth}
\centering
\begin{tikzpicture}[>=Stealth,x=1.15cm,y=1.15cm,font=\footnotesize]
\path[use as bounding box] (-3.2,-.6) rectangle (3.45,2.9);
\coordinate (O) at (0,0);
\foreach \j in {0,...,4} {
 \pgfmathsetmacro{\lo}{36*\j}
 \pgfmathsetmacro{\hi}{36*(\j+1)}
 \pgfmathsetmacro{\mid}{36*\j+18}
 \ifodd\j
  \fill[black!4] (O)--(\lo:2.55)
   arc[start angle=\lo,end angle=\hi,radius=2.55]--cycle;
 \else
  \fill[MidnightBlue!8] (O)--(\lo:2.55)
   arc[start angle=\lo,end angle=\hi,radius=2.55]--cycle;
 \fi
 \node at (\mid:2.16) {$\mathcal S_{\j}$};
}
\draw[black!35] (0:2.55) arc[start angle=0,end angle=180,radius=2.55];
\foreach \angle in {36,72,108,144}
 \draw[black!40] (O)--(\angle:2.55);
\draw[->,black!45] (-2.95,0)--(3,0)
 node[below] {$\operatorname{Re}\mathcal Z$};
\node[below] at (-2.55,0) {$\arg=\pi$};
\draw[->,thick,black] (O)--(18:1.3)
 node[below right,inner sep=2pt] {$\mathcal Z(E)$};
\draw[->,thick,black] (O)--(54:1.3)
 node[above left,inner sep=2pt] {$\mathcal Z(RE)$};
\draw[->,black] (18:1.6)
 arc[start angle=18,end angle=54,radius=1.6];
\node at (36:1.93) {$R$};
\draw[draw=black,fill=white] (O) circle[radius=2pt];
\node[below left] at (O) {$0$};
\end{tikzpicture}
\end{minipage}\hfill
\begin{minipage}[c]{.37\linewidth}
\centering
\textup{Charges in $\mathcal S_0$.}\par\smallskip
\begin{tikzpicture}[>=Stealth,x=1.25cm,y=1.25cm,font=\footnotesize]
\path[use as bounding box] (-.25,-.6) rectangle (3.2,2.1);
\coordinate (O) at (0,0);
\fill[MidnightBlue!8] (O)--(0:2.65)
 arc[start angle=0,end angle=36,radius=2.65]--cycle;
\draw[densely dashed,black] (O)--(0:2.85)
 node[below] {excluded};
\draw[black,thick] (O)--(36:2.85)
 node[above] {included};
\draw[black!60] (0:.65) arc[start angle=0,end angle=36,radius=.65];
\node at (18:.97) {$\pi/d$};
\node[align=center] at (18:2.02) {allowed\\charges};
\draw[draw=black,fill=white] (O) circle[radius=2pt];
\node[below left] at (O) {$0$};
\end{tikzpicture}
\end{minipage}
\caption{Charges in the reflection blocks for $d=5$.}
\label{fig:reflection-block-charges}
\end{figure}

We now use the nested aisles to divide the category into narrow
phase sectors.  For $j\in\Z$, put $U_j=R^jU$.
Lemma~\ref{lem:all-degree-R-aisle}
and~\eqref{eq:all-degree-aisle-periodicity} give
\[
 U_{j+1}\subset U_j,\qquad U_{j+d}=U_j[1].
\]
The exact autoequivalence $R^j$ carries the t-structure with
aisle $U$ and heart $\mathcal A'$ to one with aisle $U_j$ and
heart $R^j\mathcal A'$.  Since the heart of an aisle is its
intersection with the right orthogonal of its shift by $[1]$,
we can write this heart as
\[
 \mathcal H_j=U_j\cap(U_j[1])^\perp
             =U_j\cap U_{j+d}^\perp=R^j\mathcal A'.
\]
We also define the \emph{block}
\[
 \mathcal S_j=U_j\cap U_{j+1}^{\perp}.
\]
An object of $\mathcal S_j$
lies in both $\mathcal H_j$ and $\mathcal H_{j+1}[-1]$.
By Lemma~\ref{lem:all-degree-stability-function}
and~\eqref{eq:all-degree-charge-rotation}, the allowed half-planes
for their rotated stability functions intersect in the sector
\begin{equation}\label{eq:all-degree-block-sector}
 \mathcal Z(E)\in\R_{>0}e^{i\pi\theta},\qquad
 j/d<\theta\le(j+1)/d
 \quad(0\ne E\in\mathcal S_j).
\end{equation}
The lower endpoint is excluded because a stability function on
$\mathcal H_j$ cannot take a nonzero object to the positive real
axis after rotation by $e^{-\pi ij/d}$. The sectors are illustrated in 
Figure~\ref{fig:reflection-block-charges}.

\begin{lemma}\label{lem:all-degree-block-structure}
Every object $E\in\mathcal S_0$ belongs to $\mathcal T\subset\mathcal H$.
Write $E=(X,Y,m)$ and put $C=\operatorname{coker}(FX\to Y)$,
where $FX\to Y$ is adjoint to $m:X\to FY$.
Then there is a distinguished triangle
\begin{equation}\label{eq:all-degree-block-cohomology}
 K\longrightarrow R^{-1}E\longrightarrow(C,0)[-1]
                  \longrightarrow K[1],
\end{equation}
where $K=H^0_{\mathcal H}(R^{-1}E)\in\mathcal F$.
All HN factors of the vertex objects $C_j$ have phases in $(t,1]$,
where $t=1-\frac1d$.
\end{lemma}

\begin{proof}
Let $D=(X,Y,m)\in\mathcal H$, and consider the map $FX\to Y$
adjoint to $m:X\to FY$.  The inverse-reflection
formula~\eqref{eq:inverse-vertex-reflection}, followed by the
shift $[1]$, gives
\[
 R^{-1}D[1]\simeq\bigl(\Cone(FX\to Y),X[1]\bigr).
\]
At each source, the cone has cohomology equal to the kernel
in degree $-1$ and the cokernel in degree $0$.
At each target, $X_j[1]$ has cohomology only in degree $-1$.
Thus all vertex cohomology lies in nonpositive degrees, so
$R^{-1}D[1]\in V$, and
\[
 H^0_{\mathcal H}(R^{-1}D[1])
   =\bigl(\operatorname{coker}(FX\to Y),0\bigr).
\]
This zeroth cohomology diagram is supported only at the sources.
It belongs to $\mathcal T$, since $\operatorname{Im}\mathcal Z=p$
on it and all its quotients.  The definition of $U$ therefore
gives $R^{-1}D[1]\in U$, proving
$\mathcal H[1]\subset RU$.

Now take $E\in\mathcal S_0$.  Since $U[1]\subset RU$, we have
$E\in U\cap(U[1])^\perp=\mathcal A'$.  Its tilted-heart sequence is
\[
 0\longrightarrow F_E[1]\longrightarrow E\longrightarrow T_E
   \longrightarrow0.
\]
But $F_E[1]\in\mathcal H[1]\subset RU$, whereas
$E\in(RU)^\perp$, so the displayed injection is zero.
Thus $F_E=0$ and $E\in\mathcal T\subset\mathcal H$.
The inverse-reflection
formula~\eqref{eq:inverse-vertex-reflection} gives
\[
 R^{-1}E\simeq\bigl(\Cone(FX\to Y)[-1],X\bigr).
\]
The shifted cone has the kernel in degree $0$ and the cokernel
in degree $1$, so
\[
 H^0_{\mathcal H}(R^{-1}E)=\bigl(\ker(FX\to Y),X\bigr)=K,
 \qquad H^1_{\mathcal H}(R^{-1}E)=(C,0),
\]
with all other cohomology zero.  The cohomology truncations
therefore give the triangle~\eqref{eq:all-degree-block-cohomology}.
For $T\in\mathcal T$, applying $\Hom(T,-)$ gives
\[
 \Hom(T,K)\simeq\Hom(T,R^{-1}E)=0.
\]
The isomorphism uses $\Hom(T,(C,0)[-k])=0$ for $k=1,2$,
since $T$ and $(C,0)$ lie in the heart $\mathcal H$.
The vanishing follows from $T\in U$ and $R^{-1}E\in U^\perp$.
Thus $K\in\mathcal F$.

At each target $y_j$, the cokernel map $Y_j\twoheadrightarrow C_j$
kills the images of both incoming arrows from $X_{j-1}$ and
$X_{j+1}$.  Together with the zero maps at the sources, these
maps therefore define a quotient $E\twoheadrightarrow(0,C)$
in $\mathcal H$.  Since $E\in\mathcal T$ and $\mathcal T$ is
closed under quotients, $(0,C)\in\mathcal T$.

Suppose some $C_j$ has an HN factor of phase at most $t$, and let
$C_j\twoheadrightarrow B$ be its lowest-phase HN quotient,
of phase $0<\phi\le t$.
Sending all other vertices to zero gives a quotient of $(0,C)$
with $B$ at the target $y_j$ and zero elsewhere.
For this diagram, $p=r_j\operatorname{Im}Z(B)>0$, since
$0<\phi<1$.  Its charge is $r_jqZ(B)$, which has phase
$\phi-t\in(-t,0]$ because $q=e^{-\pi it}$.
Its imaginary charge is therefore nonpositive, so it is low.
This contradicts the defining property of $\mathcal T$ that
every nonzero quotient is high.  Thus all vertex HN phases of
$C$ lie in $(t,1]$.
\end{proof}

\begin{lemma}\label{lem:all-degree-block-support}
Recall that $c=\cos(\pi t/2)=\sin(\pi/(2d))$.
Every object $E\in\mathcal S_0$ satisfies
\begin{equation}\label{eq:all-degree-initial-block-mass}
 M(E)\le(1+\lambda)c^{-4}|\mathcal Z(E)|.
\end{equation}
\end{lemma}

\begin{proof}
Let $K$ and $C$ be as in Lemma~\ref{lem:all-degree-block-structure}.
Expanding each $Z(C_j)$ into the charges of its HN factors
expresses $z(C)=\sum_j r_jZ(C_j)$ as a sum of vectors in the
sector of phases $(t,1]$.  Their total length is
$\sum_j r_j m_\sigma(C_j)=M((C,0))$.
Since this sector has angle $\pi/d\le\pi t$, projection onto
its middle ray retains at least a factor $c$, as in the proof
of Lemma~\ref{lem:all-degree-mass}.  Consequently,
\begin{equation}\label{eq:all-degree-cokernel-mass}
 M((C,0))\le c^{-1}|z(C)|.
\end{equation}

We can now estimate the two cohomology objects in
the triangle~\eqref{eq:all-degree-block-cohomology}.
Additivity of the charge and~\eqref{eq:all-degree-charge-rotation} give
\[
 e^{-\pi i/d}\mathcal Z(E)=\mathcal Z(K)-z(C).
\]
Since $K\in\mathcal F$,~\eqref{eq:all-degree-torsion-mass}
puts $\mathcal Z(K)$ in the sector of phases $(-t,0]$ when
nonzero.  The phase bound for $C$ puts $-z(C)$ in
$(-1/d,0]\subseteq(-t,0]$ when nonzero.
Applying the same projection estimate to these two charges gives
\begin{equation}\label{eq:all-degree-block-charge-sum}
 |\mathcal Z(K)|+|z(C)|\le c^{-1}|\mathcal Z(E)|.
\end{equation}
Combining these estimates, we obtain
\[
 M(R^{-1}E)\le M(K)+M((C,0))
 \le c^{-3}\bigl(|\mathcal Z(K)|+|z(C)|\bigr)
 \le c^{-4}|\mathcal Z(E)|.
\]
The first inequality applies the mass triangle
inequality~\eqref{eq:mass-triangle-inequality} at each vertex
of~\eqref{eq:all-degree-block-cohomology}, using invariance
under shifts.  The second uses the bound for $K$
in~\eqref{eq:all-degree-torsion-mass}, the bound for $C$
in~\eqref{eq:all-degree-cokernel-mass}, and $c\le1$.
The last follows from~\eqref{eq:all-degree-block-charge-sum}.
Applying~\eqref{eq:all-degree-reflection-mass} to $R^{-1}E$
now gives the required bound for $M(E)$.
\end{proof}

\begin{lemma}\label{lem:all-degree-uniform-block-support}
Put $B_d=(1+\lambda)^dc^{-4}$.  For every $j\in\Z$ and
$E\in\mathcal S_j$, one has
\begin{equation}\label{eq:all-degree-block-mass}
 M(E)\le B_d|\mathcal Z(E)|.
\end{equation}
In particular, every nonzero block object has nonzero charge.
Fix a norm on $\Lambda_\R$ and a support constant $C_\sigma$
for $\sigma$.  With the corresponding sum norm on the vertex
lattice $\Lambda^{\oplus2(d-1)}$, one also has
\begin{equation}\label{eq:all-degree-vertex-support}
 \|v(E)\|\le C_\sigma M(E)
 \le C_\sigma B_d|\mathcal Z(E)|
 \qquad(E\in\mathcal S_j).
\end{equation}
\end{lemma}

\begin{proof}
First suppose $0\le j<d$.  Then $R^{-j}E\in\mathcal S_0$.
Applying~\eqref{eq:all-degree-reflection-mass} $j$ times and
then~\eqref{eq:all-degree-initial-block-mass}, we obtain
\[
 M(E)\le(1+\lambda)^jM(R^{-j}E)
 \le(1+\lambda)^{j+1}c^{-4}|\mathcal Z(E)|
 \le B_d|\mathcal Z(E)|.
\]
Here~\eqref{eq:all-degree-charge-rotation} gives
$|\mathcal Z(R^{-j}E)|=|\mathcal Z(E)|$, and the last inequality
uses $j+1\le d$ and $1+\lambda\ge1$.
Since $\mathcal S_{j+d}=\mathcal S_j[1]$ and shifts preserve both
$M$ and the absolute value of the charge, the bound holds for
every $j\in\Z$.  If $E\ne0$, at least one vertex object is
nonzero, so $M(E)>0$ and hence $\mathcal Z(E)\ne0$.

For any object $A\in\mathcal C$, applying the support
inequality~\eqref{eq:background-support} to its HN factors with
respect to $\sigma$
and summing gives $\|v(A)\|\le C_\sigma m_\sigma(A)$.
Summing this inequality over all vertices and using $r_k\ge1$
gives $\|v(E)\|\le C_\sigma M(E)$.
Together with~\eqref{eq:all-degree-block-mass}, this
proves~\eqref{eq:all-degree-vertex-support}.
\end{proof}

The stronger estimate~\eqref{eq:all-degree-block-mass} will also
provide the mass--Hom bound after constructing HN filtrations.

\subsection{Harder--Narasimhan filtrations}

Successive truncations along
$U\supset RU\supset\cdots\supset R^dU=U[1]$ give filtrations
by blocks.  We refine these to HN filtrations using the support
bound~\eqref{eq:all-degree-vertex-support}.  
A~similar use of nested aisles appears
in~\cite[\S4]{OAIQuintic}.

\begin{lemma}\label{lem:all-degree-block-filtrations}
Let $a<b$.  Every object $E\in U_a\cap U_b^\perp$ has a
filtration by distinguished triangles with successive factors
in $\mathcal S_{b-1},\ldots,\mathcal S_a$, in that order,
allowing zero factors.  Every object belongs to such an
intersection for suitable integers $a<b$.
For every $j\in\Z$, the pair $(\mathcal T_j,\mathcal S_j)$ is
a torsion pair in $\mathcal H_j$, where
\[
 \mathcal T_j=U_{j+1}\cap\mathcal H_j
   =\langle\mathcal S_{j+d-1},\ldots,\mathcal S_{j+1}\rangle.
\]
Moreover, $\Hom(\mathcal S_l,\mathcal S_j)=0$ whenever $l>j$.
\end{lemma}

\begin{proof}
Suppose $E\in U_a\cap U_b^\perp$, with $a<b$.
Truncation at the aisle $U_{a+1}$ gives a triangle
\[
 E'\longrightarrow E\longrightarrow E''\longrightarrow E'[1],
 \qquad E'\in U_{a+1},\quad E''\in U_{a+1}^\perp.
\]
The rotated triangle $E\to E''\to E'[1]$ shows that
$E''\in U_a$, since both outside terms lie in this aisle.
Hence $E''\in U_a\cap U_{a+1}^\perp=\mathcal S_a$.
The remaining term satisfies $E'\in U_b^\perp$ as well.
Indeed, for $A\in U_b$ we have $A[1]\in U_{a+1}$, so
$\Hom(A,E''[-1])=\Hom(A[1],E'')=0$.
Applying $\Hom(A,-)$ to $E''[-1]\to E'\to E$ then gives
$\Hom(A,E')=0$.

We can therefore repeat the construction with
$E'\in U_{a+1}\cap U_b^\perp$.
It ends when the remaining term lies in $U_b\cap U_b^\perp=0$,
giving a filtration with factors in
$\mathcal S_{b-1},\ldots,\mathcal S_a$, in that order.
Boundedness and $U_{j+d}=U_j[1]$ ensure that every object
belongs to $U_a\cap U_b^\perp$ for suitable integers $a<b$.

For $E\in\mathcal H_j$, take $a=j$ and $b=j+d$.
All the remaining terms and factors then lie in $\mathcal H_j$,
so the triangles are short exact sequences in this heart.
The first truncation gives the torsion decomposition for the pair
\[
 (\mathcal T_j,\mathcal S_j),\qquad
 \mathcal T_j=U_{j+1}\cap\mathcal H_j
   =\langle\mathcal S_{j+d-1},\ldots,\mathcal S_{j+1}\rangle,
\]
where orthogonality follows from $\mathcal S_j\subset U_{j+1}^\perp$.
In particular, $\mathcal S_j$ is closed under subobjects and
extensions in $\mathcal H_j$.
More generally, if $l>j$, then $\mathcal S_l\subset U_{j+1}$,
so $\Hom(\mathcal S_l,\mathcal S_j)=0$.
\end{proof}

By Lemma~\ref{lem:all-degree-stability-function}
and~\eqref{eq:all-degree-charge-rotation}, every nonzero
$E\in\mathcal H_j$ has a unique phase $\phi\in(j/d,j/d+1]$
such that
\[
 \mathcal Z(E)\in\R_{>0}e^{\pi i\phi}.
\]
We call $E$ semistable if no nonzero subobject in
$\mathcal H_j$ has larger phase.
The following argument is fairly standard; compare
Collins--Polishchuk~\cite[Proposition~3.1 and Lemma~3.2]{CP}.

\begin{lemma}\label{lem:all-degree-block-HN}
For every $j\in\Z$, every nonzero object $E\in\mathcal S_j$
has a finite filtration by subobjects in $\mathcal H_j$ whose
successive quotients are nonzero semistable objects of
$\mathcal S_j$ with strictly decreasing phases.
\end{lemma}

\begin{proof}
Fix $0\ne E\in\mathcal S_j$, and use the torsion pair
$(\mathcal T_j,\mathcal S_j)$ of
Lemma~\ref{lem:all-degree-block-filtrations}.
Call a subobject $A\subset E$ in $\mathcal H_j$ \emph{saturated}
if $E/A\in\mathcal S_j$.

It suffices to test semistability on saturated subobjects.
Indeed, given $0\ne G\subset E$, let $\overline G$ be the
inverse image in $E$ of the $\mathcal T_j$ part of $E/G$.
Then $E/\overline G\in\mathcal S_j$, so $\overline G$ is
saturated, and $\overline G/G\in\mathcal T_j$.
The latter has a filtration by blocks with indices greater than $j$.
Their charges have phases in $((j+1)/d,j/d+1]$, above the phase
of $G\in\mathcal S_j$.  All these rays, including that of $G$,
lie in a sector of angle strictly less than $\pi$.
Adding their charges therefore cannot lower the phase of $G$.
Thus any destabilizing subobject gives a destabilizing saturated
subobject.

There are only finitely many possible charges of saturated
subobjects of $E$.  To see this, project onto the middle ray of
the block sector:
\[
 b_j(w)=\operatorname{Re}(e^{-\pi i(j+1/2)/d}w).
\]
By~\eqref{eq:all-degree-block-sector}, the sector has angle
$\pi/d\le\pi t$.  The projection
estimate~\eqref{eq:sector-projection-estimate} therefore gives
$b_j(\mathcal Z(A))\ge c|\mathcal Z(A)|>0$ for every nonzero
$A\in\mathcal S_j$.
If $A\subset E$ is saturated, its quotient also belongs to
$\mathcal S_j$, so $b_j(\mathcal Z(A))\le b_j(\mathcal Z(E))$.
The support estimate~\eqref{eq:all-degree-vertex-support} now gives
\[
 \|v(A)\|\le C_\sigma B_d|\mathcal Z(A)|
 \le C_\sigma B_dc^{-1}b_j(\mathcal Z(E)).
\]
A bounded subset of the vertex lattice is finite, and
$\mathcal Z(A)$ depends only on $v(A)$, proving the assertion.

Choose a nonzero saturated subobject $A_1\subset E$ of maximal
phase and, among those, maximal charge modulus.
Such a choice exists because $E$ itself is saturated and only
finitely many charges occur.
A saturated subobject of $A_1$ is saturated in $E$: the quotient
of $E$ by it is an extension of two objects of $\mathcal S_j$.
Thus no saturated subobject of $A_1$ has larger phase, and the
criterion just proved shows that $A_1$ is semistable.
Every nonzero saturated subobject $B\subset E/A_1$ has smaller
phase than $A_1$.  Otherwise its inverse image in $E$ would be
saturated with larger phase, or with the same phase and larger
charge modulus, contrary to the choice of $A_1$.

Repeat this construction on $E/A_1$.
The resulting factors are semistable with strictly decreasing
phases.  Their successive inverse images
$0=E_0\subset E_1\subset\cdots\subset E$ are all saturated
in $E$, because each remaining quotient lies in $\mathcal S_j$.
Their charges therefore belong to the finite set found above,
whereas $b_j(\mathcal Z(E_r))$ strictly increases at each step.
The construction must therefore reach $E$ after finitely many
steps, giving its HN filtration.
\end{proof}

\begin{lemma}\label{lem:all-degree-grid}
The charge $\mathcal Z$ defines a stability condition
$\Sigma=(\mathcal Z,\mathcal P_\Sigma)$ with heart $\mathcal A'$
and support property on $\Lambda^{\oplus2(d-1)}$.
The condition $\Sigma$ is of Gepner type $(R,1/d)$.
For every quiver diagram $E$,
\begin{equation}\label{eq:all-degree-global-mass}
 M(E)\le B_dm_\Sigma(E),\qquad
 B_d=\bigl(1+2\cos(\pi/d)\bigr)^d\sin(\pi/(2d))^{-4}.
\end{equation}
\end{lemma}

\begin{proof}
For each $\phi\in\R$, there is a unique integer $j$ such that
$j/d<\phi\le(j+1)/d$.
Define $\mathcal P_\Sigma(\phi)$ to be the full subcategory
consisting of zero and the objects of $\mathcal S_j$ that are
semistable of phase $\phi$ in the heart $\mathcal H_j$.
The subobject and quotient criteria for semistability
\cite[Definition~2.2 and the following paragraph]{Bridgeland}
imply that extensions and finite direct sums of semistable
objects of phase $\phi$ in $\mathcal H_j$ are again semistable
of phase $\phi$.
A nonzero direct summand is both a subobject and a quotient,
so it has the same phase; its subobjects are also subobjects
of the original object, so it is semistable as well.
By Lemma~\ref{lem:all-degree-block-filtrations}, the block
$\mathcal S_j$ is closed under extensions and subobjects;
hence all these objects remain in $\mathcal S_j$.
Thus $\mathcal P_\Sigma(\phi)$ is additive and closed under
extensions and direct summands.

We next check the vanishing of morphisms in decreasing phase.
Let $A\in\mathcal P_\Sigma(\phi)$ and
$B\in\mathcal P_\Sigma(\psi)$, with $\phi>\psi$.
If $A$ and $B$ lie in the same block $\mathcal S_j$, any
morphism $f:A\to B$ has an image $I$ in the abelian heart
$\mathcal H_j$.
If $f\ne0$, then $I\ne0$ is a quotient of $A$ and a subobject
of $B$.  The quotient criterion for semistability gives
$\phi(I)\ge\phi$, while the subobject criterion gives
$\phi(I)\le\psi$, a contradiction.
Thus $f=0$.  If $A$ and $B$ lie in different blocks, the block
containing $A$ has larger index, and the vanishing follows
from Lemma~\ref{lem:all-degree-block-filtrations}.

Every object has the block filtration of
Lemma~\ref{lem:all-degree-block-filtrations}.
Refine each block factor by its HN filtration from
Lemma~\ref{lem:all-degree-block-HN}.
The phases decrease within each block.  They also decrease
when passing to a block of smaller index, since every phase
in $\mathcal S_j$ is greater than $j/d$, whereas every phase
in a lower block is at most $j/d$.
The refined filtration is therefore an HN filtration for
$\mathcal P_\Sigma$.

The shift $[1]$ is an exact equivalence from $\mathcal H_j$
to $\mathcal H_{j+d}$, taking $\mathcal S_j$ to
$\mathcal S_{j+d}$.
It identifies subobjects in the two hearts and, because
$\mathcal Z(E[1])=-\mathcal Z(E)$, increases phases by $1$.
Thus it preserves semistability and gives
$ \mathcal P_\Sigma(\phi+1)=\mathcal P_\Sigma(\phi)[1]$.
We have now verified the slicing axioms.
The condition on the central charge holds by the definition
of $\mathcal P_\Sigma(\phi)$, so
$\Sigma=(\mathcal Z,\mathcal P_\Sigma)$ is a stability condition.
The block filtration in $\mathcal H_0=\mathcal A'$ uses the
blocks $\mathcal S_0,\ldots,\mathcal S_{d-1}$.
Its refinement shows that $\mathcal A'$ is the extension
closure of the subcategories $\mathcal P_\Sigma(\phi)$ with
$0<\phi\le1$.  Hence $\mathcal P_\Sigma((0,1])=\mathcal A'$.

To check the Gepner property, use the exact equivalence
$R:\mathcal H_j\xrightarrow{\sim}\mathcal H_{j+1}$, which
carries $\mathcal S_j$ to $\mathcal S_{j+1}$.
Every subobject of $RE$ in $\mathcal H_{j+1}$ is the image
under $R$ of a subobject of $E$ in $\mathcal H_j$.
By~\eqref{eq:all-degree-charge-rotation}, the phases of $E$
and all its subobjects increase by the same amount $1/d$.
Thus $E$ is semistable if and only if $RE$ is semistable, and
\[
 R\mathcal P_\Sigma(\phi)=\mathcal P_\Sigma(\phi+1/d).
\]
Together with~\eqref{eq:all-degree-charge-rotation}, this is
the Gepner property of Definition~\ref{def:phase-shift-gepner}.

Every $\Sigma$-semistable object lies in one block, so
\eqref{eq:all-degree-vertex-support} gives the support property.
Finally, let $E_1,\ldots,E_s$ be the HN factors of an arbitrary
quiver diagram $E$.
The triangle inequality~\eqref{eq:all-degree-mass-triangle}
for $M$ and the block
estimate~\eqref{eq:all-degree-block-mass} give
\[
 M(E)\le\sum_{k=1}^s M(E_k)
 \le B_d\sum_{k=1}^s|\mathcal Z(E_k)|
 =B_dm_\Sigma(E),
\]
which proves~\eqref{eq:all-degree-global-mass}.
\end{proof}

\begin{lemma}\label{lem:all-degree-vertex-numerical}
The vertex class map
\[
 K_0(\mathcal C\otimes\widehat{\mathcal Q}_d)
 \longrightarrow\Lambda^{\oplus2(d-1)},\qquad
 [E]\longmapsto\bigl(v(E_j)\bigr)_j,
\]
is numerical: it factors through
$K_{\num}(\mathcal C\otimes\widehat{\mathcal Q}_d)$.
Here $j$ ranges over all source and target vertices, and
$v(E_j)$ is the given numerical class of the vertex object
$E_j\in\mathcal C$.
\end{lemma}

\begin{proof}
Fix a vertex $j$, and recall the projective diagrams
$P_j\otimes_{\C}B$ used in
the resolution~\eqref{eq:quiver-projective-resolution}.
Their representable-module description in the proof of
Lemma~\ref{lem:quiver-diagram-perfectness} gives the evaluation
adjunction
\[
 \Hom(P_j\otimes_{\C}B,E[n])
 \simeq\Hom_{\mathcal C}(B,E_j[n])
 \qquad(B\in\mathcal C,\ n\in\Z).
\]
Taking alternating sums of dimensions gives
\begin{equation}\label{eq:all-degree-evaluation-euler}
 \chi_{\mathcal C}(B,E_j)
 =\chi_{\mathcal C\otimes\widehat{\mathcal Q}_d}
       (P_j\otimes_{\C}B,E).
\end{equation}

Both categories have Serre functors, so their left and right
Euler radicals agree, as recalled in Section~\ref{sec:background}.
If a class $\xi\in K_0(\mathcal C\otimes\widehat{\mathcal Q}_d)$
is numerically trivial, the right-hand side
of~\eqref{eq:all-degree-evaluation-euler}, extended by linearity
to $\xi$, vanishes for every $B\in\mathcal C$.
Thus its class at vertex $j$ is numerically trivial in
$\mathcal C$ and is killed by the numerical class map $v$.
This holds at every vertex, proving the lemma.
\end{proof}

\begin{corollary}\label{cor:all-degree-mass-hom}
If $\sigma$ admits a total mass--Hom bound, then so does $\Sigma$.
More precisely, if a split generator $G$ of $\mathcal C$ satisfies
\[
 \sum_{n\in\Z}\dim_{\C}\Hom(G,E[n])\le\kappa m_\sigma(E)
 \qquad\text{for every }E\in\mathcal C,
\]
then the diagram category has a split generator $P$ satisfying
\[
 \sum_{n\in\Z}\dim_{\C}\Hom(P,E[n])\le\kappa B_dm_\Sigma(E)
 \qquad\text{for every diagram }E.
\]
\end{corollary}

\begin{proof}
Take $P=\bigoplus_jP_j\otimes G$, where $P_j$ is the projective
representation at vertex $j$.
The resolution~\eqref{eq:quiver-projective-resolution} shows that
$P$ split-generates the diagram category.
Evaluation adjunction, $r_j\ge1$, and
\eqref{eq:all-degree-global-mass} give
\[
 \begin{aligned}
 \sum_{n\in\Z}\dim_{\C}\Hom(P,E[n])
 &=\sum_j\sum_{n\in\Z}\dim_{\C}\Hom(G,E_j[n])\\
 &\le\kappa\sum_j m_\sigma(E_j)
 \le\kappa M(E)\le\kappa B_dm_\Sigma(E).
 \end{aligned}
\]
Since $P$ is a split generator,
Lemma~\ref{lem:mass-hom-generator} extends this bound to every
fixed object $A\in\mathcal C\otimes\widehat{\mathcal Q}_d$.
More precisely, there is a constant $C_A>0$, independent of
$E$, such that for every quiver diagram $E$,
$
 \sum\limits_{n\in\Z}\dim_{\C}\Hom(A,E[n])
 \le C_A m_\Sigma(E)
$.
\end{proof}

\subsection{Compatibility with autoequivalences}

Suppose that the  stability condition $\sigma$ is of
Gepner type $(\Psi,1/d)$ for an exact enhanced autoequivalence
$\Psi$ of $\mathcal C$.  Applying $\Psi$ at every vertex and
along every arrow gives an autoequivalence of quiver diagrams,
which we also denote by $\Psi$.  Its dg lift respects the sums
and cones defining reflection, so this induced functor is exact
and commutes with $R$.
We will show that $\Sigma$ is also of Gepner type $(\Psi,1/d)$.
The charge identity follows directly from
$\mathcal Z(X,Y)=z(X)+qz(Y)$:
\begin{equation}\label{eq:all-degree-Psi-charge-rotation}
 \mathcal Z(\Psi E)=e^{\pi i/d}\mathcal Z(E).
\end{equation}
To prove the corresponding statement for semistability, we first
show that $\Psi$ sends the tilted aisle $U$ into itself.  

\begin{lemma}\label{lem:all-degree-Psi-aisle}
One has $\Psi U\subset U$.
\end{lemma}

\begin{proof}
Recall that $V$ is the aisle of the  heart $\mathcal H$,
and every object of $U$ is an extension of an object of
$\mathcal T$ by one of $V[1]$.
Applying $\Psi$ increases every  phase in $\cC$ by $1/d$,
so $\Psi V\subset V$ and hence $\Psi(V[1])\subset V[1]\subset U$.
It remains to consider $E\in\mathcal T$.  Since $\Psi E\in V$,
we must show that $H^0_{\mathcal H}(\Psi E)\in\mathcal T$.

At a vertex $j$, the HN factors of $E_j$ occur in decreasing
order of phase.  Let $E_j^{>t}\subset E_j$ be the subobject
given by the initial part of this filtration consisting of all
factors of phase greater than $t$; take $E_j^{>t}=0$ if there
are no such factors.  The quotient $E_j/E_j^{>t}$ then has
all HN phases in $(0,t]$.
Every arrow $E_j\to E_k$ maps $E_j^{>t}$ into $E_k^{>t}$,
because its composition with $E_k\to E_k/E_k^{>t}$ vanishes:
there are no morphisms from objects with HN phases greater
than $t$ to objects with HN phases at most $t$.
Thus the $E_j^{>t}$ form a quiver subdiagram of $E$.
Let $E_{\mathrm{low}}$ be the quotient diagram, whose vertex
objects are $E_j/E_j^{>t}$.  It belongs to $\mathcal T$ because
$E\in\mathcal T$ and $\mathcal T$ is closed under quotients
in $\mathcal H$.
Applying $\Psi$ sends the  subdiagram $E_j^{>t}$ into $\mathcal H[1]$
and $E_{\mathrm{low}}$ into $\mathcal H$.  Taking cohomology gives
\[
           H^0_{\mathcal H}(\Psi E)=\Psi E_{\mathrm{low}}\in\mathcal H.
\]
Suppose that $\Psi E_{\mathrm{low}}$ does not belong to $\mathcal T$.
Its torsion decomposition gives an exact sequence in~$\mathcal H$
\[
 0\longrightarrow L\longrightarrow\Psi E_{\mathrm{low}}
   \longrightarrow Q\longrightarrow0,
 \qquad 0\ne Q\in\mathcal F.
\]
The  objects in $\cC$ with phases in $(1/d,1]$ are closed
under quotients in $\mathcal A$.  Thus the vertices of $Q$, like
those of $\Psi E_{\mathrm{low}}$, have phases in this interval,
and $\Psi^{-1}Q\in\mathcal H$.

Applying $\Psi^{-1}$ and taking cohomology gives an exact sequence
\[
 E_{\mathrm{low}}\longrightarrow\Psi^{-1}Q
 \longrightarrow H^1_{\mathcal H}(\Psi^{-1}L)\longrightarrow0.
\]
The image of the first map belongs to $\mathcal T$, being a
quotient of $E_{\mathrm{low}}$, so its imaginary charge is
nonnegative by Lemma~\ref{lem:all-degree-stability-function}.
The last term also has nonnegative imaginary charge.  Indeed,
its vertices come from the HN factors of $L$ with phases in
$(0,1/d]$: applying $\Psi^{-1}$ and then $[1]$ puts these phases
in $(t,1]$.  Its source contributions therefore have nonnegative
imaginary parts, and multiplication by $q=e^{-\pi it}$ puts its
target contributions in phases $(0,1/d]$.
Additivity now gives $\operatorname{Im}\mathcal Z(\Psi^{-1}Q)\ge0$.
But~\eqref{eq:all-degree-Psi-charge-rotation} gives
\[
 \mathcal Z(\Psi^{-1}Q)=e^{-\pi i/d}\mathcal Z(Q),
\]
whose phase lies in $(-1,-1/d]$
by~\eqref{eq:all-degree-torsion-mass}, since $0\ne Q\in\mathcal F$.
Its imaginary part is strictly negative, a contradiction.
Thus $H^0_{\mathcal H}(\Psi E)\in\mathcal T$, as required.
\end{proof}

\begin{lemma}\label{lem:all-degree-Psi-gepner}
If $\sigma$ is of Gepner type $(\Psi,1/d)$ for an exact enhanced
autoequivalence $\Psi$ of $\mathcal C$, then $\Sigma$ is of
Gepner type $(\Psi,1/d)$ for the induced autoequivalence of
$\mathcal C\otimes\widehat{\mathcal Q}_d$.
\end{lemma}

\begin{proof}
The functor $\Psi^d[-1]$ preserves every slice $\mathcal P(\phi)$
and the central charge $Z$ of $\sigma$ on $\mathcal C$.
Applied at every vertex, it therefore preserves
$\mathcal H$, $\mathcal Z$, $p$, the torsion pair, and $U$.
Thus $\Psi^dU=U[1]$, and Lemma~\ref{lem:all-degree-Psi-aisle} gives
\[
 U[1]=\Psi^dU\subset\Psi U\subset U,
 \qquad U\subset\Psi^{-1}U\subset U[-1].
\]
The second chain shows that, for $E\in\mathcal A'$, the object
$\Psi^{-1}E$ has $\mathcal A'$-cohomology only in degrees $0$
and $1$.  Its cohomology triangle expresses it as an extension
of an object of $\mathcal A'[-1]$ by one of $\mathcal A'$.
Its HN phases for $\Sigma$ consequently lie in $(-1,1]$.

To finish, transport the slicing by $\Psi$ and record the
expected phase shift:
\[
 \mathcal P'(\phi)=\Psi\mathcal P_\Sigma(\phi-1/d).
\]
By~\eqref{eq:all-degree-Psi-charge-rotation}, this is a slicing
with central charge $\mathcal Z$; hence
$\Sigma'=(\mathcal Z,\mathcal P')$ is a stability condition with
the same charge as $\Sigma$.
The $\mathcal P'$-phases of $E$ are the $\mathcal P_\Sigma$-phases
of $\Psi^{-1}E$, increased by $1/d$.  The preceding bound gives
\[
 \mathcal A'\subset\mathcal P'((-1+1/d,1+1/d])
              \subset\mathcal P'((-1,2]).
\]
Collins--Polishchuk's uniqueness criterion
\cite[Proposition~4.2(a)]{CP} says that two stability conditions
with the same charge coincide if the heart of the first lies
in phases $(-1,2]$ of the second.  It follows that
$\Sigma'=\Sigma$, so
\[
 \Psi\mathcal P_\Sigma(\phi)
       =\mathcal P_\Sigma(\phi+1/d).
\]
Together with~\eqref{eq:all-degree-Psi-charge-rotation}, this is
the Gepner property $(\Psi,1/d)$.
\end{proof}

This proves Theorem~\ref{thm:all-degree-coefficients}.

\section{Variations on the Thom--Sebastiani theorem}
\label{sec:thom-sebastiani}

We now transfer the stability conditions of
Section~\ref{sec:all-degree} to Fermat hypersurfaces.

\begin{theorem}\label{thm:all-degree-fermat}
For every $d\ge2$ and $M\ge1$, give each variable weight $1$ and
put $$F_{d,M}=x_1^d+\cdots+x_M^d.$$
The category $\MF^{\gr}(F_{d,M})$ has a numerical
stability condition of Gepner type $(\tau,2/d)$ satisfying a total
mass--Hom bound and the support property on
$K_{\num}(\MF^{\gr}(F_{d,M}))$.
Here $\tau$ is the internal grading shift.
\end{theorem}

We first relate matrix factorizations of a sum of polynomials to the tensor product of
the factorization categories of its summands.  We then construct
compatible stability conditions on the tensor product and prove
the theorem at the end of the section.

Let $f_1,\ldots,f_r$ be weighted homogeneous polynomials of the same degree $d$
in disjoint sets of positively weighted variables, each with an isolated
critical point at the origin.  We call each set of variables a
\emph{block}.  For Fermat hypersurfaces each block consists of one
variable with potential~$x^d$. In full generality, 
put $W=f_1+\cdots+f_r$ and define
\[
 L=\frac{\bigoplus_{i=1}^r\Z e_i}
          {\langle de_i-de_j\ :\ 1\leq i,j\leq r\rangle},
 \qquad \pi:L\longrightarrow\Z,\quad e_i\longmapsto1.
\]
A variable of weight $a$ in block $i$ has $L$-degree $ae_i$,
and $W$ has the common degree $de_i$.  Choosing the last block gives a splitting
\[
 L=\Z e_r\oplus H,\qquad
 H=\ker\pi=\bigoplus_{i=1}^{r-1}(\Z/d\Z)(e_i-e_r).
\]

\begin{proposition}
\label{prop:finite-grading-correction}
Set $\mathcal T=\bigotimes_i\MF^{\Z}(f_i)$ and
$\mathcal D=\MF^{\Z}(W)$, and let $\tau_i$ be the grading shift
in the $i$th factor of $\mathcal T$.
There is a coherent $H$-action on $\mathcal T$ in which
$e_i-e_r$ acts as $\tau_i\tau_r^{-1}$, and an enhanced equivalence
\[
 \mathcal D\simeq\mathcal T^H.
\]
Under this equivalence, the ordinary grading shift on
$\mathcal D$ is induced by $\tau_r$.
\end{proposition}

\begin{proof}
There are two steps.  The Thom--Sebastiani theorem identifies
$\mathcal T$ with the category of matrix factorizations of $W$ carrying the block grading
$L$.  We then recover ordinary integer grading by identifying the
components whose degrees differ by an element of $H$.

\smallskip
\noindent\emph{The block grading and the Thom--Sebastiani theorem.}
We recall the part of \cite[Theorem~1.1]{BFK} that we use.
Suppose that two polynomials $u$ and $v$ on separate affine spaces
are multiplied by characters $\chi:G_1\to\C^*$ and
$\psi:G_2\to\C^*$ of diagonalizable groups.  Equivariant
modules decompose into character spaces, so equivariant
factorizations are simply factorizations graded by the character
group.  The theorem identifies the Morita tensor product of their equivariant
factorization categories with the category of factorizations of the polynomial
$u+v$ equivariant for
\[
 G_1\times_{\C^*}G_2
   =\{(g_1,g_2):\chi(g_1)=\psi(g_2)\}.
\]
The equivalence is enhanced. The hypotheses are that the character
$(g_1,g_2)\mapsto\chi(g_1)\psi(g_2)^{-1}$ has infinite order and
that every singular point of $\{u+v=0\}$ lies in
$\{u=0\}\times\{v=0\}$.  

For $r\ge2$,
apply the theorem first to $f_1$ and $f_2$, and then adjoin one
polynomial at a time.  At the step adjoining $f_k$, the group
already satisfies $t_1^d=\cdots=t_{k-1}^d$.  Its product with the
new $\C^*$-factor has character quotient $t_1^dt_k^{-d}$, of
infinite order since $t_k$ varies freely.  The fiber product adds
exactly the relation $t_k^d=t_1^d$.  The resulting group is therefore
\[
 G=\{(t_1,\ldots,t_r)\in(\C^*)^r:t_1^d=\cdots=t_r^d\}
   =\Hom(L,\C^*).
\]
Its character group is $L$, so $G$-equivariance is precisely
$L$-grading.  A variable of weight
$a$ in block $i$ has character $ae_i$, and $W$ has character
$de_i=de_r$.  The resulting equivalence is thus
$\mathcal T\simeq\MF^L(W)$.  It takes the grading shift $\tau_i$
to the shift by $e_i$.  Consequently the shift by $e_i-e_r$
corresponds to $\tau_i\tau_r^{-1}$.  For any $L$-graded matrix
factorization $N$, the grading shifts satisfy
$N(h)(h')=N(h+h')$ for $h,h'\in H$.
These identities give  the required coherent action on $\mathcal T$.

Write $S$ for the polynomial ring in all the variables and use
$L=\Z e_r\oplus H$.  Starting with an ordinarily graded
factorization $M$, place a copy of $M_n$ in every degree over $n$,
$
 J(M)_{ne_r+h}=M_n$ for $n\in\Z$, $h\in H$.
If $s\in S$ has $L$-degree $ae_r+b$, multiplication by $s$ sends
the copy of $M_n$ indexed by $h$ to the copy of $M_{n+a}$ indexed
by $h+b$, using the original multiplication on $M$.  Use the
original factorization maps on these copies as well; their
compositions are multiplication by $W$, of $L$-degree $de_r$.
The identity maps between the copies give isomorphisms
$\rho_h:J(M)(h)\to J(M)$ satisfying the equivariance identities.
Moreover, a free summand $S(n)$ gives
$\bigoplus_{h\in H}S(ne_r+h)$, so this construction has finite-rank
free terms.

Conversely, let $(N,\rho)$ be an $L$-graded factorization with
isomorphisms $\rho_h:N(h)\to N$ satisfying, before passage to
homotopy,
$
 \rho_0=\operatorname{id}$,
 $\rho_{h+h'}=\rho_h\circ\rho_{h'}(h)$.
Set $M_n=N_{ne_r}$.  For
$s\in S_{ae_r+b}$, define multiplication on $M$ by
$
 M_n=N_{ne_r}\xrightarrow{\ s\ }
 N_{(n+a)e_r+b}\xrightarrow{\ \rho_b\ }
 N_{(n+a)e_r}=M_{n+a}$.
The equivariance identities make this multiplication associative.
The factorization maps restrict to $M$, since their degrees have
zero $H$-component, and their compositions remain multiplication by
$W$.  These constructions are inverse: taking the components with
$H$-index zero recovers $M$ from $J(M)$, while the maps
$\rho_h:N_{ne_r+h}\to N_{ne_r}$ identify $N$ with $J(M)$.

The same constructions apply to morphisms and commute with their
differentials, so they identify the enhanced categories after
perfect completion.  By \cite[Theorem~8.7]{Elagin}, this completion
of the strict equivariant model computes the equivariant objects
in the homotopy category.  We obtain
\[
 \mathcal D\simeq\MF^L(W)^H\simeq\mathcal T^H.
\]
Finally, increasing the integer grading by one increases the
$L$-grading by $e_r$: explicitly, $J(M(1))=J(M)(e_r)$, with the
same equivariant structure.  This proves the assertion about
$\tau_r$.
\end{proof}

\begin{lemma}\label{lem:equivariant-support}
Let a finite group $H$ act coherently on an enhanced triangulated category
$\mathcal T$, and preserve a stability condition
$\sigma=(Z,\mathcal P)$ supported on
$v:K_0(\mathcal T)\to\Lambda$.  Suppose $\mathcal T^H$ has its enhanced
equivariant triangulated structure.  With $U:\mathcal T^H\to\mathcal T$
the forgetful functor, the formulas
\[
 Z^H=Z\circ U_*,\qquad
 \mathcal P^H(\phi)=\{E:U(E)\in\mathcal P(\phi)\}
\]
define a stability condition supported on the finite-rank lattice
$
           \Lambda^H=\im(v\circ U_*)\subseteq\Lambda
$.
Furthermore, if $v$ is numerical, then so is $v\circ U_*$.
Moreover, $m_{\sigma^H}(E)=m_\sigma(U(E))$.
If $\sigma$ has a total mass--Hom bound, then so does $\sigma^H$.
\end{lemma}

\begin{proof}
The forgetful functor is exact, faithful, and conservative.
The shift and orthogonality axioms follow from those in
$\mathcal T$.  The canonical HN truncations of $U(E)$ inherit
coherent equivariant structures, as explained in
Section~\ref{sec:background}, giving the HN filtration in $\mathcal T^H$.
A nonzero equivariant semistable object has a nonzero underlying
object, hence a nonzero central charge.

Choose a norm on $\Lambda_{\R}$ and a constant $C$ satisfying
the support inequality for $\sigma$.
For an equivariant semistable object $E$, its underlying object
$U(E)$ is semistable, so
\[
 \|v(U(E))\|\le C|Z(U(E))|=C|Z^H(E)|.
\]
Restricting the norm to $\Lambda^H_{\R}$ proves the support property
on the stated image lattice.

The HN construction also proves $m_{\sigma^H}(E)=m_\sigma(U(E))$.
Induction $I(A)=\bigoplus_{h\in H}F_h(A)$, with its canonical
equivariant structure, is both left and right adjoint to $U$
\cite[Lemma~3.8]{Elagin}.  Hence
$
 \chi_{\mathcal T}(U_*x,A)=\chi_{\mathcal T^H}(x,I(A))$.
If $x$ is numerically trivial, the right-hand side vanishes for
every $A\in\mathcal T$.  Thus $U_*x$ is numerically trivial, and
$v\circ U_*$ factors through the numerical Grothendieck group of $\mathcal T^H$.

Finally, faithfulness of $U$ in every degree gives, for any fixed
$A\in\mathcal T^H$,
\[
 \begin{aligned}
 \sum_{n\in\Z}\dim_{\C}\Hom(A,E[n])
 &\le\sum_{n\in\Z}\dim_{\C}\Hom(U(A),U(E)[n])\\
 &\le C_{U(A)}m_\sigma(U(E))=C_{U(A)}m_{\sigma^H}(E).
 \end{aligned}
\]
This proves the total mass--Hom bound.
\end{proof}

\begin{corollary}\label{cor:gepner-thom-sebastiani}
With the notation of Proposition~\ref{prop:finite-grading-correction},
suppose the  category $\mathcal T$ carries a stability condition
$\sigma_{\mathcal T}=(Z_{\mathcal T},\mathcal P_{\mathcal T})$
supported on a finite-rank numerical quotient
$v:K_0(\mathcal T)\to\Lambda$.  Assume that
$\sigma_{\mathcal T}$ is of Gepner type $(\tau_i,2/d)$ for every
internal grading shift $\tau_i$.
Let $U:\mathcal D\simeq\mathcal T^H\to\mathcal T$ be the forgetful
functor.  Then
\[
 Z_{\mathcal D}=Z_{\mathcal T}\circ U_*,\qquad
 \mathcal P_{\mathcal D}(\phi)
       =\{E:U(E)\in\mathcal P_{\mathcal T}(\phi)\}
\]
define a numerical stability condition on
$\mathcal D=\MF^{\Z}(W)$, supported on
$\Lambda_{\mathcal D}=\im(v\circ U_*)\subseteq\Lambda$.
This condition is of Gepner type $(\tau,2/d)$, where $\tau$ is
the ordinary internal grading shift.
A total mass--Hom bound for $\sigma_{\mathcal T}$ is preserved.
\end{corollary}

\begin{proof}
The shifts $\tau_i$ have the same effect on the charge and phases,
so these effects cancel in $\tau_i\tau_r^{-1}$.  The generators
of the $H$-action therefore preserve $\sigma_{\mathcal T}$.
Lemma~\ref{lem:equivariant-support} gives the stated condition,
including its support lattice, numericality, and the mass--Hom bound.
By Proposition~\ref{prop:finite-grading-correction}, the ordinary
grading shift is induced by $\tau_r$, so $U\tau\simeq\tau_rU$.
The formulas for the induced charge and slicing therefore show
that the induced condition is of Gepner type $(\tau,2/d)$.
\end{proof}

To apply Corollary~\ref{cor:gepner-thom-sebastiani}, we need a
stability condition on the tensor product that is of Gepner type
$(\tau_i,2/d)$ for every factor.
Theorem~\ref{thm:all-degree-coefficients} lets us add one quiver
tensor factor at a time while preserving these Gepner properties.
Fix $d\ge2$, and write $\mathcal Q_d$ for the derived category of
the alternating quiver $Q_d$.  Recall the weights
\begin{equation}\label{eq:all-degree-weights}
 r_j=\frac{\sin(j\pi/d)}{\sin(\pi/d)},\qquad
 q=-e^{\pi i/d},\qquad
 a_j=\begin{cases}r_j,&j\text{ odd},\\q r_j,&j\text{ even}.
 \end{cases}
\end{equation}

\begin{proposition}\label{prop:all-degree-tensor}
For every $d\ge2$ and $M\ge1$, the category
$\mathcal D_{d,M}=\mathcal Q_d^{\otimes M}$ has a numerical
stability condition
$\Sigma_{d,M}=(Z_{d,M},\mathcal P_{d,M})$ supported on
$K_{\num}(\mathcal D_{d,M})\simeq\Z^{(d-1)^M}$.
It also satisfies a total mass--Hom bound.
Let $\tau_i$ act by the inverse Serre functor followed by shift,
$S_{Q_d}^{-1}[1]$, in the $i$th factor.  Then $\Sigma_{d,M}$ is
of Gepner type $(\tau_i,2/d)$ for every $i$.
The central charge is
\begin{equation}\label{eq:all-degree-charge}
 Z_{d,M}(E)=-\sum_{\boldsymbol\alpha\in\{1,\ldots,d-1\}^{M}}
       \left(\prod_{i=1}^{M}a_{\alpha_i}\right)
                    n_{\boldsymbol\alpha}(E),
\end{equation}
where $\boldsymbol\alpha$ chooses one vertex in each factor and
$n_{\boldsymbol\alpha}(E)=\sum_k(-1)^k\dim H^k(E_{\boldsymbol\alpha})$
is the Euler characteristic of the complex at that vertex.
For $d=2$, we have $\mathcal Q_2=\Perf(\C)$ and $a_1=1$.
\end{proposition}

\begin{proof}
For $d=2$, the tensor category is $\Perf(\C)$ and each $\tau_i$
is $[1]$.  Use the standard stability condition, with charge
$-\dim$ on vector spaces and phase $1$.  Its generator satisfies
\[
 \sum_{n\in\Z}\dim_{\C}\Hom(\C,E[n])
 =\sum_{n\in\Z}\dim_{\C}H^n(E)=m(E),
\]
giving the mass--Hom bound.

We now assume $d\ge3$.
Begin with this same condition on $\Perf(\C)$ and apply
Theorem~\ref{thm:all-degree-coefficients} $M$ times.  At each step,
the coefficient category is the tensor product with fewer factors,
and the new factor is the doubled quiver category
$\widehat{\mathcal Q}_d$.  These categories are smooth and proper,
as required by the theorem.  We obtain a stability condition on
$\widehat{\mathcal Q}_d^{\otimes M}$ satisfying the support property
and a total mass--Hom bound, since the coefficient theorem
preserves both properties at each step.
At each step, the new condition is of Gepner type $(R_i,1/d)$
for the reflection functor of the new factor.  Compatibility with
coefficient autoequivalences preserves the Gepner properties for
the reflection functors of all previous factors.

We next choose the orientation $Q_d$ in each factor.
Since $\widehat{\mathcal Q}_d$ is the direct sum of the categories
for $Q_d$ and $Q_d^{\rm op}$, its $M$th tensor power has $2^M$
orthogonal components.  Projection to a component computes direct
summands, so it preserves the aisle and coaisle of every phase-cut
t-structure.  The HN filtrations therefore restrict to each
component.  The support estimate restricts to its vertex lattice.
The mass--Hom bound restricts as well: project a split generator
to the chosen orthogonal component.  It still generates that
component, and all morphisms from the other components vanish.
In the component $\mathcal D_{d,M}$ where all orientations are
$Q_d$, the square $R_i^2$ preserves the component and equals
$S_{Q_d}^{-1}[1]$ on its $i$th factor by
Lemma~\ref{lem:all-degree-reflection-periodicity}.
Since $\tau_i=R_i^2$, the restricted condition is of Gepner type
$(\tau_i,2/d)$ for every $i$.

It remains to identify the charge and the lattice.
Adding a factor multiplies the coefficient at its $j$th vertex
by $a_j$.  Starting with the charge $-\dim$ therefore gives the
product of weights in~\eqref{eq:all-degree-charge}, with a single
initial minus sign.  The vertex classes identify
$K_0(\mathcal D_{d,M})$ with $\Z^{(d-1)^M}$.
In the basis of simple vertex objects, the Euler matrix is the tensor
product of the Euler matrices of $Q_d$.  Order the vertices so
that every arrow points forward.  Each matrix is triangular with
diagonal entries $1$, so the tensor product is nondegenerate.
Thus the Euler pairing has zero radical and this vertex lattice
is the full numerical Grothendieck group.
\end{proof}

\begin{proof}[Proof of Theorem~\ref{thm:all-degree-fermat}]
By Proposition~\ref{prop:one-variable-model}, the enhanced
equivalence $\MF^{\gr}(x^d)\simeq\mathcal Q_d$ identifies the
internal grading shift with $S_{Q_d}^{-1}[1]$.
Proposition~\ref{prop:all-degree-tensor} gives a numerical
stability condition on the tensor product that satisfies a total
mass--Hom bound and is of Gepner type $(\tau_i,2/d)$ for every $i$.
Corollary~\ref{cor:gepner-thom-sebastiani}, applied to $f_i=x_i^d$,
transfers both properties to $\MF^{\gr}(F_{d,M})$.
Its numerical Grothendieck group is finite rank by the one-variable
model when $M=1$, and by Orlov's semiorthogonal decompositions
when $M\ge2$.  Lemma~\ref{lem:mass-hom-generator}
therefore gives support on the full numerical lattice.
\end{proof}

\section{Restriction to weighted hyperplanes}
\label{sec:hyperplane}

We now pass from Fermat hypersurfaces to their smooth linear
sections. Derived restriction and extension of scalars give adjoint
functors between the corresponding singularity categories.

\begin{lemma}
\label{lem:dg-induction}
Let $\cC$ and $\cD$ be smooth proper dg categories, taken to be
pretriangulated and idempotent complete, and let
\[
 F:\cC\longrightarrow\cD,\qquad R:\cD\longrightarrow\cC
\]
be exact dg functors with $F\dashv R$.  Suppose that the image of
$F$ split-generates $\cD$.
Let
$\sigma=(Z,\mathcal P)$ be a stability condition on $\cC$, satisfying the
support property relative to the class map $v:\,K_0(\cC)\to\Lambda$, where $\Lambda$ is a
finite-rank free abelian group.  Assume
\begin{equation}
 (R\circ F)\bigl(\mathcal P(\phi)\bigr)
       \subseteq\mathcal P([\phi,\infty))
       \qquad\text{for every }\phi\in\R.
 \label{eq:induction-phase}
\end{equation}
Then we have a Bridgeland stability condition on $\cD$ defined by 
\begin{align}
 \mathcal P_{\cD}(\phi)
   &=\{E\in\cD:R(E)\in\mathcal P(\phi)\},
       \label{eq:induced-slicing}\\
 Z_{\cD}(E)&=Z(R(E)).
       \label{eq:induced-charge}
\end{align}
It satisfies the support property relative to the class map
\[
 v_{\cD}=v\circ R_*:K_0(\cD)\longrightarrow
 \Lambda_{\cD}:=\im(v\circ R_*)\subseteq\Lambda.
\]
Moreover, $R$ preserves Harder--Narasimhan filtrations and phase intervals: for every interval $I\subseteq\R$,
$E\in\mathcal P_{\cD}(I)$ if and only if $R(E)\in\mathcal P(I)$.
Writing $\sigma_{\cD}=(Z_{\cD},\mathcal P_{\cD})$, we have
$m_{\sigma_{\cD}}(E)=m_\sigma(R(E))$.
If $\sigma$ has a total mass--Hom bound, then so does $\sigma_{\cD}$.
\end{lemma}

\begin{proof}
We adapt Polishchuk's induction argument
\cite[Theorem~2.1.2]{Polishchuk}; compare
\cite[Theorem~A.17 and Corollary~A.23]{Families}.
First extend the adjunction to
\[
 F:\Ind(\cC)\longrightarrow\Ind(\cD),\qquad
 R:\Ind(\cD)\longrightarrow\Ind(\cC).
\]
Here $\Ind(\cC)$ and $\Ind(\cD)$ are the derived categories of
right dg modules over $\cC$ and $\cD$, as in
Section~\ref{subsec:categorical-tensor-products}.
The extended right adjoint is particularly concrete: it restricts
a dg module along $F$:
\[
                         (RM)(A)=M(F(A)).
\]
Restriction preserves quasi-isomorphisms and arbitrary coproducts
because both are computed at each object.  Its left adjoint is
derived extension of scalars, which sends the representable module
$h_A$ to $h_{F(A)}$ and hence extends the original $F$
\cite[Section~3.8]{Keller}.
For $B\in\cD$, the original adjunction gives
\[
 (Rh_B)(A)=\operatorname{Hom}_{\cD}^{\bullet}(F(A),B)
       \simeq\operatorname{Hom}_{\cC}^{\bullet}(A,R(B)).
\]
Thus the restricted module $Rh_B$ represents the original object
$R(B)$.  This justifies using the same letters for the two adjunctions.
The extended $R$ is conservative: if $R(E)=0$, then
\[
 \Hom(F(A),E[k])=\Hom(A,R(E)[k])=0
 \qquad(A\in\cC,\ k\in\Z).
\]
Since $F(\cC)$ split-generates $\cD$, the same vanishing holds
with any object of $\cD$ in place of $F(A)$.  All values of the
dg module $E$ are therefore acyclic, so $E=0$.
We also claim that $E$ is compact whenever $R(E)$ is compact.
Indeed, for every $A\in\cC$, adjunction gives
\[
 \bigoplus_k\Hom(F(A),E[k])
   \simeq\bigoplus_k\Hom(A,R(E)[k]).
\]
The right-hand side is finite dimensional by properness of $\cC$.
This finiteness is preserved under shifts, cones, finite sums,
and direct summands in the first argument.  It therefore holds
for every object of $\cD$ in place of $F(A)$.
We use the following criterion: a dg module over a smooth proper
category is perfect if and only if its value at every object has
finite-dimensional total cohomology
\cite[Lemma~2.8(3)]{ToenVaquie}.
This criterion  makes $E$ perfect, hence compact.

Fix $a\in\R$.  Cutting the slicing at phase $a$ gives a bounded
t-structure on $\cC$, with aisle $\mathcal P((a,\infty))$ and
right orthogonal $\mathcal P((-\infty,a])$.
We use the following existence theorem: in a triangulated category
with coproducts, the smallest full subcategory containing a set
of compact objects and closed under coproducts, extensions, and
translation $[1]$ is an aisle~\cite[Theorem~A.1]{AJS}.
Choose representatives of the isomorphism classes in
$\mathcal P((a,\infty))$.  They are compact in $\Ind(\cC)$,
and their images under $F$ are compact in $\Ind(\cD)$.
The theorem therefore gives two aisles:
\[
 \begin{aligned}
 U_a&\subset\Ind(\cC),
       &&\text{generated by }\mathcal P((a,\infty)),\\
 V_a&\subset\Ind(\cD),
       &&\text{generated by }F\bigl(\mathcal P((a,\infty))\bigr).
 \end{aligned}
\]
In both cases, ``generated'' means closure under the three
operations just listed.  The t-structure with aisle $U_a$ extends
the original phase cut on $\cC$.  Indeed, its original truncation
triangles still have first terms in $U_a$ and last terms in
$U_a^\perp$: orthogonality to the generators persists under
coproducts, extensions, and positive shifts in the first argument.
Thus these triangles remain truncation triangles in $\Ind(\cC)$;
compare~\cite[Lemma~2.1.1]{Polishchuk}.

Adjunction identifies the other right orthogonal explicitly:
\[
 V_a^\perp
 =\{E:\Hom(F(A),E)=0\text{ for all }A\in\mathcal P((a,\infty))\}
 =R^{-1}(U_a^\perp).
\]
Every object of $\mathcal P((a,\infty))$ is a finite extension of
semistable objects of phases greater than $a$.  Hence
\eqref{eq:induction-phase} gives
$RF(\mathcal P((a,\infty)))\subseteq\mathcal P((a,\infty))$.
Since $R$ preserves coproducts, extensions, and shifts, it follows
that $R(V_a)\subseteq U_a$.
Applying $R$ to a truncation triangle for $V_a$ therefore gives
a triangle whose first term lies in $U_a$ and whose last term
lies in $U_a^\perp$.  By uniqueness, this is the truncation
triangle of the image under $R$.  Thus $R$ commutes with these
truncations and is t-exact.

For $E\in\cD$, the object $R(E)$ belongs to $\cC$, and its
truncations for $U_a$ also belong to $\cC$.
By t-exactness, these are the images under $R$ of the truncations
of $E$ for $V_a$.  As proved above, compactness of $R(M)$ implies
compactness of $M$.  
Smoothness and properness ensure that the truncations we
construct in $\Ind(\cD)$ belong to $\cD$.
Thus both truncations of $E$ belong to
$\cD$.  The t-structure consequently restricts
from $\Ind(\cD)$ to $\cD$.
Its aisle and right orthogonal are
\[
 \begin{aligned}
 V_a\cap\cD&=\{E\in\cD:R(E)\in\mathcal P((a,\infty))\},\\
 V_a^\perp\cap\cD&=\{E\in\cD:R(E)\in\mathcal P((-\infty,a])\}.
 \end{aligned}
\]
For the first equality, if $R(E)$ belongs to the indicated aisle,
the last term of the truncation triangle of $E$ has zero image
under $R$ and hence vanishes.  The second equality follows from
adjunction above.  These formulas also prove boundedness, because
the phase-cut t-structure on $\cC$ is bounded.
The same construction works with the aisle
$\mathcal P([a,\infty))$ and right orthogonal
$\mathcal P((-\infty,a))$.  We can therefore cut the phases on
$\cD$ either above or at $a$.  These cuts are compatible, since
their aisles and right orthogonals are inverse images under $R$.
The shift rule in \eqref{eq:induced-slicing} follows because $R$
commutes with shifts.  To prove orthogonality, let
$E\in\mathcal P_{\cD}(\phi)$ and $E'\in\mathcal P_{\cD}(\psi)$
with $\phi>\psi$.  Choose $a$ with $\phi>a\geq\psi$.
Then $E$ belongs to the induced aisle for this cut and $E'$ to
its right orthogonal, so $\Hom(E,E')=0$.

For $0\ne E\in\cD$, let $\phi_1>\cdots>\phi_s$ be the HN phases
of $R(E)$.  Truncate $E$ at $[\phi_1,\infty)$.
Because $R$ commutes with this truncation, the first term maps
to the first HN factor of $R(E)$ and hence belongs to
$\mathcal P_{\cD}(\phi_1)$.  Apply the next cut,
$[\phi_2,\infty)$, to the remaining object, and continue.
At step~$j$, the removed term  lies in $\mathcal P_{\cD}(\phi_j)$
and maps to the corresponding HN factor of $R(E)$.
After $s$ steps, the remainder has zero image under $R$ and
therefore vanishes.  These triangles give the HN filtration of
$E$.  In particular, $E$ and $R(E)$ have exactly the same HN
phases, which proves the assertion about phase intervals.

If $0\ne E\in\mathcal P_{\cD}(\phi)$, then $R(E)$ is nonzero
and semistable of phase $\phi$.  Thus
$Z_{\cD}(E)=Z(R(E))\in\R_{>0}e^{i\pi\phi}$.
Choose a norm and support constant $C$ for the original stability condition.
Then
\[
 \lVert v_{\cD}(E)\rVert
 =\lVert v(R(E))\rVert
 \leq C|Z(R(E))|
 =C|Z_{\cD}(E)|,
\]
proving the support inequality on the finite-rank image lattice
$\Lambda_{\cD}$.

Preservation of HN filtrations and equality of charges give
$m_{\sigma_{\cD}}(E)=m_\sigma(R(E))$.
Let $G$ split-generate $\cC$ and satisfy
\[
 \sum_{n\in\Z}\dim_{\C}\Hom(G,A[n])\le C_Gm_\sigma(A)
 \qquad\text{for every }A\in\cC.
\]
Since $F(\cC)$ split-generates $\cD$, the object $F(G)$ is a
split generator.  Adjunction gives
\[
 \begin{aligned}
 \sum_{n\in\Z}\dim_{\C}\Hom(F(G),E[n])
 &=\sum_{n\in\Z}\dim_{\C}\Hom(G,R(E)[n])\\
 &\le C_Gm_\sigma(R(E))=C_Gm_{\sigma_{\cD}}(E),
 \end{aligned}
\]
and Lemma~\ref{lem:mass-hom-generator} gives the total mass--Hom bound.
\end{proof}

Let $S=\C[x_1,\ldots,x_N]$ with variables having positive weights, and let $W\in S$
be a weighted homogeneous polynomial of degree $d$, with an isolated critical
point at the origin.  Choose a coordinate $z$ of weight $a$ and put
\[
 A=S/(W),\qquad B=A/(z),\qquad w=W|_{z=0}.
\]
Assume that $w$ is nonzero and still has an isolated critical point
at the origin.  The condition $w\ne0$ means that $z$ does not divide
$W$, so multiplication by $z$ is injective on $A$:
\begin{equation}
 0\longrightarrow A(-a)\xrightarrow{z}A
   \longrightarrow B\longrightarrow0.
 \label{eq:coordinate-resolution}
\end{equation}

Since the origin is a critical point of $W$, all its first partial
derivatives vanish there.  Thus $W\in(x_1,\ldots,x_N)^2$, and
its restriction $w$ also has no linear terms in the remaining
variables.  The coordinates therefore give bases of the respective
homogeneous maximal ideals modulo their squares.  We will use
this minimality in the residue-field calculation below.
Set
\[
 \cC=\Dsg^{\mathrm{gr}}(A)\simeq\MF^{\mathrm{gr}}(W),
 \qquad
 \cD=\Dsg^{\mathrm{gr}}(B)\simeq\MF^{\mathrm{gr}}(w).
\]
The weighted form
of Orlov's theorem shows that these categories are smooth and proper:
they are either admissible components of the derived categories
of smooth proper hypersurface Deligne--Mumford stacks
or are obtained
from those derived categories by adjoining finitely many
exceptional objects.  See
\cite[Theorem~2.5 and Remark~3.12]{Orlov} and
\cite[Theorem~6.13]{BFKHodge}.
Denote internal grading shift by $\tau$; in either category $\cC$ or $\cD$,
$\tau^d\simeq[2]$.

Derived extension of scalars $F=B\otimes_A^{\mathbf L}-$ is left
adjoint to restriction of scalars $R:\,\cD\to\cC$.
To see that these functors preserve the bounded derived categories
and descend to the singularity categories, use the exact sequence \eqref{eq:coordinate-resolution}.
This finite free resolution of $B$ over $A$ shows that $F$
preserves boundedness and that $R$ sends perfect $B$-complexes
to perfect $A$-complexes.  Restriction also preserves bounded
complexes of finitely generated modules, since $B$ is finite over
$A$.  Finally, $F$ sends graded free $A$-modules to graded free
$B$-modules and hence preserves perfect complexes.
The adjunction therefore descends to
\[
 F:\cC\longrightarrow\cD,\qquad
 R:\cD\longrightarrow\cC,\qquad F\dashv R
\]
by~\cite[Lemma~1.1]{Orlov}, and extends to their idempotent
completions.  Both functors commute with internal grading shifts.

\begin{lemma}
\label{lem:residue-field-base-change}
Let $k_A$ and $k_B$ be the residue fields, placed in internal
degree zero.  Then
\[
 F(k_A)\simeq k_B\oplus k_B(-a)[1]
 \qquad\text{in}\quad \Db_{\gr}(B).
\]
Here $\Db_{\gr}(B)$ is the bounded derived category of finitely
generated graded $B$-modules, with module maps preserving the
internal degree.  Thus the splitting holds before taking the
quotient by perfect complexes.  Moreover, $F(\cC)$ split-generates $\cD$.
\end{lemma}

\begin{proof}
Resolution~\eqref{eq:coordinate-resolution} computes the cohomology of
$F(k_A)$: it is $k_B$ in degree zero and $k_B(-a)$ in degree $-1$.
The obstruction to splitting these  cohomology objects is the
connecting morphism in
\begin{equation}
 k_B(-a)[1]\longrightarrow F(k_A)\longrightarrow k_B
   \xrightarrow{e}k_B(-a)[2].
 \label{eq:residue-triangle}
\end{equation}
We show that $e$ vanishes by applying $\Hom(-,k_B(-a)[1])$ and using adjunction.  The resulting exact
sequence contains
\begin{equation}
\ 0\longrightarrow\Ext^1_{\mathrm{gr},B}(k_B,k_B(-a))
 \longrightarrow\Ext^1_{\mathrm{gr},A}(k_A,k_A(-a))
 \longrightarrow\C
 \xrightarrow{\ e\ }\Ext^2_{\mathrm{gr},B}(k_B,k_B(-a)).
\label{eq:change-rings-ext}
\end{equation}
The last map sends $1$ to $e$, up to sign.
For either ring, the degree-$1$ Ext group is dual to the weight-$a$
part of the homogeneous maximal ideal modulo its square.  Its dimension
therefore counts minimal generators of weight $a$.  The hypotheses on
$W$ and $w$ ensure that the coordinates give these generators.
Passing from $A$ to $B$ removes precisely the generator $z$ of
weight $a$.  Thus the middle $\Ext^1$ group has dimension one
larger than the first.  Exactness forces the map to $\C$ to be
surjective, and therefore $e=0$.  This proves the splitting.

The graded residue fields $k_B(j)$ split-generate $\cD$.
For $\dim B>0$, the isolated-singularity assumption implies that
$B$ is reduced, and generation follows from
\cite[Corollary~4.14]{BFKHodge}, applied to the grading action
of $\mathbf G_m$ on $\operatorname{Spec}B$.
In dimension zero, it follows from finite-length filtrations.
\footnote{
We recall a direct argument that covers both cases.
For a matrix factorization of $w$ with differential $\delta$,
differentiating $\delta^2=w\operatorname{id}$ gives
$
 \delta(\partial_i\delta)+(\partial_i\delta)\delta
   =(\partial_iw)\operatorname{id}$.
Thus multiplication by each partial derivative is nullhomotopic,
with the appropriate grading shift.  Because the critical point
is isolated, a sufficiently high power of each coordinate belongs
to the Jacobian ideal.
Now represent an object of $\Dsg^{\mathrm{gr}}(B)$ by a bounded
complex of finitely generated graded modules and tensor it with
the Koszul complex of these powers.  In the singularity category,
each step takes the cone of a zero multiplication map, so the
original object is a direct summand of the result.  The resulting
complex has finite-length cohomology supported at the origin;
its cohomology modules have filtrations with factors $k_B(j)$.
It is therefore built from the graded residue fields.
The same conclusion holds after idempotent completion.}
Since $(d)=[2]$,
the finitely many objects $k_B(j)$ with $0\leq j<d$ suffice.
Each is a direct summand of $F(k_A(j))$ by the splitting above.
Hence the image of $F$ split-generates $\cD$.
\end{proof}

\begin{proposition}\label{prop:weighted-restriction}\label{prop:graded-restriction}
Suppose a stability condition $(Z,\mathcal P)$ on $\cC$, supported
on $v:K_0(\cC)\to\Lambda$, satisfies
\begin{equation}
 \tau^{-a}[1]\mathcal P(\phi)
   \subseteq\mathcal P([\phi,\infty))
   \qquad(\phi\in\R).
 \label{eq:weighted-phase-bound}
\end{equation}
Then restriction of scalars induces a stability condition on $\cD$ given by 
\eqref{eq:induced-slicing}--\eqref{eq:induced-charge}.
Its support lattice is the image $\im(v\circ R_*)\subset\Lambda$.
A total mass--Hom bound is preserved.
If the stability condition $(Z,\mathcal P)$ is of Gepner type $(\tau,2/d)$ and $a\leq d/2$,
the induced stability condition is again of the same Gepner type.
\end{proposition}

\begin{proof}
Tensoring~\eqref{eq:coordinate-resolution} with an object gives the
functorial triangle
\begin{equation}
 \tau^{-a}E\xrightarrow{z}E\longrightarrow RF(E)
   \longrightarrow\tau^{-a}E[1].
 \label{eq:restriction-monad}
\end{equation}
Rotating this triangle expresses $RF(E)$ as an extension of
$\tau^{-a}E[1]$ by $E$.  For $E\in\mathcal P(\phi)$, both terms
belong to $\mathcal P([\phi,\infty))$ by
\eqref{eq:weighted-phase-bound}.  Lemma~\ref{lem:residue-field-base-change}
shows that $F(\cC)$ split-generates $\cD$, so
Lemma~\ref{lem:dg-induction} applies.
If $(Z,\mathcal P)$ is of Gepner type $(\tau,2/d)$, then
$\tau^{-a}E[1]$ has phase $\phi+1-2a/d$, which is at least $\phi$
when $a\leq d/2$.  Equality is allowed: if $a=d/2$, the two
terms lie in the same extension-closed slice $\mathcal P(\phi)$.
Since $R$ commutes with $\tau$, the formulas for the induced charge
and slicing show that the induced condition is also of Gepner type
$(\tau,2/d)$.

The containment~\eqref{eq:weighted-phase-bound} also holds for the
induced slicing $\mathcal P_{\cD}$, even without the Gepner assumption.
Indeed, if $E\in\mathcal P_{\cD}(\phi)$, then $R(E)\in\mathcal P(\phi)$,
and hence
\[
 R(\tau^{-a}E[1])\simeq\tau^{-a}R(E)[1]
       \in\mathcal P([\phi,\infty)).
\]
The characterization of phase intervals in
Lemma~\ref{lem:dg-induction} now gives
$\tau^{-a}E[1]\in\mathcal P_{\cD}([\phi,\infty))$.
\end{proof}

\begin{lemma}
\label{lem:waring-smooth-flag}
Let $f\in\operatorname{Sym}^d(V^*)$, $d\geq2$, define a smooth
hypersurface $X\subset\Pj(V)$.  For all sufficiently large $M$,
there is a linear embedding $\Pj(V)\hookrightarrow\Pj^{M-1}$
identifying $X$ with a linear section of the Fermat hypersurface
\[
 Y_M=\{y_1^d+\cdots+y_M^d=0\}\subset\Pj^{M-1}.
\]
Moreover, there is a flag of projective subspaces from the embedded $\Pj(V)$ down to
$\Pj^{M-1}$, with each inclusion of codimension one, such that
every section of $Y_M$ along the subspace in the flag is smooth.
\end{lemma}

\begin{proof}
A Waring decomposition realizes a hypersurface as a linear section
of a Fermat hypersurface; see~\cite[\S~1]{DimcaSticlaruWaring}.
We recall the construction and choose smooth intermediate sections.
Writing $N=\dim V$, we seek a flag
$X=Y_M\cap L_N\subset Y_M\cap L_{N+1}\subset\cdots
          \subset Y_M\cap L_M=Y_M$,
          $L_k\simeq\Pj^{k-1}$,
where each $L_k$ is a hyperplane in $L_{k+1}$ and
every intersection $Y_M\cap L_k$ is smooth.

Powers of linear forms span $\operatorname{Sym}^d(V^*)$: polarization expresses
a product of linear forms as
\[
 m_1\cdots m_d=\frac1{d!}
 \sum_{I\subseteq\{1,\ldots,d\}}(-1)^{d-|I|}
        \left(\sum_{i\in I}m_i\right)^d.
\]
The coefficients can be absorbed into $d$th roots, giving
$f=\sum_{i=1}^{M_0}\ell_i^d$.  The $\ell_i$ span $V^*$: otherwise a
nonzero vector $v$ would satisfy both $f(v)=0$ and
$
 df(v)=d\sum_i\ell_i(v)^{d-1}\ell_i=0$,
contradicting smoothness.  Thus $v\mapsto(\ell_i(v))$ is a linear
embedding.  For $M\geq M_0$, add zero linear forms
to the expression. The projectivization realizes $X$ as the required linear section.

Write $L=\Pj(V)$ for its image.  On $L$, the Fermat equation
restricts to $f$, whose differential is nonzero at every point of
$X$.  Thus every linear subspace containing $L$ cuts out a
hypersurface smooth along $X$.
If $M=N$, the flag is trivial.  Otherwise, start with
$L_M=\Pj^{M-1}$.  Suppose $Y_M\cap L_k$ is smooth
and $k>N+1$.  The hyperplanes of $L_k$ containing $L$ have base
locus exactly $X$ on this hypersurface.  By Bertini's theorem
\cite[III, Corollary~10.9 and Remark~10.9.2]{Hartshorne}, a general
such hyperplane $L_{k-1}$ gives a section smooth away from $X$;
the differential calculation gives smoothness along $X$.
Continuing down to $L_{N+1}$ and taking $L_N=L$ as the final
hyperplane gives the required flag, since its last section is $X$.
\end{proof}

\begin{proposition}
\label{prop:fermat-reduction}
Fix $d\geq2$.  Suppose, for arbitrarily large integers $M$, that
$\MF^{\mathrm{gr}}(\sum_{i=1}^M y_i^d)$ has a stability condition
satisfying the support property with respect
to a specified finite-rank lattice and satisfying
\eqref{eq:weighted-phase-bound} with $a=1$.
Then $\MF^{\mathrm{gr}}(f)$ has a stability condition
for every smooth degree-$d$ projective hypersurface $\{f=0\}$.
If the source conditions are of Gepner type $(\tau,2/d)$, so are the
resulting conditions.  Their support lattices are the successive
images specified in Lemma~\ref{lem:dg-induction}.
If the source conditions have total mass--Hom bounds, so do the
resulting conditions.
\end{proposition}

\begin{proof}
Choose $M$ large enough to apply
Lemma~\ref{lem:waring-smooth-flag}.  Along its smooth flag, each
hyperplane becomes a coordinate hyperplane after a linear change
of variables, and every defining equation has an isolated critical
point at the origin.  Apply
Proposition~\ref{prop:weighted-restriction} successively with $a=1$:
the phase bound is preserved at each step, so the construction
can be iterated.  The same proposition gives the asserted support
lattices, preserves mass--Hom bounds, and preserves Gepner type
$(\tau,2/d)$ since $1\leq d/2$.
\end{proof}

The following theorem implies Theorem~\ref{thm:main}.

\begin{theorem}
\label{thm:all-degrees}
Give each variable $x_i$ degree one, and let
$W\in\C[x_1,\ldots,x_N]$ be homogeneous of degree $d\geq2$,
with an isolated critical point at the origin.
The category $\MF^{\gr}(W)$ admits a numerical Bridgeland stability
condition of Gepner type $(\tau,2/d)$, with a total mass--Hom bound
and support on the full numerical lattice $K_{\num}(\MF^{\gr}(W))$.
Here $\tau$ is the internal grading shift.
For a smooth hypersurface $X=(W=0)\subset\Pj^{N-1}$ with $N>d$,
this gives a stability condition on $\Ku(X)$.  For $N=d$, it gives
a stability condition on $\Db(X)$.
\end{theorem}

\begin{proof}
If $N=1$, rescaling identifies $W$ with $x^d$, so the assertion is
already contained in Theorem~\ref{thm:all-degree-fermat}.
For $N\geq2$, the hypersurface $X=(W=0)\subset\Pj^{N-1}$ is smooth.
Theorem~\ref{thm:all-degree-fermat} supplies a stability condition
of Gepner type $(\tau,2/d)$ with the support property and a total mass--Hom bound
on every Fermat category.  For such a condition,
$
 \tau^{-1}[1]\mathcal P(\phi)
     =\mathcal P(\phi+1-2/d)
     \subseteq\mathcal P([\phi,\infty))$,
since $d\geq2$.  Proposition~\ref{prop:fermat-reduction} therefore
gives a stability condition of Gepner type $(\tau,2/d)$ on $\MF^{\gr}(W)$,
with support on the successive image lattice.
These class maps are numerical.  Indeed, for the adjoint functors
$F:\cC\to\cD$ and $R:\cD\to\cC$ at each restriction step,
$
 \chi_{\cC}(A,R(E))=\chi_{\cD}(F(A),E)$.
Thus $R_*$ sends numerically trivial classes to numerically trivial
classes, and the numerical class map on $\cC$ induces one on $\cD$.
The mass--Hom bound is preserved at each restriction step.
For $N\ge2$, Orlov's semiorthogonal decompositions express
$\MF^{\gr}(W)$ either as an admissible component of $\Db(X)$
or as a category obtained from $\Db(X)$ by adjoining finitely
many exceptional objects.  Its numerical group is therefore a
finite-rank lattice, and Lemma~\ref{lem:mass-hom-generator}
gives the full numerical support property.
Finally, Orlov's theorem~\cite[Theorem~3.11]{Orlov} identifies
$\MF^{\gr}(W)$ with $\Ku(X)$ if $N>d$, and with $\Db(X)$ if
$N=d$.  
\end{proof}

\subsection{Proper moduli spaces}

\begin{proof}[Proof of Corollary~\ref{cor:proper-gepner-moduli}]
By Orlov's equivalence, $\mathcal C$ is smooth, proper, and
idempotent complete.  Theorem~\ref{thm:all-degrees} gives the full
numerical support property and a total mass--Hom bound.
Since the class map is numerical, $\sigma$ is locally constant
in the sense of~\cite[Example~2.29]{HLR}.  Over $\C$, the
mass--Hom bound therefore allows us to apply the moduli theorem
\cite[Theorem~2.35 and Remark~2.36]{HLR}.

Choose $a$ with $\phi\in(a,a+1]$.  In the heart
$\mathcal P((a,a+1])$, semistable objects of class $v$ have
phase $\phi$.  By~\cite[Proposition~2.41]{HLR}, their stack is
algebraic of finite type with affine diagonal and admits a
proper good moduli space.
\end{proof}

Recall that a morphism $\pi:\mathfrak M\to M$ from a stack to an
algebraic space is a \emph{good moduli space} if $\pi_*$ is exact
on quasi-coherent sheaves and the natural map
$\OO_M\to\pi_*\OO_{\mathfrak M}$ is an isomorphism.
In our case, its closed points classify semistable objects up to
S-equivalence: two objects are S-equivalent when their
Jordan--H\"older graded objects are isomorphic.  Equivalently,
these points classify polystable objects up to isomorphism.
Fixing the phase is essential, since even shifts preserve both
the numerical class and the mass.  If all semistable objects
of the prescribed class and phase are stable, the resulting
space is their proper coarse moduli space.

\end{document}